\documentclass[11pt,a4paper,reqno]{amsart}

\usepackage[colorlinks,bookmarksopen, bookmarksnumbered,citecolor=red,urlcolor=red]{hyperref}
\usepackage{amsaddr}
\usepackage{epsfig,amsmath,amssymb,amsfonts,mathrsfs}
\usepackage{graphicx}
\usepackage{float} 
\usepackage{textcomp, eufrak}
\usepackage{verbatim}
\usepackage{fancyhdr}
\usepackage[numbers]{natbib}
\usepackage{enumitem}
\usepackage{setspace}
\usepackage{caption}
\usepackage{scalerel}[2016/12/29]

\usepackage{xcolor}
\newcommand{\PP}{\hspace*{.000cm}\mathord{\raisebox{-0.095em}{\scaleobj{.86}{\includegraphics[width=1em]{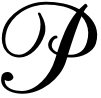}}}}\hspace*{.02cm}}

\numberwithin{equation}{section}

\newtheorem{theorem}{Theorem}[section]
\newtheorem{definition}[theorem]{Definition}
\newtheorem{lem}[theorem]{Lemma}
\newtheorem{prop}[theorem]{Proposition}
\newtheorem{rem}[theorem]{Remark}
\newtheorem{cor}[theorem]{Corollary}
\newtheorem{exa}[theorem]{Example}

\numberwithin{equation}{section}
\newtheorem{Assum}[theorem]{Assumption}
\newcommand{\rd}{\textrm{d}}

\newcommand{\frot}{{\textstyle \frac{1}{2}}}
\newcommand{\R}{\mathbb{R}}

\newcommand{\N}{\mathbb{N}}

\newcommand{\Rp}{\mathbb{R}^+}
\newcommand{\F}{\mathcal{F}}

\newcommand{\PS}{(\Omega, \mathcal{F}, \mathbb{P})}

\newcommand{\I}{\mathbb{I}}

\newcommand{\Bb}{\mathcal{B}}

\newcommand{\D}{\mathcal{D}}

\newcommand{\p}{\mathbb{P}}

\newcommand{\bn}{\begin{definition}\rm}
\newcommand{\en}{\end{definition}} 
\newcommand{\bt}{\begin{theorem}}                
\newcommand{\et}{\end{theorem}}
 \newcommand{\bnm}{\begin{enumerate}}              
\newcommand{\enm}{\end{enumerate}}
\newcommand{\br}{\begin{rem}} 
\newcommand{\er}{\end{rem}}

\newcommand{\om}{\omega}

\newcommand{\Om}{\Omega}

\newcommand{\btm}{\begin{itemize}}
 \newcommand{\etm }{\end{itemize}}

\newcommand{\E}{\mathbb{E}}

\newcommand{\Rd}{\R^d}

\newcommand{\ccdot}{\,\cdot\,}
\newcommand{\frat}{{\textstyle \frac{1}{2}}}
\newcommand{\mmod}{\hspace{-0.2cm}\mod}

\newcommand{\Ic}{\mathcal{I}}
\newcommand{\Jc}{\mathcal{J}}
\newcommand{\Ac}{\mathcal{A}}

\newcommand{\Sc}{[0, \tau)}

\newcommand{\Rm}{\R^m}

\newcommand{\LG}{\mathcal{L}}

\newcommand{\Pt}{\mathcal{P}}

\newcommand{\XX}{\mathcal{X}}

\newcommand{\Mod}{\;\mathrm{mod}\;}
 \newcommand{\BT}{\Sc\times\Rd}
 
 \binoppenalty=\maxdimen
\relpenalty=\maxdimen

\renewcommand\leftmark{{Time-periodic measures, random periodic orbits, and the linear response for dissipative non-autonomous SDE's}}

\renewcommand\rightmark{{Time-periodic measures, random periodic orbits, and the linear response for dissipative non-autonomous SDE's}}
\begin{document}
\author{Micha\l \;Branicki and Kenneth Uda}

\address{\vspace*{-0.2cm} Department of Mathematics, University of Edinburgh, Scotland, UK}
\email{M.Branicki@ed.ac.uk, K.Uda@ed.ac.uk}

\vspace*{-.3cm}\title{Time-periodic measures, random periodic orbits, and  the linear response for dissipative non-autonomous stochastic differential equations}
\begin{abstract}
We consider a class of dissipative stochastic differential equations (SDE's) with time-periodic coefficients in finite dimension, and the response of time-asymptotic probability measures induced by such SDE's to sufficiently regular, small perturbations of the underlying dynamics. Understanding such a response provides a systematic way to  study changes of  statistical observables in response to perturbations, and  it is often very useful for sensitivity analysis, uncertainty quantification, and for improving probabilistic predictions of nonlinear dynamical systems, especially in high dimensions. Here, we are concerned with the linear response to small perturbations in the case when the time-asymptotic probability measures are time-periodic. First, we establish sufficient conditions for the existence of  stable random time-periodic orbits generated by the underlying SDE. Ergodicity of time-periodic probability measures supported on these random periodic orbits is subsequently discussed. Then, we derive the so-called fluctuation-dissipation relations which allow to describe the linear response of statistical observables to small perturbations away from the time-periodic ergodic regime in a manner which only exploits the unperturbed dynamics. The results are formulated in an abstract setting but they apply to problems ranging from aspects of climate modelling, to molecular dynamics, to the study of approximation capacity of neural~networks and robustness of their estimates.

\end{abstract}

\date{}
\maketitle

\tableofcontents

\newpage
\section{Introduction}\label{intro}

In many scientific applications a systematic determination of the  response of a complex nonlinear dynamical system to  time-dependent perturbations is of key importance;  topical examples in high-dimensional, non-autonomous and/or stochastic settings include climate models (e.g., \cite{ Abramov07, Majda05, Madja10, Madja10b,grit99,Bran12}), statistical physics and non-equilibrium thermodynamics (e.g., \cite{kubo57, ruel98, ruel99, Johnson, Nyquist}),  and even neural networks (e.g., \cite{cessac19,park18,Coolen01}). The sought response is usually quantified in terms of a change in an `observable' expressed as a statistical/ensemble average of some  functional defined on the trajectories of the underlying dynamical system.  The classical theory of linear response (e.g., \cite{ruel98b,Majda05}) is concerned with capturing changes in observables to sufficiently small perturbations of the original dynamics close to its statistical equilibrium. It turns out that in such a setting the response can be expressed, with some caveats,  through formulas linking the external perturbations to spontaneous fluctuations and dissipation in the unperturbed time-asymptotic dynamics (e.g., \cite{ruel98, ruel99,lucar08}). The classical fluctuation-dissipation theorem (FDT) is of fundamental importance in  statistical physics (e.g., \cite{kubo66,Agarwal,Einstein}), and it roughly states that for systems of identical particles in statistical equilibrium, the average response to small external perturbations can be calculated through the knowledge of suitable correlation functions of the unperturbed time-asymptotic dynamics; see, for example, \cite{kubo85,bal97} for some of the many applications of the FDT in the statistical physics setting.

The validity of the linear response and fluctuation-dissipation relationships for more general dynamical systems encountered, for example, in climate modelling (e.g.,~\cite{Majda05})  is an important topic which is particularly relevant  for uncertainty quantification in reduced-order predictions and reduced model tuning (e.g.,~\cite{Madja10b,MAG10, Bran12,Madja12}).  In an early influential work  Leith \cite{Leith75} suggested that if the climate dynamics satisfied a suitable FDT,  the climate response to small external forcing could  be calculated by estimating suitable statistics in the unperturbed  climate\footnote{\,The meaning of the term `climate' used in most theoretical work in atmosphere-ocean science is loosely related to properties of the probability measure induced by the time-asymptotic dynamics.}.  Climate dynamics  is modelled as a forced dissipative chaotic or stochastic dynamical system which is arguably rather far from the statistical physics' setting for FDT. Nevertheless, Leith's conjecture stimulated a lot of activity in generating new theoretical formulations (e.g., \cite{Hairer10,Madja10}) and in designing approximate algorithms for FDT to study the  climate response (e.g.,~\cite{ Abramov07, Abramov08,Abramov09, Majda05, Madja10, Madja10b,grit99,grit02,grit07,grit08,majdaqi19}). However, despite numerous applications in autonomous and non-autonomous settings,  there is little rigorous evidence supporting the validity of the linear response and FDT  in the non-autonomous setting beyond the formal derivation of FDT for time-dependent stochastic systems \cite{Madja10}.

The goal here is to provide a more rigorous justification of the linear response theory  for a class of forced dissipative stochastic differential equations (SDE's) with time-periodic coefficients which induce  time-periodic probability measures.  Our objective is twofold: 
\begin{itemize}[leftmargin=0.8cm]  
\item[(i)] Establish sufficient conditions for the existence and ergodicity  (in an appropriate sense) of time-periodic measures associated with time-asymptotic dynamics  for a  class of {\it `dissipative'} SDE's (defined later in (\ref{diss_linear_growth})) with time-periodic coefficients in finite dimensions. 

\vspace{.2cm}\item[(ii)] Analyse the linear response of such SDE's in the time-periodic regime to small perturbations, and express the change in the statistical observables based  on time-periodic ergodic measures via  fluctuation-dissipation type relations. 
\end{itemize}

The results derived in the sequel will concern SDE's whose time-periodic measures are supported on certain stable  random periodic solutions.  In principle, the results discussed in the context of the linear response apply to a wider class of dynamical systems generating time-periodic measures; however, establishing conditions for the existence and ergodicity of such measures in a more general setting (for SDE's or otherwise) is not trivial and is beyond the scope of this work.

\smallskip
Time-periodic probability measures associated with the time-asymptotic dynamics are arguably ubiquitous in many mathematical models. In particular, seasonal and diurnal cycles in climate models due to time-periodic forcing or retarded self-interactions in neural networks provide  some of the obvious candidates, and highlight the need for developing the linear response framework in the time-periodic setting. It is worth stressing that the need for rigorous formulation of the linear response and FDT for dissipative stochastic dynamical systems (in line with, e.g., \cite{ris89,Majda05, majda08, Pal01, Pal10, majdaqi19}) is justified by contemporary approaches to  the simulation and  reduced-order modelling of high-dimensional, multi-scale dynamical phenomena. For example,  comprehensive models for climate change prediction  or molecular dynamics simulations involve stochastic components (e.g., \cite{Pal01,MTV03,Pal10,majda08, Zwanzig73, Abramov10,Erban19}) to mimic the effects of unresolved dynamics, while reduced-order models typically involve stochastic noise terms (e.g., \cite{Chorin00,MTV01,DKelly17,Mourat16,Bella16}. Here, similar to \cite{Hairer10,Madja10}, the presence of noise leads to improved regularity of the problem which simplifies key aspects of the analysis  compared to deterministic, dissipative  nonlinear  systems (e.g., \cite{ruel97,bal08,bal08b,bal93,Gouz08}). As a consequence, we are able to focus on systems that have other important features of realistic dynamics, namely a lack of ellipticity, non-compactness of state space, and a lack of global Lipschitz continuity of the coefficients in the underlying SDE. The results established below apply to a broad class of nonlinear functionals which include common quantities of interest, such as the mean and covariance of subsets of~variables.

 \section{General setup}\label{stp_sec}
 \newcommand{\Xt}{\mathcal{M}}
 Our framework relies on the theory of Markovian\footnote{\,Here, the notion of a {\it `Markovian RDS'} means that there exists a version of the RDS which has the Markov property w.r.t.~the canonical filtration generated on the Wiener space by the Brownian motion discussed later.} random dynamical systems (RDS), which provides a  geometric link between stochastic analysis and dynamical systems.  This relationship  was established through the discovery (e.g.,~\cite{Kunita,Arnold}) that for sufficiently regular coefficients $b,\sigma$ the stochastic differential equation (SDE)
\begin{align}\label{NSDE}
dX_t = b(t,X_t)dt+\sigma(t,X_t)dW_{t-s}, \qquad X_{s} \in \Rd, 
\end{align} 
 generates a stochastic flow  $\{\phi(t,s,\ccdot, \cdot\,)\!:\, s,t\in \Ic\subseteq\R, \,s\leqslant t \}$ of homeomorphisms on $\Rd$ such~that 
  $$X_t^{s,x}(\om) = \phi(t,s,\om,x)\qquad \p \,\textrm{-\,a.s.}$$ 
  for $x=X_s(\omega)$,  $\om\in \Omega$ in the Wiener space $\PS$ with $W_t$ an  $m$-dimensional Brownian motion. For $b = b(x)$, $\sigma=\sigma(x)$, the SDE will be called {\it autonomous}, and {\it non-autonomous} otherwise. It turns out (e.g., \cite{Arnold,Arnoldp}) that, for an autonomous SDE (in the above sense) with sufficiently regular coefficients  there exists essentially a one-to-one correspondence between the SDE and an RDS; a rough but convenient interpretation (skipping the filtration) is that in the autonomous case there exists an RDS generating the SDE, which in turn generates the stochastic flow and vice-versa.

  One of the key concepts relevant for the analysis of the long-time behaviour of RDS is the extension of the notion of ergodicity to  the random setting (e.g.~ \cite{Arnold, Bax2, Bax3, Hans3, Has12, Mattingly, Meyn93, Meyn931}). These important results are established in the regime of (random) stationary measures and (random) stationary processes,  in the case when the source of time-dependence is only due to  the noise process (i.e., $b=b(x)$, $\sigma = \sigma(x)$ in (\ref{NSDE}) and the SDE is autonomous in the jargon established above). Over the last decade significant progress has been made in the  study of the long-time behaviour of SDE's generated by time-dependent vector fields (e.g., \cite{Feng11, Feng12, Feng18, Uda16, Uda18, Wan14, Zhao09, Cherub17}). Based on the insight from the latter results, we shall study the ergodicity of SDE's with time-periodic coefficients in order to  establish fluctuation-dissipation formulas through the linear response in the random periodic regime.  Our strategy is to first prove the existence of a unique time-periodic measure under certain  {\it `dissipative'} assumptions on the SDE via a version of Lyapunov second method and coupling.
  The standard Lyapunov second method is a well-known and powerful technique for the investigation of stability of solutions of nonlinear dynamical systems in finite and infinite dimensions. An extension of this method to an RDS generated by an autonomous SDE is essentially due to Has\'minskii (e.g.,~\cite{Has12}); subsequent extensions  include applications to SDE's with random switching (e.g.,~\cite{Mao94}) and to the case  of non-trivial random stationary solutions and random attractors by Schmalfuss \cite{Schmal01}.  Importantly, this method involves the study of random invariant sets (under the considered dynamics) without the need for the explicit knowledge of solutions of the underlying SDE, and it is based solely on the vector fields encoded in the coefficients  of the SDE  even when the drift term, i.e., $b$ in (\ref{NSDE}), is only locally Lipschitz continuous. However, in the present non-autonomous, time-periodic setup, the lack of stationarity and the unavoidable skew-product structure of the underlying dynamics pose additional challenges when dealing with ergodicity of time-asymptotic probability measures.  The main issue which prevents one from using the `classical' methods (e.g., \cite{Doob48,Has12}) for proving ergodicity the random periodic regime stems from the fact that these probability measures are defined on the skew-product fibre bundle in the space of measures on the time-extended state space and that they are not mixing. 
 Here, this complication is overcome  by employing an extension of the  Krylov-Bogolyubov procedure (e.g., \cite[\S 1.5]{Arnold}) which allows for dealing with ergodicity of probability measures on appropriate Poincar\'e sections in the narrow topology generated by the dual of an appropriate  discrete-time transition semigroup, and then `linking' the results via the continuous-time transition  semigroup induced by the SDE dynamics on the space of skew-product probability measures.  In the present case, the properties of the time-periodic measures established with the help of the Lyapunov's second method for dissipative SDE's allows us to dispense with explicit assumptions on the ergodicity in the Poincar\'e sections which is otherwise required.

  The rest of the article is organised as follows. In the remainder of this section we outline the frequently used notation. In Section \ref{Random pr}, we recap some basic results and definitions, including  the notion of a Random Dynamical System (RDS) generated by an SDE in finite dimensions, and we outline the notion of a random periodic process.  In Section \ref{Erg_S}, we first prove the existence of stable random time-periodic solutions for a class of dissipative SDE's with time-periodic coefficients,  and the existence of the associated time-periodic measures (\S\ref{permeas}); sufficient conditions for ergodicity of such measures (in an appropriate sense) are established in~\S\ref{permes_erg}. Section \ref{Linear response_per} deals with the linear response theory in the above setting. The derivation of  the linear response formula in the time-periodic setting is followed by the derivation of  two classes of  fluctuation-dissipation relationships:  the first one applies to perturbations of dynamics with time-periodic ergodic probability measures required only to exist in the unperturbed dynamics; the second fluctuation dissipation formula involves simpler formulas but it requires persistence of time-periodicity in the perturbed~measures.

\subsection{Function spaces}\label{gennot} 
\addtocontents{toc}{\protect\setcounter{tocdepth}{1}}
Below, we outline function spaces which are used in the sequel. 

Let $(\XX, \text{d})$ be a complete separable metric space. We consider either $\XX= \Rd$ or $\XX= \R\times \Rd$, or a flat cylinder $\XX = \Sc\times \R^d$, $0<\tau<\infty$, $\Sc\simeq \R\Mod\tau$, which arises when `lifting' the dynamics generated by a non-autonomous SDE. In this section, we use $\XX$ for all these spaces  to unify the notation. Throughout, we set $\N_0:=\{0,1,2, \dots\}$ and $\N_1:=\{1,2, \dots\}$, while $\PP(\Sc)$, $\PP(\Rd)$  are the spaces of Borel probability measures on, respectievely, $\Sc$ and $\Rd$.

\vspace{.1cm}
\noindent \raisebox{.05cm}{\tiny$\bullet$}  $\PS$ is  the Wiener probability space where $\Om:=\mathcal{C}_{0}(\R;\R^m)$, $m\in \N_1$; i.e., the sample space $\Om$ is identified with a linear subspace of continuous functions $\mathcal{C}(\R;\R^m)$ which vanish at zero. $\F$ is the  Borel $\mathfrak{S}$-algebra on $\Om$ generated by open subsets  in the compact-open topology defined~via
\begin{align*}
\varrho(\om,\hat{\om})= \sum_{n=0}^\infty\frac{1}{2^n}\frac{\Vert \om -\hat{\om}\Vert_n}{1+\Vert \om -\hat{\om}\Vert_n}, \qquad \Vert \om -\hat{\om}\Vert_{n} := \sup_{t\in [-n, n]}\vert \om(t) -\hat{\om}(t)\vert, \quad \om,\hat\om\in \Omega, 
\end{align*}
with $|{\cdot}|$ the Euclidean norm. Finally, $\p$ is the Wiener measure on $\mathcal{F}$.  In such a setup the canonical Wiener process (with two-sided time) on $(\Om,\mathcal{F})$ with values in the Borel-measurable space $(\R^m,\mathcal{B}(\R^m))$ is  defined as $W_t(\om) = \om(t)$, $t\in \R$, via the identification of $\om\in \Om$ with  functions $\om(\ccdot)\in \mathcal{C}_{0}(\R;\R^m)$; see, e.g., \cite[Appendix A.2]{Arnold} and references therein for details.

\vspace{.1cm}
\noindent \raisebox{.05cm}{\tiny$\bullet$} Given the probability space $\PS$ and $\mathcal{G}\subseteq\F$,  $L^p(\Om, \mathcal{G}, \p)$, $1\leqslant p <\infty$, is the space of $\mathcal{G}$-measurable random variables $X:\Om\rightarrow\Rd$ such that $\E\vert X\vert^p:=\int_\Omega\vert X(\om)\vert^p\,\p(d\om)<\infty, $ and equipped with the norm 
$\Vert X\Vert_p : = \left(\E\vert X\vert^p\right)^{1/p}.$

\vspace{.1cm}
\noindent \raisebox{.05cm}{\tiny$\bullet$} Given the (Borel) measurable space $\big(\XX,\mathcal{B}(\XX)\big)$,  where $\mathcal{B}(\XX)$ denotes the Borel $\mathfrak{S}$-algebra over~$\XX$. 
(For $\XX = \XX_1\times \XX_2$  we consider the product algebras $\mathcal{B}(\XX_1)\times\mathcal{B}(\XX_2)$.)  
\\ 
\hspace*{.3cm}-\,$\mathbb{M}(\XX)$ denotes the space of measurable functions on $\XX$, \\
\hspace*{.3cm}-\,$\mathcal{C}(\XX)$ denotes the space of continuous functions on $\XX$, \\
\hspace*{.3cm}-\,$\mathbb{M}_\infty(\XX)$ denotes  the space of bounded, measurable, real-valued, scalar functions on $\XX$, i.e.,  
 \begin{align*} 
 \mathbb{M}_\infty(\XX) :=\big\{f: \XX\rightarrow \R, \;\,f\in \mathbb{M}(\XX): \;\Vert f\Vert_{\infty}<\infty\big\}, \quad \|f\|_\infty := \sup_{x\in \XX} |f(x)|.
\end{align*}
\hspace*{.3cm}-\,$ \mathcal{C}_\infty(\XX)$  denotes  the space of  bounded, real-valued, scalar, continuous functions on $\XX$.

\noindent \raisebox{.05cm}{\tiny$\bullet$} The space $\mathcal{C}^l(\XX)$ contains  $l$-times continuously differentiable real-valued functions.   

\noindent \raisebox{.05cm}{\tiny$\bullet$} The space $\mathcal{C}_\infty^l(\XX)$ contains those functions in $\mathcal{C}^l(\XX)$ which are bounded.  

\vspace{.1cm}
\noindent \raisebox{.05cm}{\tiny$\bullet$} The space $\mathcal{C}^{l,k}(\XX_1\times\XX_2)$ denotes the space of functions which are $\mathcal{C}^l$ on $\XX_1$, and $\mathcal{C}^k$ on $\XX_2$. The space $\mathcal{C}^{l,k}_\infty$ contains bounded functions in $\mathcal{C}^{l,k}$.

\noindent \raisebox{.05cm}{\tiny$\bullet$}   The space of bounded, real, scalar  Lipschitz continuous functions on $(\XX,\rd)$ is denoted by
 \begin{align*}
  \text{Lip}_\infty(\XX)&:=\big\{f\in\mathcal{C}_\infty(\XX): \;\Vert f\Vert_{\textsc{bl}}<\infty\big\},\\
  \Vert f\Vert_\textsc{bl} &:= \max\big\{\Vert f\Vert_{\infty},\text{Lip}(f)\big\},\quad \text{and} \quad \text{Lip}(f) := \sup\bigg\{\frac{\vert f(y)-f(z)\vert}{\text{d}(y,z)}: y\neq z, \;y,z\in \XX\bigg\}.
 \end{align*}

\noindent \raisebox{.05cm}{\tiny$\bullet$}  $\tilde{\mathcal{C}}^{l,\delta}(\XX)$, $l\in \N_0$, $0<  \delta \leqslant 1$,  is the Fr\'echet space of functions $f\,{:}\; \XX\rightarrow \XX$, whose  $l$-th derivatives  are $\delta$-H\"older continuous, and which is furnished with the countable family of semi-norms 
\begin{align*}
&\Vert f\Vert_{l,0,N}:=  \sup_{x\in \XX}\frac{|\langle  f(x),x\rangle|}{1+\vert x\vert^2}+\sum_{1\leqslant  \vert \beta\vert \leqslant  l}\,\sup_{x\in B_N} \vert D^{\beta}f(x)\vert, \\[.0cm]
& \Vert f\Vert_{l,\delta;N}: = \Vert f\Vert _{l,0;N}+\sum_{\vert \beta\vert =l}\,\sup_{x,y\in B_N, x\neq y}\frac{\vert D^{\beta}f(x)-D^{\beta}f(y)\vert}{\vert x-y\vert^{\delta}}, 
\end{align*} 
where $|{\cdot}|$ is the Euclidean norm and $\langle \cdot,\cdot\rangle$ the dot product on $\XX$, and $B_N\,{=}\,\{x \,{\in}\,\XX \!: \vert x\vert\,{\leqslant}\, N\}$, $N\in \N_1$, is a closed ball in $\XX$  with radius $N$, and 
\begin{align*}
D^{\beta}f(x): = \frac{\partial^{\vert \beta\vert}f}{(\partial x_1)^{\beta_1}\cdots(\partial x_d)^{\beta_n}}, \quad \vert \beta\vert := \beta_1+\cdots+\beta_n\;, \quad \beta_i\in \mathbb{N}_0, \;\;\;i=1,\dots,n,\quad D^0 = \I,
\end{align*} 
denotes the Fr\'echet derivative; $f\in \tilde{\mathcal{C}}^{l,\delta}(\XX)$ if  $\Vert f\Vert_{l,\delta;N}<\infty$ for $N\in \mathbb{N}_1$, $N<\infty$.

\smallskip
\noindent \raisebox{.05cm}{\tiny$\bullet$}  $\bar{\bar{\mathcal{C}}}^{l,\delta}(\XX)$, $l\in \N_0$, $0< \delta \leqslant 1$,  is the space of functions $f\,{:}\; \XX\rightarrow \XX$, whose  $l$-th derivatives  are $\delta$-H\"older continuous with the norm 
\begin{align*}
&\Vert f\Vert_{l,0}:=  \sup_{x\in \XX}\frac{|f(x)|}{1+|x|}+\sum_{1\leqslant  \vert \beta\vert \leqslant  l}\,\sup_{x\in \XX} \vert D^{\beta}f(x)\vert, \\[.1cm]
& \Vert f\Vert_{l,\delta}: = \Vert f\Vert _{l,0}+\sum_{\vert \beta\vert =l}\,\sup_{x,y\in \XX, x\neq y}\frac{\vert D^{\beta}f(x)-D^{\beta}f(y)\vert}{\vert x-y\vert^{\delta}}. 
\end{align*} 
Functions $f\in \bar{\bar{\mathcal{C}}}^{l,\delta}(\XX)$ are such that  $\Vert f\Vert_{l,\delta}<\infty$.


\section{Random periodic processes} \label{Random pr}
In order to facilitate subsequent derivations, we recall  definitions of random dynamical systems (RDS) generated by SDE's (see,~e.g.,~\cite{Hkunita, Kunita, Arnold,Arnoldp}), random periodic processes (see,~e.g.,~\cite{Feng11, Feng12, Feng18, Uda16, Uda18, Wan14, Zhao09, Cherub17} and transition evolutions generated by SDE's (see,~e.g.,~\cite{Arnold, Da Prato, Da Prato2, Kunita}).  We also  provide an intuitive example of a random periodic solution arising in the stochastic dynamics of periodically forced FitzHugh--Nagumo model. All definitions below are restricted to $\Rd$ but a number of them are subsequently extended (explicitly or otherwise) to the skew-product representation in $\Sc\times\Rd$ which is used to deal with the non-autonomous dynamics. 
\begin{definition}[{\bf Stochastic flow} \cite{Hkunita, Kunita}]\label{Flow} \rm 
 Let $\phi(t,s,\om,x)\in \Rd$, $s,t\in \Ic\subseteq\R$, $ x \in \Rd$,  be a random field on a  probability space $(\Om, \F, \p)$. The two-parameter family $\{\phi(t,s,\ccdot, \ccdot)\!:\, s,t\in \Ic\subseteq\R \}$ is called a {\it stochastic flow of homeomorphisms} if there exists a null set $\mathcal{N}\subset \Om$ such that for any $\om\,{\notin}\, \mathcal{N},$ there exists a family of continuous maps $\{\phi(t,s,\om,\ccdot)\!:\, s,t\in\Ic\}$ on $\Rd$ satisfying
\begin{itemize}[leftmargin=0.9cm]
\item[(i)] $\phi(t,s,\om,\ccdot) = \phi(t,u,\om,\phi(u,s,\om,\ccdot))$ holds for any $s,t,u\in \Ic,$
\item[(ii)] $\phi(s,s,\om,\ccdot) = \textup{id}_{\XX}$, for all $s\in\Ic,$
\item[(iii)] the map $\phi(t,s,\om,\ccdot): \Rd\rightarrow \Rd$ is a homeomorphism for any $t,s\in \Ic$.
\end{itemize}
The map $x\mapsto \phi(t,s,\om,x)$ is a {\it stochastic flow of $\mathcal{C}^l$-diffeormorphisms}, if it is a homeomorphism and $\phi(t,s,\om,x)$ is $l$-times continuously differentiable with respect to $x\in \Rd$ for all $s,t\in \Ic\subseteq\R$ and the derivatives are continuous in $(s,t,x)\in\Ic\times \Ic\times\Rd.$ The stochastic flow is referred to as {\it `forward'} for $s\leqslant t$, and as {\it `backward'} for $t\leqslant s$. In the sequel, we will confine the discussion to  $(\Om, \F, \p)$ which is  the Wiener space defined in \S\ref{gennot}.
\end{definition}

\begin{definition}[{\bf Filtration generated by a stochastic flow}] \rm Given a probability space $(\Om, \F, \p)$, let $\F_s^{\hspace{.02cm}t}\subseteq \mathcal{F}$ be the smallest $\mathfrak{S}$-algebra on $\Om$ generated by  $$\cap_{\varepsilon>0}\,\mathfrak{S}\big(\phi(u,v,\ccdot,\ccdot)\!:\; s-\varepsilon\leqslant u, v\leqslant t+\varepsilon\big), $$ and containing all null sets of $\F$. The {\it two-parameter filtration} $\{\F_s^{\hspace{.02cm}t}: s\leqslant t\}$  is the filtration generated by the forward stochastic flow $\big\{\phi(t,s,\ccdot,\ccdot)\!: \;s,t\in \Ic\subseteq\R, \,s\leqslant t\big\}$ and  the  filtered probability space is denoted by $\big(\Om, \F, (\F_s^{\hspace{.02cm}t})_{s\leqslant t},\p\big)$.
\end{definition}

\begin{definition}[{\bf Transition kernel}]\label{Mkernel}\rm
Consider the stochastic flow $\big\{\phi(t,s,\ccdot,\ccdot)\!: \,s,t\in \Ic;\; s\leqslant t\big\}$ induced by the SDE \eqref{NSDE} for some fixed $s\in \Ic$.  Given the Borel-measurable space $\big(\Rd,\mathcal{B}(\Rd\hspace{.03cm})\big)$, the {\it transition probability kernel} $P(s,x; t, \ccdot)$ induced by solutions of \eqref{NSDE} is defined~by
 \begin{equation}\label{phiP}
 P(s,x; t, A) = \p\big(\{\om\in\Om: \phi(t,s,\om,x )\in A\}\big) \quad \forall\;s,t\in\Ic,\;s\leqslant t,\; A\in\Bb\big(\Rd\hspace{.03cm}\big).
 \end{equation}
 The transition kernel satisfies the Chapman-Kolmogorov equation  
 \begin{equation}\label{ChK}
P(s,x;t,A) = \int_{\Rd} P(u,y;t,A)P(s,x;u,dy), 
\end{equation}
for any $s,t,u\,{\in}\, \Ic$, $s\leqslant u\leqslant t$, and for all $x\,{\in}\, \Rd$, $A\in \Bb(\Rd)$.
 \end{definition}
 \begin{definition}[{\bf Transition evolution and its dual}]\label{trans_ev}\rm
 Given the forward stochastic flow $\big\{\phi(t,s,\ccdot,\ccdot)\!: \,s,t\in \Ic;\; s\leqslant t\big\}$ and the transition kernel (\ref{phiP}) induced by the  solutions of~\eqref{NSDE}, the operator  $\mathcal{P}_{s,t}: \mathbb{M}_\infty(\Rd)  \rightarrow\mathbb{M}_\infty(\Rd)$ called the {\it transition evolution} is defined by 
 \begin{equation}\label{calP}
 \mathcal{P}_{s,t}\varphi(x) = \int_{\Rd}\varphi(y)P(s,x; t, dy) = \E\big[\varphi(\phi(t,s,x))\big] \quad \quad \forall\,s,t\in\Ic,\;s\leqslant t,\;x\in \Rd,
 \end{equation}
 where we use the shorthand notation $\E\big[\varphi(\phi(t,s,x))\big]:=\int_\Om \varphi\big(\phi(t,s,\om,x)\big)\p(d\om)$. The action of transition evolutions to  arbitrary measurable functions 
 is extended in a standard way. 
 
\smallskip
 For any probability measure $\mu_s \in \PP(\Rd)$, $s\in \Ic$, on $\big(\Rd,\Bb(\Rd)\big)$, the $L^2(\Rd;\mu_s)$ dual $\mathcal{P}_{s,t}^{*}$ of the transition evolution $\mathcal{P}_{s,t}$ is defined by
 \begin{equation}\label{P*mu}
 (\mathcal{P}_{s,t}^{*}\mu_s)(A) = \int_{\Rd}P(s,x; t, A)\mu_s(dx)  \quad  \quad \forall\;s,t\in\Ic,\;s\leqslant t,\;A\in \Bb\big(\Rd\big).
 \end{equation}
 Consequently, with the help of (\ref{ChK}), we have for any $s,u,t\,{\in}\, \Ic$, $s\leqslant u\leqslant t$, and for all $A\in \Bb(\Rd)$
 \begin{equation}
 \mu_t(A) = (\mathcal{P}_{s,t}^{*}\mu_s)(A) = (\mathcal{P}_{u,t}^{*}\mathcal{P}_{s,u}^{*}\hspace{.03cm}\mu_s)(A) = (\mathcal{P}_{u,t}^{*}\hspace{.03cm}\mu_u)(A).
 \end{equation}
 \end{definition}

\begin{theorem}[\textit{\textbf{Stochastic flows generated by solutions of SDE's}}]\label{SDEflow}
Suppose that the coefficients  of the SDE \eqref{NSDE} are such that $b(\,\cdot\,, x), \,\sigma(\,\cdot\,,x)$ are continuous for all $x\in \Rd$, and for all $t\in \R$,  $b(t,\ccdot),\,\sigma_k(t,\ccdot)\in \tilde{\mathcal{C}}^{l,\delta}(\Rd)$, $l\in \mathbb{N}_0$, $0<\delta\leqslant 1$, where $\{\sigma_k\}_{k=1}^m$, denote the columns of~$\sigma$. If the initial condition $X_s$ in (\ref{NSDE}) is independent of the $\mathfrak{S}$-algebra  generated by $W_{t-s}(\ccdot)$, $t\geqslant s$,  and $\E\big[|X_s|^2\big]<\infty$, there exist unique global solutions of \eqref{NSDE} which generate a forward stochastic flow of homeomorphisms $(l=0)$ or ${\mathcal{C}}^{l}$-diffeomorphisms $(l\geqslant 1)$ on $\Rd$
, $\big\{\phi(t,s,\ccdot,\ccdot)\!: \, s\leqslant t\big\}$ 
such that 
\begin{equation}
 X^{s,x}_{t}(\om)=\phi(t,s,\om,x) \qquad \p \,\textrm{-\,a.s.} \quad  \forall\,s,t\in \R, \;s\leqslant t, \;\; x\in \Rd ,
\end{equation}
and which are adapted to the filtration $(\mathcal{F}^{\hspace{.02cm}t}_s)_{s\leqslant t}$ on $(\Om,\mathcal{F},\p)$, see,  e.g., \cite[Thm 3.4.6 and \S 4.7]{Kunita} with slight modifications. If,  in addition,  $\E|X_s|^p<\infty$, $2\leqslant p<\infty$, then $\E\big[|X_t|^p\big]<\infty$, $s\leqslant t<\infty$. Stronger (e.g., dissipative)  growth conditions may have to be imposed on the coefficients $(b,\,\sigma)$ in the SDE \eqref{NSDE} to guarantee the existence of the absolute moments of the solutions for all time (see, e.g.,  Remark \ref{rem_hasm_lyap} and Appendix~\ref{app_hasminski}).
\end{theorem}

\begin{definition}{\bf \!(Infinitesimal generators).}\label{inf_gnrt} \rm
Let $f\in \mathcal{C}^{1,2}(\R\times\Rd,\R)$, and $t\mapsto\phi(t,s,\om,x)$, $s \leqslant t$, $\om\in \Om$, be a solution of the SDE in \eqref{NSDE}. Considering the evolution of $f\big(t,\phi(t,s,\om,x)\big)$  allows one to represent the infinitesimal generator  of solutions of \eqref{NSDE} through the second-order operator\footnote{\,Strictly, $\mathcal{L}_t$ in (\ref{Lf}) coincides with the generator of \eqref{NSDE} on $f\in \mathcal{C}^{1,2}_c(\R\times\Rd,\Rp)$ but it  is well-defined for $f\in \mathcal{C}^{1,2}(\R\times\Rd,\Rp)$ and we refer to $\mathcal{L}_t$ as the generator throughout; the same applies to $\mathcal{L}^{(2)}_t$ in (\ref{2pp}). } \textup{(e.g., \cite{Kunita})};
\begin{equation}\label{Lf}
\mathcal{L}_tf(t,x) =\partial_t f(t,x)+ \sum_{i=1}^{d}b_i(t, x)\partial_{x_i} f(t,x)+ \frat\sum_{i,j=1}^{d}\sum_{k=1}^m\sigma_{ik}(t,x)\sigma_{jk}(t,x)\partial_{x_ix_j}^2f(t,x),
\end{equation}
where $(b, \sigma)$ are sufficiently regular  drift and diffusion coefficients in the SDE \eqref{NSDE}. Analogously, for $g\in \mathcal{C}^{1,2}(\R\times\Rd\times\Rd,\R)$, the  infinitesimal generator of the {\it two-point motion} \textup{\cite[\S 4.2]{Kunita}}, $t\mapsto\big(\phi(t,s,\om,x), \phi(t,s,\om,y)\big)$,  of the flow $\big\{\phi(t,s,\ccdot, \ccdot)\!:\; \,s\leqslant t \big\}$ can be represented through the second-order differential operator 
\begin{align}
\notag \mathcal{L}^{(2)}_tg(t,x,y) &= \partial_t g(t,x,y)+\sum_{i=1}^d\Big( b(t,x)\partial_{x_i}g(t,x,y)+b(t,y)\partial_{y_i}g(t,x,y)\Big)\\
\notag &\hspace{.5cm} + \frac{1}{2}\sum_{i,j=1}^d\sum_{k=1}^{m}\Big(\sigma_{ik}(t,x)\sigma_{jk}(t,x)\partial^2_{x_ix_j}g(t,x,y)+\sigma_{ik}(t,x)\sigma_{jk}(t,y)\partial^2_{x_i y_j}g(t,x,y)\\
&\hspace{1.4cm} + \sigma_{ik}(t,y)\sigma_{jk}(t,x)\partial^2_{y_ix_j}g(t,x,y)+\sigma_{ik}(t,y)\sigma_{jk}(t,y)\partial^2_{y_iy_j}g(t,x,y)\Big). \label{2pp}
\end{align}
\end{definition}

\begin{definition} \!\!{\bf (Random Dynamical System} {\rm \cite{Arnold, Arnoldp}}\,{\bf ).}\label{rds}\rm
\, Given a probability space  $(\Om,\F,\p)$,  a measurable {\it random dynamical system} ({\bf RDS}) on $\big(\Rd, \mathcal{B}(\Rd)\big)$ over a  measurable dynamical system ({\bf DS}),  $\Theta:=\big(\Om,\F,\p,(\theta_t)_{t\in \R}\big)$, satisfying\hspace{.02cm}\footnote{\,Here, the notation $\theta\hspace{.02cm}\p=\p$ means that $\p(\{\om\in \Om: \theta_{t}\om \in A\}) = \p(\{\om\in \Om: \om \in A\}),  \;\forall \, A\in \mathcal{F}, \,t\in \Ic$; i.e., the semigroup $(\theta_t)_{t\in \Ic}$, $\theta_t:\Om\rightarrow \Om$, preserves the measure $\p$; we restrict the definition of the RDS to $\Ic=\R$.} $\theta\hspace{.02cm}\p=\p$, is a mapping $\Phi: \Ic\times\Om\times\Rd\rightarrow \Rd$ such that the following hold:
\begin{itemize}[leftmargin=0.9cm]
\item[(a)] $(t,\om,x)\mapsto \Phi(t,\om,x)$ is measurable for all $t\in \Ic\subseteq\R$,
\item[(b)] $\Phi(0,\om,\ccdot) = \textup{id}_{\Rd}$ for all $\om\in\Om$,   
\item[(c)] $\Phi(t+s, \om,\ccdot)= \Phi(t,\theta_s\hspace{.02cm}\om,\Phi(s,\om,\ccdot))$ for all $s,t\in\Ic$, $\om\in\Om$ \;\textup{(}{\it cocycle property}\textup{)},
\item[(d)]  $\Phi$ is continuous if $(t,x)\mapsto\Phi(t,\om,x)$ is continuous for all $t\in \Ic$, $x\in \Rd$,
\item[(e)]  $\Phi$ is smooth of class $\mathcal{C}^l,$ if \,$\Phi(t,\om,x)$ is $l$-times differentiable w.r.t.~ $x\in \Rd$, and the derivatives are continuous w.r.t.\! $(t,x)\in \Ic\times\Rd$.  
\end{itemize}
The canonical filtration on $(\Om,\F,\p)$ for the RDS is generated by $(\theta_t)_{t\in \R}$.
\end{definition}

\begin{definition}\rm
{\bf (Canonical DS for processes with stationary increments).} Consider a probability space $(\Om,\F,\p_\xi)$ with the measure $\p_\xi$ on $(\Om,\mathcal{F})$ induced by the law of a stochastic process with continuous time $\xi = (\xi_t)_{t\in \R}$, $\xi_t:\Om\rightarrow \Rd$. A process $\xi$ is said to have stationary increments if for any $t_0\leqslant \dots\leqslant t_{n}$,   $n\in \mathbb{N}_1$, the distribution of $(\xi_{t_1+t}-\xi_{t_0+t}, \dots, \xi_{t_n+t}-\xi_{t_{n-1}+t} )$ is independent of $t\in \R$; i.e.,
\begin{equation}
\theta(t)\p_\xi = \p_\xi \quad \textrm{for all} \;\;t\in \R,
\end{equation}
where $(\theta_t)_{t\in \R}$ is a semigroup of time shifts. The corresponding measurable dynamical system $\Theta:=(\Om,\F,\p_\xi,(\theta_t)_{t\in \R})$ is called the {\it canonical dynamical system} for the process with stationary increments; see, e.g. \cite[Appendix A.3]{Arnold} for details.
\end{definition}

\begin{prop}[\textit{\textbf{Canonical DS for Brownian motion/Wiener process}}]\label{canW}
For the Wiener probability space $(\Om,\F,\p)$ defined in \S\ref{gennot}, the canonical dynamical system   $\Theta=\big(\Om,\F,\p,(\theta_t)_{t\in \R}\big)$ for a stochastic process with stationary increments is given by 
\begin{equation}\label{shift}
\theta_t: \Om\rightarrow\Om, \qquad \theta_s\hspace{.02cm}\om(t) = \om(t+s)-\om(s) \qquad \forall\; s,t\in \R,  
\end{equation}
so that the set $\Omega=\mathcal{C}_0(\R,\R^m)$ is invariant w.r.t.~the shifts $(\theta_t)_{t\in \R}$. The canonical stochastic process $W_t(\om)=\om(t)$, $t\in\R$, with stationary independent increments is the Wiener process/Brownian motion $($with two-sided time$)$  which satisfies identically  
 \begin{equation}\label{tht_base} 
W_{t}(\theta_s\om) = W_{t+s}(\om)-W_s(\om) \qquad \forall\; s,t\in \R.
\end{equation}
\end{prop}
\noindent {\it Proof:} See \cite[Appendix A.3]{Arnold} for an outline or, e.g., \cite{Nelson67}.

\begin{rem} \rm In the sequel it will be more convenient to use (\ref{tht_base}) in the alternative form 
\begin{equation}\label{tht_base2}
W_{t}(\om) = W_{t+s}(\theta_{-s}\hspace{.02cm}\om)-W_s(\theta_{-s}\hspace{.02cm}\om) \qquad \forall\; s,t\in \R.
\end{equation}
\end{rem}

\smallskip
Assuming suitable regularity of the coefficients of autonomous SDE's, such as those in Theorem~\ref{SDEflow},  together with  adoption of two-sided stochastic calculus, the solutions of autonomous SDE's generate\footnote{\,The generation of an RDS from an SDE requires a `perfection of the  crude cocycle' associated with the SDE (see, e.g., \cite[Theorem 2.3.26]{Arnold}); here, this important technical nuance does not require an explicit discussion.}\label{perf_coc} an RDS over $\Theta$ (e.g.,~\cite{Arnold, Arnoldp, Elw78, Ikeda81, Kunita}). 
We will consider the non-autonomous dynamics of the SDE (\ref{NSDE}) with time-periodic coefficients as an RDS on a suitably extended~space.

\subsection{Time-periodic setting}
In the sequel,  we consider non-autonomous SDE's (\ref{NSDE}) on $\Rd$  with time-periodic coefficients; i.e., $b(t+\tau, \ccdot) = b(t,\ccdot), \;\sigma(t+\tau, \ccdot) = \sigma(t,\ccdot)$, $0<\tau<\infty$, satisfying the conditions in Theorem \ref{SDEflow}  so that (\ref{NSDE}) has global solutions generating the forward stochastic flow $\big\{\phi(t+s,s,\ccdot,\ccdot)\!: \;s\in \R,\, t\in \Rp\big\}$  such that, for all $s\in \R, \,t\in \Rp$,  
\begin{align}\label{Periodic_FL}
\phi(t+s+\tau, s+\tau, \om, \ccdot) = \phi(t+s,s, \theta_{\tau}\hspace{.02cm}\om,\ccdot) \quad \p\,\text{-\,a.s.}
\end{align}
The above property follows from the time-periodicity of the coefficients and the uniqueness of solutions of~(\ref{NSDE}). The relationship in (\ref{Periodic_FL}) is essential for constructing an RDS on $\Sc\times\Rd$ from solutions of (\ref{NSDE}) with time-periodic coefficients, which is important for asserting the existence and ergodicity of time-periodic measures supported on random time-periodic paths defined below.

\begin{definition}[{\bf Random periodic path of a stochastic flow} \cite{Feng11, Feng12, Zhao09}]\label{Chu}\;\rm
 A {\it random periodic path} of period $0<\tau<\infty$ {\it generated by a stochastic flow} $\big\{\phi(t+s,s,\ccdot,\ccdot)\!: \;\,s \in \R,\; t\in \Rp\big\}$ is a measurable function $S\,{:}\;\R\times\Om\rightarrow\Rd$ such that  for any  $s\in \R$ the following holds 
 \begin{equation}\label{per_sol_sde}
 S(\tau+s,\om) = S(s,\theta_{\tau}\hspace{.02cm}\om)\quad \text{and} \quad \phi(t+s,s,\om,S(s,\om)) = S(t+s,\om) \quad  \p\,\text{-\,a.s.} \quad \forall \; t\in \Rp.
 \end{equation} 
 \end{definition}

\medskip
\begin{definition}[{\bf Random periodic path of RDS}\label{rp_Phi} \cite{Feng18, Zhao09}]\rm
A {\it random periodic path} of period $0<\tau<\infty$ generated by an RDS, $\Phi: \Rp\times\Om\times\Rd\rightarrow \Rd$, is a measurable function $S: \R\times\Om\rightarrow \Rd$ such that for any $s\in \R$ and almost all $\om\in \Om$ the following holds 

\vspace{-.6cm}
\begin{align}\label{per_sol_rds}
S(\tau+s, \om) = S(s, \theta_\tau\hspace{.02cm}\om) \quad \text{and} \quad \Phi(t, \theta_s\hspace{.02cm}\om, S(s,\om)) = S(t+s, \om)\qquad \forall \;t\in \Rp.
\end{align}
\end{definition}

\vspace*{-.0cm}

\begin{exa}\rm
Let $ b: \Rd\rightarrow\Rd, \; d\geqslant 2$, be a globally Lipschitz vector field, and consider the deterministic flow $\{\psi(t,\,\cdot\,)\!: \;t\in\Rp\}$, defined via $\psi(t,\,\cdot\,)\equiv\phi(t,0,\,\cdot\,)$ and generated by the autonomous ODE
 \begin{align}\label{ODE1}
  \frac{dY_t}{dt} =b(Y_t), \qquad Y_0 = y\in \Rd, \quad   t\in \Rp.
 \end{align}
Assume that there exists a  periodic solution $\mathfrak{Y}: \R\rightarrow\Rd$ of  the ODE (\ref{ODE1}) of period $0<\tau<\infty,$ 
\begin{align*}
\mathfrak{Y}(\tau+s) = \mathfrak{Y}(s) \quad \text{and}\quad \psi(t, \mathfrak{Y}(s)) = \mathfrak{Y}(t+s), \quad s\in\R, \;\; t\in \Rp.
\end{align*}
 Consider the stochastic process $X_t(\om) = \mathfrak{Y}(t)+Z_t(\om),$  where $Z_t$ solves the following SDE
\begin{align}\label{sde111}
dZ_t   = \hat{b}(t, Z_t)dt+ \hat{\sigma}(t, Z_t)dW_t, \quad Z_0 = 0, \quad   t\in \Rp,
\end{align}
with time-periodic coefficients $$\hat{b}(t,z) := b(\mathfrak{Y}(t)+z)-b(\mathfrak{Y}(t)), \qquad  \hat{\sigma}(t,z) := \sigma(\mathfrak{Y}(t)+z).$$ If $\mathfrak{Z}(t,\om)$ is a random $\tau$-\,periodic solution of  (\ref{sde111}), then $S(t,\om) = \mathfrak{Y}(t)+\mathfrak{Z}(t,\om)$ is a random $\tau$-\,periodic solution of the autonomous SDE:
\vspace{.1cm}$$dX_t = b(X_t)dt+\sigma(X_t)dW_t, \; X_0 = \mathfrak{Y}(0), \;\;t\in \Rp.$$
\end{exa}

\vspace{.3cm}
\begin{exa}[Stochastic FitzHugh-Nagumo model with periodic current]\;\rm
Consider the following SDE with nontrivial random periodic solutions (see \cite{Feng12}) which has less restrictive conditions on the drift than those considered in the sequel: 
\begin{align}\label{FHNEX}
dX_t = AX_tdt+ b(t,X_t)dt+ \sigma(t) dW_{t-s}, \qquad X_{s} =x\in \R^2, \quad s,t\in \R, \;s\leqslant t, 
\end{align}
where 
\vspace{.2cm}
$$A = \begin{pmatrix}
1 & -1\\ a& -1
\end{pmatrix}, \quad b(t,x,y) = \begin{pmatrix}
-\frac{1}{3}x^3+B_1\sin (\tau \,t)\\ c
\end{pmatrix},  \quad \sigma(t) = \begin{pmatrix}
\sqrt{2\beta^{-1}}+B_2\cos(\tau\,t) & 0 \\ 0 & 0
\end{pmatrix}, 
$$ 

\vspace{.3cm}
\noindent with $ a<1$, $\beta>0$, $B_1,B_2, c\in \R$, $0<\tau<\infty$,  and $W_t = (W^1_t, \;0)^T$, where $W^1_t$ is a two-sided Wiener process on $\R$. Let $X_{t}^{s,x}(\om) = \phi(t,s,\om,x)$, $s\leqslant t$, be the solution of (\ref{FHNEX}) represented via 
\vspace{.2cm}
\begin{align*}
\phi(t,s,\om,x) = e^{A(t-s)}x+\int_{s}^te^{A(t-\zeta)}b\big(\zeta, \phi(\zeta,s,\om,x)\big)d\zeta+\int_{s}^t e^{A(t-\zeta)}\sigma(\zeta)dW_{\zeta-s}(\om),
\end{align*}

\vspace{.1cm}
\noindent where $x\mapsto e^{A(t-s)}x$ is the solution of the linear ODE
\vspace{.3cm}
\begin{align*}
\frac{dY_t}{dt} = AY_t, \qquad Y_{s} = y\in\R^2, \quad s,t\in \R, \;s\leqslant t.
\end{align*}

\medskip
\noindent Now, consider the projections $P^-: \R^2\rightarrow E^-, \; P^{+}: \R^2\rightarrow E^+$, where the linear subspaces are 
\begin{align*}
E^- = \text{span}\{y\in \R^2: Ay= -\lambda y\}, \; \;  E^+ = \text{span}\{y\in \R^2: Ay= \lambda y\}, \; \lambda: = \sqrt{1-a}.
\end{align*}

\smallskip
\noindent The process $S(t,\om)$ defined by
\vspace{.2cm}
\begin{align*}
\notag S(t,\om) &= \int_{-\infty}^te^{A(t-\zeta)}P^-b(\zeta,S(\zeta,\om))d\zeta-\int_{t}^{\infty}e^{A(t-\zeta)}P^{+}b(\zeta,S(\zeta,\om))d\zeta
\\[.3cm] &\hspace*{.6cm}+\int_{-\infty}^te^{A(t-\zeta)}P^{-}\sigma(\zeta)dW_{\zeta-s}(\om) - \int_{t}^{\infty}e^{A(t-\zeta)}P^{+}\sigma(\zeta)dW_{\zeta-s}(\om),
\end{align*}
is a random $2\pi/\tau$-periodic solution of the flow generated by the SDE (\ref{FHNEX}); see, e.g.,~\cite{Feng12,Cherub17}. 
\end{exa}


\newpage
\section{ Time-periodic ergodic measures for dissipative  SDE's}\label{Erg_S}
In this section we consider a class of non-autonomous SDE's (\ref{NSDE}) which generate stable random periodic paths. First, in \S\ref{permeas} we prove the existence of a unique stable random periodic solution for a class of {\it `dissipative\hspace{-.03cm}'}\hspace{.03cm}\hspace{.06cm}\footnote{\,See (\ref{diss_linear_growth}) for one such class which we focus on in this work.} SDE's with time-periodic coefficients, and we assert the existence of  time-periodic measures  induced by such dynamics (Theorem~\ref{Rand_SOL}).  Ergodicity (in an appropriate sense, and under typical regularity conditions) of these time-periodic measures are established in Theorem \ref{Ps_erg} of \S\ref{permes_erg}. We conclude  with an example of a periodically forced stochastic Lorenz model, which is then used in \S\ref{Linear response_per} to illustrate the utility of  fluctuation-dissipation formulas for time-periodic measures when considering the linear response of the dynamics to small perturbations.

\subsection{Preliminaries, definitions, and assumptions}\label{permes_ass}
\addtocontents{toc}{\protect\setcounter{tocdepth}{2}}
First, we recall the notion of a time-periodic probability measure which will be needed throughout the reminder of this paper.

\vspace{-0.1cm}
\begin{definition}[{\bf Time-periodic probability measure} \cite{Feng18}]\;\label{Rand_pm}\rm
 A  measure-valued map given by  $t\mapsto\mu_{s+t}\in \PP(\Rd)$ and induced by the family $(\mathcal{P}^{*}_{s,s+t})_{t\in \Rp}$, $s\in \R$,  defined in  (\ref{P*mu}) is referred to as a {\it time-periodic probability measure} of period $0<\tau<\infty,$ if the following holds for any  $s\in \R$
 \begin{equation}\label{10.11}
  \mu_{s+t} =\mathcal{P}_{s,\,s+t}^*\,\mu_s \quad \text{and}\quad \mu_{s+\tau} = \mu_s \qquad \forall\;  t\in \Rp.
 \end{equation}
Furthermore,  $\mu_{s+t}\in \PP(\Rd)$, $s\in \R$, $t\in \Rp$, is called a {\it time-periodic measure with the minimal $($or fundamental$\,)$ period $\tau,$} if $\tau$ is the smallest strictly positive number such that \eqref{10.11} holds\footnote{\,Sufficient conditions for the existence of the minimal period, which are satisfied here, are established in \cite[\S 5]{Feng18}.}\label{mintau}. 
 \end{definition}

  \begin{prop}\label{Srem}
 Let  $ S:\,\R\times\Om\rightarrow \Rd$  be a random periodic path (\ref{per_sol_sde}) of a stochastic flow $\big\{\phi(t+s,s,\ccdot,\ccdot)\!:\; s \in \R,\; t\in \Rp\big\}$ on $\big(\Rd,\mathcal{B}(\Rd)\big)$ and consider a family  of probability measures 
\begin{align*}
  \mu_{s+t}(A) :=  \p\big(\{\om: S(s+t,\om)\in A\}\big) \qquad \forall\; s\in \R,\,t\in \Rp,\;A\in \Bb(\Rd).
 \end{align*} 
Then, the family $(\mu_{s+t})_{s\in \R, t\in \Rp}$ consists of $\tau$-periodic probability measures on $\Rd$.  
\end{prop}
 \noindent {\it Proof:} This follows by a direct calculation combined with the properties of a random periodic path (\ref{per_sol_sde}), since  for all  $s\in \R,\,t\in \Rp,\;A\in \Bb(\Rd)$, we have 
 \begin{align*}
 \hspace{2cm}\mu_{s+\tau}(A)= \p\big(\{\om: \;S(s+\tau,\om)\in A\}\big)
&= \p\big(\{\om: \;S(s, \theta_{\tau}\hspace{.03cm}\om)\in A\}\big) \notag\\ &= \p\big(\{\om: \;S(s,\om)\in A\}\big)=\mu_s(A). \hspace{2cm}\qed
\end{align*}

The above results will be generalised to the dynamics in the extended state space in \S\ref{extS}.

\subsubsection{\bf Dynamics on the extended state space}\label{extS} A useful way of examining ergodicity of time-periodic measures induced by  non-autonomous SDE's with time-periodic coefficients of period $\tau$ is to lift the original dynamics from $\Rd$ to the extended state space  $\Sc\times \Rd$, $\Sc\simeq \R\Mod\tau$, so that the resulting `lifted' SDE is autonomous. Such a representation of the original dynamics does not necessarily simplify the formulation of the problem but the flows of the lifted solutions generate a cocycle\footnote{\,See Definition~\ref{rds}.}  in the  skew-product variables on   $\Sc\times \Rd$; we refer to this extended state space as  the   {\it `flat cylinder'}.  Then,  the lifted random periodic paths (\ref{per_sol_sde}) of the stochastic flow induced by the non-autonomous SDE (\ref{NSDE}) can be associated with random periodic paths (satisfying (\ref{per_sol_rds})) of an RDS (see Definition \ref{rds}) generated by the lifted flow in  the skew-product representation  on  $\Sc\times \Rd$; this fact  allows to prove ergodicity (in an appropriate sense) of time-periodic measures supported on the random periodic paths on the fibre bundle\footnote{\,See, e.g., \cite{sin} for a detailed description of such structures on spaces of probability measures.  }  ${\PP}\big(\Sc\big)\times \PP\big(\Rd\big)$. 
 
 \smallskip
 To this end,  consider the solutions of the SDE (\ref{NSDE}) satisfying the conditions  of Theorem~\ref{SDEflow} and assume that  the coefficients of (\ref{NSDE}) are time-periodic with period $0<\tau<\infty$; we recast the solutions of (\ref{NSDE}) as an extended process $\tilde X_t(\om) = \big(t,X^{s,x}_t(\om)\big)^{T}$ in the {\it skew-product} representation on $\R\times\Rd$ satisfying  
 \begin{align}\label{NSDE1}
 d\tilde X_t &= \tilde b\big(\tilde{X_t}\big)dt+ \tilde{\sigma}\big(\tilde{X_t}\big) d\tilde{W}_{t-s}, \qquad \tilde X_s = (s,x)\in \R\times\Rd,\;\; s\leqslant t,
 \end{align}
 where $\tilde{W}_{t-s}(\om) = \big(0, W_{t-s}(\om)\big)$, $\om\in \Om$, and $W_{t-s}$ is the $m$-dimensional Brownian motion for the two-sided time~(see~\S\ref{gennot} or \cite{Arnold}), and $\tilde b: \R^{d+1}\rightarrow\R^{d+1}$, $\tilde \sigma: \R^{d+1}\rightarrow\R^{(d+1)\times(m+1)}$,  
so that 
 \begin{align}\label{LiftExp}
  d\begin{pmatrix} \zeta_{\hspace{.02cm}t} \\ X_t\end{pmatrix} &= \begin{pmatrix}
  1\\ b\big(\zeta_{\hspace{.02cm}t}, X^{s,x}_t\big)
  \end{pmatrix}dt + \begin{pmatrix}
 0& 0\\ 0 &\sigma\big(\zeta_{\hspace{.02cm}t}, X_t^{s,x}\big)
  \end{pmatrix} d\tilde{W}_{t-s}, \;\; \;s\leqslant t.
  \end{align}
The dynamics in (\ref{NSDE1}) or (\ref{LiftExp}) can be represented in a more convenient form for the subsequent derivations by setting $t\rightarrow t+s$, so that 
 \begin{align}\label{NSDE11}
 d\tilde X_{t+s} &= \tilde b(\tilde{X}_{t+s})dt+ \tilde{\sigma}(\tilde{X}_{t+s}) d\tilde{W}_{t+s}, \quad \tilde X_s = (s,x)\in \R\times\Rd, \,t\in \Rp,
 \end{align}
where $\tilde{W}_{t+s} = \tilde{W}_{t+s}(\theta_{-s}\hspace{.02cm}\om)$ is the Brownian motion satisfying  (\ref{tht_base2}). Finally, given the form of the coefficients $\tilde b$, $\tilde \sigma$, it is convenient to consider the dynamics induced by (\ref{NSDE11}) on the flat cylinder  $\Sc\times \Rd$, where $\Sc\simeq \R\Mod\tau$.

\smallskip
The RDS   associated with the lifted dynamics~(\ref{NSDE11}) is generated in the skew-product representation (see, e.g., \cite{sin,Arnold}) on  $\Sc\times\Rd$ via 
\begin{align}\label{Lift}
\tilde\Phi\big(t, \om, \tilde x\big) := \big( t+s \Mod\tau, \;\;\phi(t+s,s,\theta_{-s}\hspace{.02cm}\om,x)\big) \quad \;\forall\;\tilde x:=(s,\,x)\in \Sc\times \Rd, \,t\in \R^+.
\end{align}
 The cocycle property\footnote{\,To be more accurate, the so-called `crude' cocycle property can be easily verified from the flow induced by the SDE, and the crude cocycle needs to be  `perfected' in order to generate an RDS over the DS for the Brownian motion  (see, e.g., \cite[Theorem 2.3.26]{Arnold}); here, this important technical nuance does not require an explicit discussion.}  
of $\tilde\Phi$ in (\ref{Lift}), i.e., $\tilde{\Phi}(t+r, \om,\ccdot)= \tilde{\Phi}\big(t,\theta_r\hspace{.02cm}\om,\tilde \Phi(r,\om,\ccdot)\big)$ for all $r,t\in \Rp$, and a.a.~$\om\in \Om$, can be verified
 by recalling that $t+r\,\mmod\tau = t+r-k\tau,$ where $k=\lfloor\frac{t+r}{\tau}\rfloor$, and utilising (\ref{Periodic_FL}). Note that, unless (\ref{NSDE}) is autonomous,  $\{\phi(t+s,s,\theta_{-s}\hspace{.02cm}\om,\ccdot)\!: \,s\in \R, \, t\in \Rp \}$  does not have the cocycle property, and hence it does not generate an RDS on $\Rd$. 
 The RDS representation of the non-autonomous dynamics of the SDE (\ref{NSDE}) will be useful in \S\ref{permes_erg} when considering the ergodicity of measures supported on random periodic paths, and in the discussion of the linear response in \S\ref{LRsp}.

\smallskip
The transition kernel and transition evolutions (see Definitions \ref{Mkernel} and \ref{trans_ev}) on $\Sc\times \Rd$
are constructed as follows.  For any $\tilde x:=(s,\,x)\in \Sc\times \Rd$, $t\in \Rp$,  and $\tilde{A}\in \Bb\big(\Sc\big)\times\Bb\big(\Rd\big)$,  the {\it transition kernel} $\tilde{P}(\tilde{x};t,\ccdot)$ associated with  $\tilde \Phi$   is given~by (see Definition~\ref{Mkernel})
 \begin{align}\label{PPhi}
\notag \tilde{P}\big(\tilde{x}; t, \tilde{A}\hspace{.04cm}\big) :\!&= \p\big(\{\om: \tilde\Phi(t,\om,\tilde{x})\in \tilde{A}\,\}\big)\\ \notag &= \p\big( \{\om: (t+s\Mod\tau, \;\phi(t+s,s,\theta_{-s}\hspace{.02cm}\om,x))\in \Jc\times A\,\}\big)\\  &=  \delta_{(t+s\Mod\tau)}(\Jc)\times P\big(s,x;t+s,A\big), 
 \end{align}
for all $\tilde{A} \equiv \Jc\times A\in \Bb\big(\Sc\big)\times\Bb\big(\Rd\big)$ and the transition kernel $P$ defined in (\ref{phiP}).

 The {\it transition evolution} $(\tilde{\mathcal{P}}_t)_{t\in \Rp}$ induced by $\tilde\Phi$ and its dual $(\tilde{\mathcal{P}}^*_{t})_{t\in \Rp}$  are given by\footnote{\,We abuse the notation whereby $\mathbb{M}_\infty\big(\Sc\times\Rd\big)$ denotes functions measurable w.r.t.~$\mathcal{B}(\Sc)\times\mathcal{B}(\Rd)$; see~\S\ref{gennot}.} 
\begin{align}
\tilde{\mathcal{P}}_t\varphi(\tilde{x}) &:= \int_{\Sc\times\Rd}\varphi(\tilde{y})\tilde{P}(\tilde{x}; t, d\tilde{y}) \hspace{.8cm} \forall \;\varphi\in \mathbb{M}_\infty\big(\Sc\times\Rd\big),\label{tldP}\\[.2cm]
\tilde\mu_{t+r}(\tilde A)=\big(\tilde{\mathcal{P}}^*_t\hspace{.03cm}\tilde{\mu}_r\big)(\tilde{A}) &:= \int_{\Sc\times\Rd}\tilde{P}(\tilde{x}; t, \tilde{A})\tilde{\mu}_r(d\tilde{x}) \quad\; \forall\;\tilde{\mu}_r\in \PP\big(\Sc)\otimes\PP(\Rd\big),\;r\in \Rp,\label{tldP*}
\end{align}
with the short-hand notation  $\tilde\mu_r(d\tilde x) = \delta_{(r\Mod \tau)}(s)ds\otimes \mu_{r}(dx)$ for skew product probability measures in the fibre bundle  $\PP(\Sc)\otimes\PP(\Rd)$, where 
$\tilde\mu_r\in \PP(\Sc)\otimes\PP(\Rd)$ and $\mu_r\in \PP(\Rd)$; see, e.g.,~\cite{sin} for more details on  skew-product fibre bundles on spaces of  probability measures. Extension  of (\ref{tldP}) to arbitrary functions $\mathbb{M}\big(\Sc\times\Rd\big)$ can be carried out in a standard way.

Given (\ref{Lift}), if the RDS $\big\{\tilde\Phi(t,\ccdot,\,\cdot\,)\!:\, t\in \Rp\big\}$ on $\Sc\times\Rd$ has a random periodic path $t\rightarrow \tilde S(t,\om)$ of period $0<\tau<\infty$, then it has the form  $\tilde S(t,\om) = \big(t\Mod\tau,\;S(t,\om)\big)$, 
 and $t\rightarrow S(t,\om)$ is a random periodic path of  $\big\{\phi(t+s,s,\ccdot,\,\cdot\,)\!:\, t\in \Rp\big\}$ on $\Rd$ (see Definition \ref{rp_Phi}).

\begin{lem}\label{semi_P} The transition evolutions $(\tilde{\mathcal{P}}_t)_{t\in \Rp}$ and $(\tilde{\mathcal{P}}^*_{t})_{t\in \Rp}$
have a semigroup structure. In particular,  for  $\tilde\mu_t = \delta_{(t\Mod \tau)}\otimes \mu_{t}$ in the fibre bundle $\PP\big(\Sc)\otimes\PP(\Rd\big)$ the following holds:  
\begin{equation*}
\tilde\mu_{t+r+u}=\tilde{\mathcal{P}}^*_{t+r}\hspace{.02cm}\tilde{\mu}_u = \tilde{\mathcal{P}}^*_{t}\big(\tilde{\mathcal{P}}^*_{r}\tilde{\mu}_u\big)=\tilde{\mathcal{P}}^*_{t}\hspace{.02cm}\tilde{\mu}_{r+u} \qquad\;\; \forall\; r,t,u\in \Rp.
\end{equation*}
If the  RDS $\big\{\tilde\Phi(t,\ccdot,\,\cdot\,)\!:\, t\in \Rp\big\}$ on $\Sc\times\Rd$ has a random $\tau$-periodic path $t\rightarrow \tilde S(t,\om)$, then all skew-product probability measures in the family  $(\tilde \mu_t)_{t\in \Rp}$, $\tilde\mu_t\in \PP\big(\Sc\big)\otimes\PP\big(\Rd\big)$, supported on such a path are $\tau$-periodic, i.e., 
\begin{equation*}
\tilde \mu_{t+r} =\tilde{\mathcal{P}}_{r}^*\,\tilde\mu_t,\qquad \tilde\mu_{t+\tau} = \tilde\mu_t \qquad \forall \,t\in \Rp,
\end{equation*}

\vspace{-.2cm}
\noindent and 
\vspace{-.3cm}
\begin{equation*}
\hspace{1.5cm}\tilde\mu_t(\tilde A) = \mu_t(A_t),\hspace{1cm} \forall\; \tilde A\in \Bb\big(\Sc\big)\times\Bb\big(\Rd\big), \quad A_t =\big\{x\in \Rd: (t\Mod\tau,\;x)\in \tilde{A}\,\big\}. 
\end{equation*}
Every such $\tau$-periodic measure is invariant under the discrete dynamics induced by $(\tilde{\mathcal{P}}^*_{n\tau})_{n\in \mathbb{N}_0}$,~i.e., 
\begin{equation*} 
   \tilde{\mathcal{P}}^*_{n\tau}\,\tilde\mu_t = \tilde\mu_t  \qquad \forall\; n\in \mathbb{N}_0, \;t\in \Rp. 
\end{equation*}
\end{lem}

\noindent {\it Proof}: The first claim is a direct consequence of (\ref{tldP*}), and the proof follows either by using the cocycle property of $\tilde\Phi$ in the first line  of (\ref{PPhi}) or by utilising the Chapman-Kolmogorov equation~(\ref{ChK}) for $P$ in the last line of (\ref{PPhi}). 

Regarding the second claim, consider  measures supported on the random periodic~path  $\tilde S$ 
\begin{equation}\label{tilmu2mu}
 \tilde{\mu}_{t+r}(\tilde{A}\hspace{.03cm}) := \p\big(\{\om: \; \tilde S(t+r,\om)\in \tilde{A}\,\}\big), \qquad \tilde{A}\in \Bb\big(\Sc\big)\times\Bb\big(\Rd\big). 
 \end{equation}
Since $\tilde S$ is a random periodic path of the RDS $\tilde \Phi$, we have for all $\tilde{A}\in \Bb\big(\Sc\big)\times\Bb\big(\Rd\big)$ that 
 \begin{align*}
\tilde \mu_{t+r}(\tilde A\hspace{.03cm}) &= \p\big(\{\om: \;\tilde S(t+r,\om)\in \tilde{A}\,\}\big)=
\p\big(\{\om: \;\tilde \Phi(r, \theta_t\hspace{.02cm}\om, \tilde S(t,\om))\in \tilde{A}\,\}\big) = \tilde{\mathcal{P}}^*_r\tilde\mu_t(\tilde A\hspace{.03cm}), 
 \end{align*}
 for all $r,t\in \Rp$ by the general properties the random periodic path (\ref{per_sol_rds});
this could also be obtained directly from (\ref{tldP*}) by using the invariance of $\tilde S$ under the action of $\tilde \Phi$. 
Moreover,  
 \begin{align*}
\tilde \mu_{t+\tau}(\tilde A\hspace{.03cm}) &= \p\big(\{\om: \;\tilde S(t+\tau,\om)\in \tilde{A}\,\}\big)
=\p\big(\{\om: \;\tilde S(t,\theta_\tau\hspace{.02cm}\om)\in \tilde{A}\,\}\big) = \tilde\mu_t(\tilde A\hspace{.03cm}),
 \end{align*}
 by the property (\ref{per_sol_rds}). Thus  $\tilde{\mu}_t$ is a $\tau$-periodic measure for the RDS $\big\{\tilde\Phi(t,\ccdot,\ccdot)\!:\, t\in \Rp\big\}$ on $\Sc\times\Rd$ which is supported on the random periodic path $\tilde S$.  
 The last two claims are simple consequences of the properties established above and the skew-product structure of probability  measures supported on random periodic paths. 
 \qed

\smallskip
   In the following sections, after outlining some general assumptions, we will investigate the existence and uniqueness of stable random periodic paths of the RDS  $\big\{\tilde\Phi(t,\ccdot,\ccdot)\!:\, t\in \Rp\big\}$, and we will prove the ergodicity of  probability measures associated with the dynamics of the skew-product lift (\ref{NSDE11}) of  the  dynamics in (\ref{NSDE}) under some standard regularity assumptions.

\subsubsection{{\bf Assumptions}} Throughout, we assume that the SDE (\ref{NSDE}) with time-periodic coefficients of period $0<\tau<\infty$ satisfies the conditions of Theorem \ref{SDEflow}, so that  (\ref{NSDE}) has global solutions.

\smallskip
In order to establish the existence of stable random periodic paths in \S\ref{permeas}, we will require the following assumption: 

\vspace{-.1cm}\begin{Assum}\label{A2.1}\rm
Let $V\in \mathcal{C}^{1,2}\big(\R{\times}\,\Rd;\Rp\big)$ s.t. \!$V(t,0) = 0$ for all $t\in \R$,  satisfy the following:   
\begin{itemize}[leftmargin = 0.9cm] 
\item[(i)]  There exist $\lambda\in L^1(\R;dt),$ and a constant $\mathfrak{C}\geqslant 1$, such that for some $1<p<\infty$ and all $\xi, \eta\in L^{p}(\Om, \F_{-\infty}^{\hspace{0.03cm}t}, \p)$, we have 
\begin{align}\label{Lyap_f}
\begin{cases}
\E\vert \xi\vert^p\leqslant \E\big[V(t,\xi)\big]\leqslant \mathfrak{C}\,\E\vert \xi\vert^p<\infty,\\[.2cm]
\E\big[\mathcal{L}^{(2)}V(t,\xi-\eta)\big]\leqslant\lambda(t) \E\big[V(t,\xi-\eta)\big], 
\end{cases}
\end{align} 
where $\LG^{(2)}$ is the two-point generator defined in (\ref{2pp}) and associated with the SDE (\ref{NSDE}). 

\item[(ii)] There exists $\bar{\lambda}>0$ such that 
\begin{align}\label{Rate1}
\limsup_{(t-s)\rightarrow\infty}\frac{1}{t-s}\int_{s}^t\lambda(u)du<-\bar{\lambda}<0.
\end{align}

\item[(iii)] For the one-point motion  $t\mapsto\phi(t,s,\om,\xi)$  induced by (\ref{NSDE}) for $\om\in \Om$, $\xi\in \Rd$, and $ s\leqslant t$,  there exists $0<\mathfrak{D}<\infty$ independent of $s,t\in \R$ such that\footnote{\,This condition can be replaced by a stronger but a more concrete constraint on the global existence of the $p$-th absolute moment of~$\phi$; see Lemma \ref{lem_app1} in Appendix \ref{app_hasminski}.}  for all $\xi\in L^{p}(\Om, \F_{-\infty}^{s}, \p)$
\begin{align} \label{Temp1}
\limsup_{(t-s)\rightarrow \infty}\E\big[V\big(t,\phi(t,s,\xi) -\xi\big)\big]\leqslant \mathfrak{D}, 
\end{align}
where $\E\big[V(\phi(t,s,\xi))\big]:=\int_\Om V\big(\phi(t,s,\om,\xi)\big)\p(d\om)$. 
 
\end{itemize}
\end{Assum}

\noindent As pointed out later (Remark \ref{disc_ass} in \S\ref{permes_erg}), this assumption is not strictly required for proving  ergodicity of $\tau$-periodic probability measures. However, without showing the existence of random periodic paths (in this case, stable random periodic paths), the existence of the skew-product $\tau$-periodic measures would have to be assumed a priori alongside the ergodicity of $\tilde\mu_t$ for all fixed $t\in \Sc$ with respect to the discrete transition evolution $(\tilde{\mathcal{P}}^*_{n\tau})_{n\in \mathbb{N}_0}$, as done in~\cite{Feng18}.

\begin{rem}\label{rem_hasm_lyap}\rm\mbox{}
\begin{itemize}[leftmargin=0.7cm]
\item[(a)] An important class of coefficients satisfying Assumption \ref{A2.1},  which yield  global solutions of  (\ref{NSDE}) are specified in Appendix \ref{app_hasminski}. In particular, we might take $b(t,\ccdot)\in \tilde{\mathcal{C}}^{1,\delta}(\Rd)$ and $\sigma_k(t,\ccdot)\in \bar{\bar{\mathcal{C}}}^{1,\delta}(\Rd)$, $0<\delta\leqslant 1$, $k=1,\dots,m$, satisfying  the following {\it `dissipative'} condition
\begin{align}\label{diss_linear_growth}
\langle b(t,x), x\rangle \leqslant L_{b_1}(t) -L_{b_2}(t)\vert x\vert^2, \quad \Vert \sigma(t,x)\Vert^2_{\textsc{hs}}\leqslant L_\sigma(t)\big( 1+\vert x\vert^2\big), 
\end{align}
where $L_{b_1}, L_{b_2}, L_\sigma\in \mathcal{C}_\infty(\R, \Rp)$. Here, $\langle \,\cdot,\cdot\,\rangle$ denotes the dot product on $\R^d$ and $\| \cdot\|_{\textsc{hs}}$ denotes the Hilbert--Schmidt norm (aka Frobenius norm) defined by $\| A\|^2_{\textsc{hs}} = \text{trace}(AA^T)$.
Condition (\ref{Temp1}) is satisfied for (\ref{diss_linear_growth}) when (see Lemma \ref{lem_app1} in Appendix \ref{app_hasminski})
\begin{align}\label{diss_linear_growth_mom}
\inf_{t\in \R} \Big(L_{b_2}(t)-2^{\frac{p}{2}-1}L_{b_1}(t)-\frot (2^{\frac{p}{2}-1}+1)L_\sigma(t)(p-1)\Big)>0,
\end{align}
and it also leads to the global existence of the $p$\,-th absolute moment of the law of the associated SDE; tighter bounds can be obtained for $p=2,3$ as shown in Propositoin \ref{shrpbnd} in Appendix~\ref{app_hasminski}. Condition (\ref{Temp1}) is reminiscent of the Ha\'sminskii-type regularity condition \cite{Has12} for the existence and uniqueness of global solutions of SDE's; sufficient conditions for verification of Ha\'sminskii's conditions require the existence of real-valued functions $L_b(\cdot), L_\sigma(\cdot)\in \mathcal{C}_\infty\big(\R;\Rp\big)$  such that 
\begin{align}\label{Linear growth}
\big\langle b(t,x), x\big\rangle\leqslant L_b(t)\big(1+|x|^2\big), \quad \|\sigma(t,x)\|^2_{\textsc{hs}}\leqslant L_\sigma(t)\big(1+|x|^2\big).
\end{align}
Coefficients satisfying (\ref{diss_linear_growth}) also satisfy (\ref{Linear growth}), since for some $L_b\in \mathcal{C}_\infty\big(\R, \Rp\big)$ we have 
$$ L_{b_1}(t) -L_{b_2}\vert x\vert^2\leqslant L_b\big(1+|x|^2\big).$$

\vspace{.2cm}\item[(b)] Construction of the Lyapunov function $V$ satisfying Assumption \ref{A2.1} is often not straightforward. However, one can construct (e.g., \cite{Has12,Hutz15, Majka17}) a polynomial Lyapunov function growing at infinity as $\vert x\vert^{2N}, \; N\in\N_1$, for a broad class of SDE's whose coefficients $b(\ccdot,x),\,\sigma(\ccdot,x)$ are continuous, and  $b(t,\,\cdot\,)\in \tilde{\mathcal{C}}^{1,\delta}\big(\Rd\big)$, $\{\sigma_{k}(t,\,\cdot\,)\}_{1\leqslant k \leqslant m}\in \bar{\bar{\mathcal{C}}}^{1,\delta}\big(\Rd\big)$ are such that 
\vspace{.1cm}
\begin{align}\label{Dep1}
\begin{cases}
\big\langle b(t,x)-b(t,y), x-y\big\rangle\leqslant -K_t\vert x-y\vert^2,\\[.2cm]
\|\sigma(t,x)-\sigma(t,y)\|_\textsc{hs} \leqslant L_t\vert x-y\vert,\\[.2cm]
 \sup_{\,t\in\R}\big\{\vert b(t,0)\vert +\|\sigma(t,0)\|_\textsc{hs}\big\}<\infty,
\end{cases}
\end{align}

\vspace{.2cm}
\noindent where $0< L_t, K_t<\infty$,  and  
\begin{align}\label{Dep2}
\limsup_{(t-s)\rightarrow\infty}\frac{1}{t-s}\int_s^t\lambda(u)du<0,
\end{align}
with  $\lambda(t) = -K_t+\frac{(p-1)}{2}\hspace{.03cm}p\hspace{.03cm}L^2_t$ for some $1<p<\infty.$
The function $K_t$ is defined by 
\begin{align*}
K_t = \liminf_{R\rightarrow \infty}K_t(R),
\end{align*}
where $K_t: \R\rightarrow\R$ is a Borel function defined by 
\begin{align*}
K_t(R) = \inf\Big\{-\frac{\big\langle b(t,x) - b(t,y), x-y\big\rangle}{\vert x-y\vert^2}: \;\;\vert x-y\vert=R\Big\}.
\end{align*}
Many important classes of SDE's driven Levy processes (including the Brownian motion) satisfy the dissipative conditions (\ref{Dep1}) - (\ref{Dep2}); see \cite{Has12,Hutz15, Majka17} for more details. 

\end{itemize}
\end{rem}

\medskip
\noindent In order to study  the ergodicity of  $\tau$-periodic measures, we will require variants of the following standard conditions (e.g.,~\cite{Hairer11}) to be satisfied:
\begin{itemize}[leftmargin=0.85cm]
\item[(i)] Relative compactness property of the transition kernel $P$ in (\ref{Mkernel}).

\vspace{.1cm}\item[(ii)] Irreducibility of the transition kernel. 

\vspace{.1cm}\item[(iii)] Strong Feller property\footnote{ \;The transition evolution $(\mathcal{P}_{s,t})_{t\geqslant s}$ on a complete separable metric space $\XX$ has strong Feller property if for $\varphi\in \mathbb{M}_\infty(\XX)$, one has $\mathcal{P}_{s,t}\varphi\in \mathcal{C}_\infty(\XX), \; \forall s\leqslant t$, i.e., $\mathcal{P}_{s,t}: \mathbb{M}_\infty(\XX)\rightarrow\mathcal{C}_\infty(\XX)$, $\forall s\leqslant t$. } of the transition evolution  $({\mathcal{P}}_{s,t})_{t\geqslant s}$ (\ref{calP}).

\end{itemize}

\medskip
Thus,  we will require the following version of the H\"ormander condition (e.g., \cite{Nulart,Malliavin}) in \S\ref{permes_erg} in addition to Assumption \ref{A2.1}:

\newpage
\vspace{-0.0cm}\begin{Assum}\label{Homa} \rm 
Denote by $\sigma_k$, $1\leqslant k\leqslant m$, the columns of $\sigma$ in (\ref{NSDE}), and assume that the following are satisfied for all $t\in \R $:
\begin{itemize}[leftmargin=0.8cm]
\item[(i)] $b(t,\ccdot)\in \tilde{\mathcal{C}}^{\infty}(\Rd)$ and  $t\mapsto b(t,\ccdot)$ is differentiable. 
\item[(ii)] $\sigma_k(t,\ccdot)\in \bar{\bar{\mathcal{C}}}^{\infty}(\Rd)$, $t\mapsto \sigma_k(t,x)$ is  differentiable, and  
\begin{align}\label{Hom3}
\big| \partial_tD_x^\beta\sigma_k(t,x)\big| \leqslant C<\infty, \qquad (t,x)\in \R\times\Rd.
\end{align}
for every multi-index $\beta$.
 
\item[(iii)]  $\text{Lie}\big(\sigma_1(t,\,\cdot\,),\cdots, \sigma_m(t,\,\cdot\,)\big) = \Rd$, for all $t\in \Ic$,  where 
 $$\text{Lie}\big(\sigma_1(t,x),\cdots, \sigma_m(t,x)\big):= \textrm{span}\big\{\sigma_i, [\sigma_i, \sigma_j],[\sigma_i,[\sigma_j,\sigma_k]],\cdots, \; 1\leqslant i,j,k\leqslant m\big\},$$
and $[F, G\hspace{0.04cm}]$ is the Lie bracket between the vector fields $F$ and $G$ defined by 
\begin{align*}
 [F,G\hspace{0.04cm}](t,x):= D_xG(t,x)F(t,x)-D_xF(t,x)G(t,x).
\end{align*}

\end{itemize}
\end{Assum}

\subsection{Existence and uniqueness of  time-periodic measures on stable random periodic paths} \label{permeas} 
Given the preliminary results and assumptions outlined in~\S\ref{permes_ass}, we have the following result on the existence of a $\tau$-periodic measure (Definition \ref{Rand_pm}) for the lifted SDE~in~(\ref{NSDE1}).

\begin{theorem}\label{Rand_SOL}
Consider the forward stochastic flow $\{\phi(t,s,\ccdot,\ccdot)\!\!: \,s,t\in \R, \,s\leqslant t\}$ generated by the SDE in  (\ref{NSDE}) with time-periodic coefficients of period $0<\tau<\infty$, and  satisfying the conditions of Theorem \ref{SDEflow}. If Assumption~\ref{A2.1} holds,
 there exists a family $(\tilde{\mu}_t)_{t\in \Rp}$ of $\tau$-periodic skew-product probability measures, $\tilde \mu_t = \delta_{(t\Mod\tau)}\otimes\mu_t$, $\tilde{\mu}_t\in \PP\big(\Sc\big)\otimes\PP\big(\Rd\big)$, given by 
\begin{equation}\label{permes_rds}
\tilde{\mu}_t(\tilde{A}) := \p\big(\{\om: \;\tilde{S}(t,\om)\in \tilde{A}\,\}\big), \qquad t\in \Rp, \; \tilde{A}\in \Bb\big([0,\tau)\big)\times\Bb\big(\Rd\big),
\end{equation}
which are supported on a unique random periodic path $\tilde S$ of the RDS $\{\tilde\Phi(t,\ccdot,\ccdot)\!: \;t\in \Rp\}$  on $\Sc\times\Rd$ and defined  in (\ref{Lift})  . 
\end{theorem}
\noindent {\it Proof.}  First, for $\xi\in L^{p}(\Om, \F_{-\infty}^{s}, \p),\; 1< p<\infty$, where $\displaystyle \F_{-\infty}^s:=\textstyle{\bigvee_{r\leqslant s}}\,\F_{r}^s,$ we show that $\{\phi(t,s,\om,\xi)\!:\; s, t\in \R,\,s\leqslant t\}$ converges to a random process $S(t,\om)\in \Rd$ almost surely as $s\rightarrow -\infty$, and that $S(t,\om)$ is bounded and independent of $\xi$.  Next, we show that $t\mapsto S(t,\om)$ is a unique stable random periodic path of period $0<\tau<\infty$ for $\{\phi(t,s,\om,\ccdot)\!: \,s,t\in \R, \;s\leqslant t\}$. Finally, we conclude that the law of the random periodic  path $\tilde{S}(t,\om) = \big( t\!\!\mod\tau, \;S(t,\om)\big)$ generates  a $\tau$-periodic measure for the skew-product RDS generated by $\tilde\Phi$ on the flat cylinder $\Sc\times\Rd$, $\Sc\simeq \R\Mod\tau$. \\[.2cm]
\underline{Existence of random periodic paths for the stochastic flow $\phi$}. Set $\xi,\eta\in \Rd$ to be  random variables on the filtered probability space $(\Om,\F_{-\infty}^s, \p)$, s.t.~$\xi,\eta \in L^{p}(\Om, \F_{-\infty}^{s}, \p)$. Then, by It\^o formula (e.g., Theorem 4.2.4 in \cite{Kunita} or Theorem 8.1 in \cite{Hkunita})  and Assumption~\ref{A2.1} 
we have for $s\leqslant t$
\begin{align*}
d\,\E\Big[V\big(t,\phi(t,s,\xi)-\phi(t,s,\eta)\big)\Big] &=  \E\Big[ \mathcal{L}^{(2)}V\big(t,\phi(u,s,\xi)-\phi(t,s,\eta)\big)\Big] dt\\[.1cm]
&\hspace{-0cm}\leqslant \lambda(t)\,\E\big[ V\big(t,\phi(t,s,\xi)-\phi(u,s,\eta)\big)\big] dt,
\end{align*}
where $\E\big[ V\big(t,\phi(t,s,\xi)-\phi(u,s,\eta)\big)\big]:=\int_\Omega V\big(t,\phi(t,s,\omega,\xi)-\phi(u,s,\omega,\eta)\big)\mathbb{P}(d\om)$. Thus, by the first part of  (\ref{Lyap_f}) and Gronwall's inequality, we arrive at 
 \begin{align}\label{EA1_m}
\E\vert \phi(t,s,\xi)-\phi(t,s,\eta)\vert^{p}&\leqslant  \E\Big[ V\big(t,\phi(t,s,\xi)-\phi(t,s,\eta)\big)\Big]\notag\\[.1cm]
&\leqslant  \E\Big[ V(s,\xi-\eta)\Big] \exp\left(\int_s^t\lambda(u)du\right).
\end{align}
Finally, given the bound (\ref{EA1_m}), for $r<s<t,$ we have 
\begin{align*}
 \E\big| \phi(t,r,\xi)-\phi(t,s,\xi)\big|^{p}&=\E\big| \phi(t,s,\phi(s,r,\xi))-\phi(t,s,\xi)\big|^{p}\\[.1cm]
 &\leqslant \E\Big[ V\big(s,\phi(s,r,\xi)-\xi\big)\Big]\exp\left(\int_s^t\lambda(u)du\right), 
\end{align*}
 and, utilising the above with Assumption \ref{A2.1}(iii),  yields 
 \begin{equation}
  \limsup_{r<s,\,(t-s)\rightarrow \infty} \E\big| \phi(t,r,\xi)-\phi(t,s,\xi)\big|^{p}= 0 .
 \end{equation}

\noindent Thus, for $\xi \in L^{p}(\Om,\F_{-\infty}^s,\p), \; 1< p<\infty$, the above bound implies that the $L^{p}$ limit of the flow $\{\phi(t,s,\ccdot,\xi)\!:\, s\leqslant t\}$  exists as $s\rightarrow-\infty$. Note that this limit is independent of the initial condition $\xi$ by (\ref{Temp1}). We denote this limit by the random process $S: \R\times\Om\rightarrow \Rd,$ so that  
\begin{align*}
\E\vert S(t)-\phi(t,s,\xi)\vert^{p}\rightarrow 0 \quad \text{as}\;\;s\rightarrow -\infty,
\end{align*}
for $\xi\in L^{p}(\Om,\F_{-\infty}^s,\p)$, where $S(t) := S(t,\ccdot)$. Then, by Chebyshev's first inequality (aka Markov's inequality; e.g.,~\cite{AMB00}), for any $\varepsilon>0,$ we have
\begin{align}\label{cheb}
\p\big(\{\om\in \Om: \vert S(t,\om)-\phi(t,s,\om,\xi)\vert \geqslant \varepsilon\}\big)\leqslant \varepsilon^{-p}\,\E\vert S(t)-\phi(t,s,\xi)\vert^{p},
\end{align}
which  implies that the convergence is also in probability. Thus, there exists a subsequence $(s_k)_{k\in \N_1}$ in $ \R$ with $s_k\rightarrow -\infty$ as $k\rightarrow \infty$ such that 
\begin{align*}
S(t,\om) = \lim_{k\rightarrow\infty}\phi(t, s_k,\om,\xi) \quad \p\;\text{-\,a.s.}
\end{align*}
To simplify notation, we write 
\begin{align}
S(t,\om) = \lim_{s\rightarrow -\infty}\phi(t,s,\om,\xi) \quad \p\;\text{-\,a.s.}
\end{align}
Note that for $\xi\in L^{p}(\Om, \F_{-\infty}^s, \p)$ with the norm $\Vert \cdot \Vert_{p}:= (\E\vert \cdot\vert^{p})^{{1}/{p}}$ we have 
\begin{align*}
\Vert \phi(t,s,\xi)\Vert_{p} &\leqslant \Vert \phi(t,s,\xi)-\xi\Vert_{p}+\Vert \xi\Vert_{p}\\
&\hspace{.0cm} \leqslant\Big(\E\Big[V\big(t,\phi(t,s,\xi)-\xi\big)\Big]\Big)^{\frac{1}{p}}+\Vert \xi\Vert_{p}\\
&\hspace{.0cm} \leqslant\Big(\sup_{s\leqslant t}\E\Big[ V\big(t,\phi(t,s,\xi)-\xi\big)\Big]\Big)^{\frac{1}{p}}+\Vert \xi\Vert_{p} <\infty, 
\end{align*}
by condition (\ref{Temp1}) of Assumption \ref{A2.1}. Consequently, for any $t\in \R$, we have 
\begin{align}
\Vert S(t)\Vert_{p}\leqslant \limsup_{s\rightarrow -\infty}\Vert \phi(t,s,\xi)\Vert_{p}<\infty,
\end{align}
implying  that  $S(t,\om)$ is bounded in $L^{p}(\Om, \F_{-\infty}^{\hspace{.03cm}t}, \p)$.

\smallskip
Next, we show that $t\rightarrow S(t,\om)$ is a random periodic path of period $0<\tau<\infty$ for the stochastic flow~$\{\phi(t,s,\ccdot,\ccdot)\!:\, s\leqslant t\}$ using its  $\tau$-periodic property  (see equation (\ref{Periodic_FL}) with appropriately changed variables); namely
\begin{align}\label{Qr}
\notag S(t+\tau,\om) &= \lim_{s\rightarrow -\infty}\phi(t+\tau,s,\om,\xi)\\
\notag &= \lim_{s\rightarrow -\infty}\phi(t+\tau,s-\tau+\tau,\om,\xi)\\
\notag &= \lim_{s\rightarrow -\infty}\phi(t,s-\tau,\theta_{\tau}\om,\xi)\\
& = S(t, \theta_{\tau}\om) \hspace{.5cm} \p\;\text{-\,a.s.}
\end{align}
Then, by the continuity of $(t,s, x)\mapsto \phi(t,s,\ccdot, x)$ and the flow property, we have 
\begin{align}\label{Qrp}
\notag \phi\big(t+s,s,\om, S(s,\om)\big) &= \lim_{r\rightarrow -\infty}\phi\big(t+s, s, \om, \phi(s, r, \om,\xi)\big)\\
\notag &=\lim_{r\rightarrow -\infty}\phi(t+s,r,\om,\xi) \\ & = S(t+s, \om) \hspace{2cm}  \forall \,t\in \Rp,\, s\in\R \qquad  \p\;\text{-\,a.s.}
\end{align}
The equalities (\ref{Qr}) and (\ref{Qrp}) imply that $S(t,\om)$ is a random periodic path (\ref{per_sol_sde}) of period $0<\tau<\infty$ of the stochastic flow  $\big\{\phi(t+s,s,\ccdot,\ccdot\,)\!:\; s\in \R, \;t\in \Rp\big\}$ on $\Rd.$\\[.2cm]
\underline{Uniqueness}: Let $S_1(t,\om)$ and $S_2(t,\om)$ be two random periodic paths of the forward stochastic flow $\big\{\phi(t+s,s,\ccdot,\ccdot\,)\!:\; s\in \R, \;t\in \Rp\big\}$ on $\Rd.$ We know from (\ref{Qrp}) that for $s, t\in \R$ with $s\leqslant t$,
\begin{align*}
S_1(t,\om) &= \phi\big(t,s,\om, S_1(s,\om)\big) \qquad \p\;\text{-\,a.s.,}\\
S_2(t,\om) &= \phi\big(t,s,\om, S_2(s,\om)\big) \qquad \p\;\text{-\,a.s.}
\end{align*}
Then, for $1<p<\infty,$ we have 
\begin{align*}
\big\Vert S_1(t)- S_2(t)\big\Vert^{p}_{p}&= \big\Vert \phi\big(t,s,S_1(s))-\phi(t,s,S_2(s)\big)\big\Vert^{p}_{p}\\
&\hspace{0cm} \leqslant \exp\left(-\bar{\lambda}(t-s)\right)\E\Big[ V\big(s, S_1(s)-S_2(s)\big)\Big]\underset{s\rightarrow -\infty}{\longrightarrow 0}.
\end{align*}
Thus, $S_1(t,\om) = S_2(t,\om)$ for all $t\in \R \;\; \;\p\,\text{-\,a.s.}$\\[.2cm]
\underline{Construction of $\tau$-periodic measure for  the RDS $\tilde\Phi$}:
Let $\tilde{S}: \R\times\Om\rightarrow \Sc\times\Rd$, $\Sc\simeq \R\Mod\tau$, be defined by 
\begin{align*}
\tilde{S}(r,\om) = \big( r\Mod\tau, \;\,S(r,\om)\big) \qquad \forall \,r\in \Rp,
\end{align*}
or, alternatively $\tilde{S}(r,\om) = \big( |r|\Mod\tau, \;\,S(r,\om)\big), \;\forall \,r\in \R$. Then 
\begin{align}\label{RQA1}
\tilde{S}(r+\tau, \om) =  \big( r+\tau\mmod\tau,\;\,S(r+\tau,\om)\big)= \big( r\mmod\tau, \;\,S(r, \theta_\tau\hspace{.02cm}\om)\big),
\end{align}
and from (\ref{Lift}) and (\ref{Qrp}) we have 
\begin{align}\label{RQA2}
\notag \tilde\Phi(t,\theta_r\hspace{.02cm}\om, \tilde{S}(r,\om)) &= \tilde\Phi\big(t,\theta_r\hspace{.02cm}\om, \big( r\Mod\tau, \;\,S(r, \om)\big)\big)\\
\notag &=\big( t+r\Mod\tau, \;\,\phi(t+r,r,\om, S(r,\om))\big)\\
\notag &= \big( t+r\Mod\tau, \;\,S(t+r,\om)\big)\\
&= \tilde{S}(t+r,\om) \hspace{3cm} \forall \,t,r\in \Rp \quad  \p\;\text{-\,a.s.}
\end{align}
The equalities (\ref{RQA1})--(\ref{RQA2}) and the lifted version of  (\ref{per_sol_rds}), imply that $\tilde{S}(t,\om)$ is a random $\tau$-periodic path  of the skew-product RDS generated by $\tilde\Phi$ (\ref{Lift}) on the flat cylinder $\Sc\times\Rd.$

Finally, let $(\tilde{\mu}_t)_{t\in \Rp}$, $\tilde\mu_t\in \PP\big(\Sc\big)\otimes\PP\big(\Rd\big)$ be defined by 
\begin{align*}
\tilde{\mu}_t(\tilde{A}) = \p\big(\{\om: \tilde{S}(t,\om)\in \tilde{A}\,\}\big) \quad \forall\;t\in \Rp, \; \tilde{A}\in \Bb\big([0,\tau)\big)\times\Bb\big(\Rd\big).
\end{align*}
It follows from  (\ref{RQA1})\,-\,(\ref{RQA2}) and Lemma \ref{semi_P} that the probability measure $\tilde{\mu}_t$ is  $\tau$-periodic  under the action of the transition  evolution $(\tilde{\mathcal{P}}^*_{t})_{t\in \Rp}$ which  is  induced by the RDS $\{\tilde\Phi(t,\ccdot,\ccdot)\!:\; t\in \Rp\}$ on $\Sc\times\Rd$. The skew-product structure of these measures in $\PP\big(\Sc\big)\otimes\PP\big(\Rd\big)$ arises from Lemma~\ref{semi_P}, or directly from (\ref{PPhi}), so that  for any $\mathcal{J}\in \Bb\big(\Sc\big), \; A\in \Bb\big(\Rd\big)$
\begin{align*}
\hspace{1.5cm}\tilde{\mu}_t(\Jc\times A) &= \delta_{(t \mmod\tau)}(\mathcal{J})\otimes\p\big(\om: S(t,\om)\in A\big)=\delta_{(t\mmod\tau)}(\mathcal{J})\otimes\mu_t(A). \hspace{2cm}\qed 
\end{align*}

\subsection{Ergodicity of time-periodic measures}\label{permes_erg}
In this section, we turn to establishing  ergodicity of the  $\tau$-periodic measures $(\tilde{\mu}_t)_{t\in \Rp}$, $\tilde\mu_t\in \PP\big(\Sc\big)\otimes\PP\big(\Rd\big)$,  
 generated by the Markovian\footnote{\,Here, the notion of a {\it `Markovian RDS'} means that there exists a version of the RDS which has the Markov property w.r.t~the filtration generated on the Wiener space by the canonical DS for the Wiener process with $\tilde W_{t+s}(\theta_{-s}\om)$ for all $s\in \Sc$, $t\in \Rp$; see Proposition~\ref{canW}.} RDS $\big\{\tilde\Phi(t,\ccdot,\ccdot\,):\,t\in \Rp\big\}$ which was constructed in~(\ref{Lift}) 
in the skew-product representation on  the flat cylinder $\Sc\times\Rd$ from the lifted flow of solutions of the SDE (\ref{NSDE}) with time-periodic coefficients. The existence of  $\tau$-periodic measures supported on stable random periodic paths was established in Theorem \ref{Rand_SOL}. The lack of stationarity and the unavoidable skew-product structure of the underlying dynamics pose additional challenges when dealing with ergodicity of  $\tilde{\mathcal{P}}^*_t$-\,invariant measures, as outlined below. The main theorem of this section (Theorem~\ref{Ps_erg}) is preceded by some preparatory results and definitions. 

\begin{definition}[{\bf Ergodic periodic measure} \cite{Feng18}]\label{PS-eR}\;\rm
 A family of $\tau$-periodic measures $(\tilde{\mu}_t)_{t\in\Rp}$ on the extended state space $\big(\Sc\times\Rd,\,\Bb(\Sc)\times\Bb(\Rd)\big)$ is said to be {\it ergodic} if 
 \begin{align}\label{ERm}
\bar{\tilde{\mu}} = \frac{1}{\tau}\int_0^\tau \tilde{\mu}_t \hspace{0.03cm}dt,
\end{align} 
 is ergodic with respect to the transition semigroup $(\tilde{\mathcal{P}}^*_t)_{t\in \Rp}$ in (\ref{tldP*}).
 \end{definition}

One can check by the linearity of $\tilde\mu_0\mapsto\tilde{\mathcal{P}}^*_{t}\tilde\mu_0$  and Fubini's theorem that $\bar{\tilde\mu}$ is an invariant measure  for the transition semigroup $(\tilde{\mathcal{P}}^*_t)_{t\in \Rp}$ defined in (\ref{tldP*}); i.e., $\tilde{\mathcal{P}}^*_t$-\,invariance of $\bar{\tilde\mu}$ implies $\tilde{\mathcal{P}}^*_t\bar{\tilde\mu}=\bar{\tilde\mu}$, for all $t\in \Rp$.  Moreover, from the definition of a $\tau$-periodic measure $\tilde\mu_t$ in ~(\ref{permes_rds}), induced by the RDS $\{\tilde\Phi(t,\ccdot,\,\cdot\,): \,t\in \Rp\}$ on $\Sc\times\Rd$, we have for any $\tilde{A}\in \Bb\big(\Sc\big)\times\Bb\big(\Rd\big)$ that 
 \begin{align*}
 \bar{\tilde{\mu}}(\tilde{A})  = \frac{1}{\tau}\int_0^\tau \tilde{\mu}_t(\tilde{A})dt &= \frac{1}{\tau}\int_0^\tau\p\big(\big\{\om: \tilde{S}(t,\om)\in \tilde{A}\big\}\big)dt=\frac{1}{\tau}\E\left[\int_0^\tau\I_{\tilde{A}}\big(\tilde{S}(t, \ccdot)\big)dt\right]\\[.1cm]
&=\E\left[\frac{1}{\tau}m_1\big(\big\{t\in [0, \tau): \tilde{S}(t,\ccdot)\in \tilde{A}\big\}\big)\right],
 \end{align*}
where  $t\rightarrow \tilde S(t,\om)=\big(t\Mod\tau,\,S(t,\om)\big)$, $t\in \Rp$, is a random periodic path  (\ref{per_sol_rds}) of an RDS generated by the lifted dynamics of the SDE (\ref{NSDE}) via $\tilde\Phi$ in (\ref{Lift}),
and $m_1$ is the Lebesgue measure on $\R$. Thus, given the invariance of  $\bar{\tilde{\mu}}$ under the action of the transition semigroup $(\tilde{\mathcal{P}}^*_t)_{t\in \Rp}$~in~(\ref{tldP*}), and the $\tau$-periodicity  of $\tilde \mu_t$ (see Definition~\ref{Rand_pm}), one has 
 \begin{align*}
 \E\left[\frac{1}{\tau}m_1\big(\{t\in [0, \tau): \tilde{S}(t,\ccdot)\in \tilde{A}\}\big)\right]&= \big(\tilde{\mathcal{P}}_u^{*}\bar{\tilde{\mu}}\big)(\tilde{A})=\frac{1}{\tau}\int_0^\tau(\tilde{\mathcal{P}}^{*}_u\tilde{\mu}_t)(\tilde{A}) dt\\[.0cm]
 &=\frac{1}{\tau}\int_0^\tau \tilde{\mu}_{t+u}(\tilde{A}) dt=\frac{1}{\tau}\int_u^{u+\tau}\tilde{\mu}_t(\tilde{A}) dt\\[.1cm]
 &=\E\left[\frac{1}{\tau}m_1\big(\{t\in [u, u+\tau): \tilde{S}(t,\ccdot)\in \tilde{A}\}\big)\right],
 \end{align*}
for any $\tilde{A}\in \Bb\big(\Sc\big)\times\Bb\big(\Rd\big)$ and any $u\in \Rp$. This implies that the expected time spent by the random periodic path $t\mapsto\tilde{S}(t,\om)$  in  any set $\tilde{A}\in\Bb\big(\Sc\big)\times\Bb\big(\Rd\big)$ over a time interval of exactly one period is independent of the starting point. 
 
\medskip
Verification of ergodicity (in the sense of  Definition~\ref{PS-eR}) of $\tau$-periodic measures $(\tilde{\mu}_t)_{t\in \Rp}$ supported on the random periodic paths of $\tilde\Phi$ requires one to assert  that the time-averaged skew-product measure $\bar{\tilde{\mu}}$ in the fibre bundle  $\PP\big(\Sc\big)\otimes\PP\big(\Rd\big)$  is $\tilde{\mathcal{P}}^*_t$-\,ergodic.  This setup arises from the need to deal with the random periodic, skew-product nature of the underlying dynamics, and it prevents a direct application of the classical tools for asserting ergodicity in the (asymptotically) stationary case.   In particular, it is well-known (e.g., \cite[Theorem 3.2.4]{Da Prato}) that the following are equivalent\footnote{\,These statements are not restricted to the skew-product representation of time-periodic measures.}:
\begin{itemize}[leftmargin=.8cm]
\item[(i)] A probability measure $\bar{\tilde \mu}$ is weakly mixing. 

\vspace{.2cm}\item[(ii)] There exists $\Ic\subset[0, \;\infty)$ of relative measure 1 such that $\lim_{t\rightarrow \infty, t\in \Ic} \tilde P(t,\tilde x, \ccdot)\rightarrow \bar{\tilde \mu}$ weakly.
\end{itemize} 
Thus, given the form of the transition kernel $\tilde P$ in (\ref{PPhi}) and the underlying skew-product structure, it is clear that one cannot establish the mixing property in the random periodic regime\footnote{\,As before, we exclude the stationary regime from the random periodic regime by requiring that fundamental period $0<\tau<\infty$; see Definition \ref{Rand_pm}.}. Thus, this key condition in Doob's Theorem \cite{Doob48} does not hold in the random periodic regime which, alongside the lack of  irreducibility of the transition kernel, renders the Hasminskii's Theorem \cite{Has12}  for asserting regularity of the transition kernel (needed in Doob's Theorem)~inapplicable.

Instead, the $\tilde{\mathcal{P}}^*_t$-\,ergodicity of $\bar{\tilde\mu}$  can be verified by means of a proposition which was proved in \cite[Lemma~2.18]{Feng18}; we repeat its statement below with a  concise proof to make this section self-contained. The main benefit of utilising the proposition below when dealing with $\bar{\tilde\mu}$ is that it essentially relies on ergodicity of $\tau$-periodic measures $\tilde\mu_t$ for any fixed $t\in \Sc$ with respect to the discrete dynamics induced by $(\tilde{\mathcal{P}}^*_{n\tau})_{n\in \mathbb{N}_0}$;  the subsequent use of the semigroup property of $(\tilde{\mathcal{P}}^*_{t})_{t\in \Rp}$ allows one to show the ergodicity of $\bar{\tilde\mu}$. Importantly,  the $\tilde{\mathcal{P}}^*_{n\tau}$\,-\,ergodicity of $\tilde\mu_t$ on the respective Poincar\'e sections  with a fixed $t\in \Sc$ turns the problem into a stationary one which  can be dealt with using the standard methods. The result below provides an extension of the classical Krylov-Bogolyubov procedure (see, e.g., \cite[\S 1.5]{Arnold}).

\begin{prop}\label{Ergequi}
 Consider a  family of $\tau$-periodic measures $(\tilde{\mu}_t)_{t\in\Rp}$ on the extended state space $\big(\Sc\times\Rd,\Bb(\Sc)\times\Bb(\Rd)\big)$. The $\tilde{\mathcal{P}}_t^*$-\,invariant measure $\bar{\tilde{\mu}}$ in (\ref{ERm}) is ergodic if and only if the following holds for any $\tilde{A}\in \Bb([0,\tau))\times\Bb(\Rd)$
 \begin{align}\label{Krylov}
\lim_{N\rightarrow\infty}\int_{\Sc\times\Rd}\left\vert \int_{0}^{\tau}\bigg\{\frac{1}{N}\sum_{n=0}^{N-1}\tilde{P}(\tilde{x};t+n\tau, \tilde{A}) - \tilde{\mu}_t(\tilde{A})\bigg\}dt\right\vert\bar{\tilde{\mu}}(d\tilde{x}) =0.
\end{align}
\end{prop}
\smallskip
\noindent {\it Proof.} Recall from (e.g.,~\cite{Arnold}) that $\bar{\tilde{\mu}}$ is ergodic if $\tilde{\mathcal{P}}_t\,\I_{\tilde{A}} = \I_{\tilde{A}}, \; \bar{\tilde{\mu}}\,\text{-\,a.e.}$  $\tilde{A}\in \Bb([0, \tau))\times\Bb(\Rd)$ implies  that either $\bar{\tilde{\mu}}(\tilde{A}) = 0$ or $\bar{\tilde{\mu}}(\tilde A) =1.$ First, we assume that (\ref{Krylov}) holds for any $\tilde{A}\in \Bb([0,\tau))\times\Bb(\Rd)$ with $\tilde{P}(\tilde{x};t,\tilde{A}) = \tilde{\mathcal{P}}_t\,\I_{\tilde{A}}(\tilde{x}) = \I_{\tilde{A}}(\tilde{x})$.  Then, it follows from (\ref{Krylov}) that 
\begin{align*}
\int_{\Sc\times\Rd}\Big\vert \I_{\tilde{A}}(\tilde{x})-\bar{\tilde{\mu}}(\tilde{A})\Big\vert \bar{\tilde\mu}(d\tilde{x}) = \int_{\Sc\times\Rd}\left\vert \frac{1}{\tau}\int_{n\tau}^{(n+1)\tau}\tilde{P}(\tilde{x};t,\tilde{A}) dt - \bar{\tilde{\mu}}(\tilde{A})\right\vert \bar{\tilde{\mu}}(d\tilde{x}) =0.
\end{align*}
This implies that $\I_{\tilde{A}}(\tilde{x})$ is a constant for $\bar{\tilde{\mu}}\;\text{-\,a.e.}\; \tilde{x}\in \Sc\times\Rd.$ Thus, either $\bar{\tilde{\mu}}(\tilde{A})=0$ or $\bar{\tilde{\mu}}(\tilde{A}) =1$.  Conversely, assume that $\bar{\tilde{\mu}}$ is ergodic, then for any $\tilde{A}\in \Bb([0,\tau))\times\Bb(\Rd)$
\begin{align*}
\lim_{T\rightarrow\infty}\frac{1}{T}\int_0^T\tilde{P}(\tilde{x};t,\tilde{A})dt = \bar{\tilde{\mu}}(\tilde{A}) \quad \text{in}\; L^2(\Rd;\bar{\tilde \mu}).
\end{align*} 
Therefore,
\vspace{-.3cm}
\begin{align}\label{IIF}
\lim_{N\rightarrow\infty}\frac{1}{N\tau}\sum_{n=0}^{N-1}\int_0^\tau \tilde{P}(\tilde{x};t+n\tau,\tilde{A})dt = \bar{\tilde{\mu}}(\tilde{A}) \quad \text{in}\; L^2(\Rd;\bar{\tilde{\mu}}),
\end{align}
and  (\ref{Krylov}) follows from (\ref{IIF}) and from the Cauchy--Schwartz inequality. \qed

\medskip
Consequently, the subsequent  verification of  the ergodicity of the $\tilde{\mathcal{P}}_t^*$-\,invariant measure $\bar{\tilde{\mu}}$  on $\mathcal{B}(\Sc)\times\mathcal{B}(\Rd)$, 
relies (explicitly or otherwise) on the semigroup property and periodicity of the transition semigroup $(\tilde{\mathcal{P}}_t^*)_{t\in \Rp}$, and on  proving the strong Feller property of the transition evolution $(\mathcal{P}_{s,t})_{t\geqslant s}$ in~(\ref{calP}). Recall that the transition evolution $(\Pt_{s,t})_{t\geqslant s}$ has the strong Feller property (i.e., $\Pt_{s,t}\varphi \in \mathcal{C}_\infty(\Rd)$ for any $\varphi\in \mathbb{M}_\infty(\Rd)$) if and only if
\begin{itemize}[leftmargin =0.9cm]
\vspace{.2cm}\item[(i)] $(\Pt_{s,t})_{t\geqslant s}$ is Feller; i.e., $\Pt_{s,t}: \mathcal{C}_\infty(\Rd)\rightarrow\mathcal{C}_\infty(\Rd)$. 

\vspace{.2cm}\item[(ii)] For any $\varphi\in \mathcal{C}_\infty(\Rd)$ the family $(\Pt_{s,t}\varphi)_{t\geqslant s}$ is equicontinuous.

\end{itemize}  

\smallskip
\noindent The first condition follows from the existence of the stochastic flow (see, e.g.,~\cite{Kunita, Has12});  thus, we only derive the second item in Proposition \ref{StrongF} below.

\begin{prop}\label{StrongF} Suppose that Assumption \ref{Homa} holds. Then, for any $t\in [s, s+T),$ there exist $0<C_T<
  \infty$ such that, for any $x,y\in \Rd$ and any $\varphi\in \mathcal{C}_\infty(\Rd),$ we have
  \begin{align*}
  \vert \mathcal{P}_{s,t}\varphi(x)-\mathcal{P}_{s,t}\varphi(y)\vert \leqslant C_T\Vert \varphi\Vert_{\infty}\vert x-y\vert.
\end{align*}   
\end{prop}
\noindent {\it Proof.} The proof consists of a tedious but relatively straightforward extension of results which are well known in the autonomous case;  for detailed derivations, involving some Malliavin calculus estimates;  see  Theorem \ref{AStrongF} in Appendix~\ref{SfA}.

\bigskip
Given the above setting, we have the following main result of this section:

\begin{theorem}\label{Ps_erg}
Suppose that Proposition \ref{StrongF} and Assumption \ref{A2.1} hold.~Then, the family of $\tau$-periodic measures $(\tilde{\mu}_t)_{t\in \Rp}$, \,$\tilde\mu_t\in \PP\big(\Sc\big)\otimes\PP\big(\Rd\big)$, in (\ref{permes_rds}) is~ergodic in the sense of Definition~\ref{PS-eR}.
\end{theorem}

\begin{rem}\label{disc_ass}\rm The requirement in the above theorem that Assumption \ref{A2.1} holds is inherited from the conditions required in Theorem \ref{Rand_SOL} for the existence of stable random periodic paths on which the $\tau$-periodic skew-product measures $(\tilde{\mu}_t)_{t\in \Rp}$ are supported; hence, the only additional condition in Theorem \ref{Ps_erg} is introduced by imposing  Assumption~\ref{Homa} which is required in Proposition \ref{StrongF} to assert the strong Feller property of $(\mathcal{P}_{s,t})_{t\geqslant s}$. If one dropped Assumption \ref{A2.1}, the existence of $\tau$-periodic skew-product probability measures would have to be assumed a priori alongside the ergodicity of $\tilde\mu_t$ for all fixed $t\in \Sc$ w.r.t.~the discrete transition evolution $(\tilde{\mathcal{P}}_{n\tau})_{n\in \mathbb{N}_0}$, as done in~\cite{Feng18}.   In the present case, the properties of the $\tau$-periodic measures derived explicitly in the previous section allow us to dispense with such assumptions. 
\end{rem}

\smallskip
\noindent {\it Proof of Theorem \ref{Ps_erg}.} The proof is relatively long and we divide it into four steps. \\Throughout, we skip the dependence on $\om\in \Om$ in all quantities involving expectations w.r.t.~$\p$.

\medskip
\noindent \underline{\underline{Step I}}: First, we show that for a random periodic path $S: \R\times\Om\rightarrow \Rd$ of the stochastic flow $\phi$ on $\Rd$, and  $\eta\in L^{p}(\Om, \F_{-\infty}^s, \p)$, $1< p<\infty$, there exists $0<\tilde{\mathfrak{C}}<\infty$ such that 
\vspace{.2cm}
\begin{align}\label{ExpPER}
\Vert \phi(s+n\tau, s,\eta)-S(s+n\tau)\Vert_{p}\leqslant \tilde{\mathfrak{C}}\exp\left(\frac{1}{{p}}\int_s^{s+n\tau}\lambda(u)du\right), \qquad n\in \N_0. 
\end{align}

\vspace{.2cm}
\noindent To see this, note that from the definition of the random periodic path of a stochastic flow (\ref{per_sol_sde}) we have  $S(s+n\tau,\om) = \phi(s+n\tau,s,\om, S(s,\om)) \; \,\p\,\text{-a.s.},$ so that

\newpage
\begin{align}\label{D2E3}
\notag  \Vert \phi(s+n\tau, s,\eta)-S(s+n\tau)\Vert_{p}  & =  \Vert \phi(s+n\tau, s,\eta)-\phi(s+n\tau, s,S(s))\Vert_{p} \\[.2cm]
&\hspace{-.0cm} \leqslant \Big(\E\Big[ V(s,\eta-S(s)) \Big]\Big)^{\frac{1}{p}}\exp\left(\frac{1}{p}\int_s^{s+n\tau}\lambda(u)du\right)\notag\\[.2cm]
& \hspace{-0.cm}= \tilde{\mathfrak{C}}\,\exp\left(\frac{1}{p}\int_s^{s+n\tau}\lambda(u)du\right) , \quad \qquad n\in \N_0,
\end{align}
by  Assumption \ref{A2.1}(i) and the fact that $S(s)\in L^{p}(\Om, \F_{-\infty}^s, \p)$, $1< p<\infty$, which was shown in the proof of Theorem \ref{Rand_SOL}.

\medskip
\noindent \underline{\underline{Step II}}: We show that for $1< p<\infty$ there exists $0<\mathfrak{C}_{\,\tau}<\infty$ such that for $n\in \N_0$
\begin{align}\label{D2E4}
\bigg\vert \Pt_{s,s+n\tau}\varphi(x)-\int_{\Rd} \varphi(y)\mu_s(dy)\bigg\vert \leqslant \mathfrak{C}_{\,\tau}\Vert \varphi\Vert_\infty\exp\left(\frac{1}{p}\int_s^{s+n\tau}\hspace{-.2cm}\lambda(u)du\right) \qquad  \forall\;\varphi\in \mathcal{C}_\infty(\Rd),
\end{align}
where $\mu_s(A)= \p\big(\{\om: S(s,\om)\in A \}\big), \; A\in \Bb(\Rd).$

\smallskip
To see this, we note that from the definition of the periodic measure $\mu_s,$ we have that 
$$ \int_{\Rd}\Pt_{s,s+n\tau}\varphi(y)\mu_s(dy) = \int_{\Rd}\varphi(y)\mu_s(dy) \qquad \forall\;\varphi\in \mathcal{C}_\infty(\Rd);$$ 
i.e., $\mu_s$ is invariant under the action of the dual of the discrete transition evolution $(\Pt^*_{s,s+n\tau})_{n\in \N_0}$. Thus, for $\psi\in \text{Lip}_\infty(\Rd),$ we have for $1< p <\infty,$
\begin{align}\label{LipCH}
\notag \bigg\vert \Pt_{s,s+n\tau}\psi(x)-\int_{\Rd} \psi(y)\mu_s(dy)\bigg\vert &= \bigg\vert \int_{\Rd}\Big( \Pt_{s,s+n\tau}\psi(x)-\Pt_{s,s+n\tau}\psi(y)\Big)\mu_s(dy)\bigg\vert\\[.1cm]
\notag &\hspace{.0cm} \leqslant \Vert \psi\Vert_{\textsc{bl}}\int_{\Rd}\E\big| \phi(s+n\tau,s,x)-\phi(s+n\tau,s,y)\big|\mu_s(dy)\\[.1cm]
\notag &= \Vert \psi\Vert_{\textsc{bl}}\,\E\big| \phi(s+n\tau,s,x)-\phi(s+n\tau,s,S(s))\big|\\[.1cm]
\notag & \leqslant \Vert \psi\Vert_{\textsc{bl}}\Big(\E\vert \phi(s+n\tau, s,x)-S(s+n\tau)\vert^{p}\Big)^{\frac{1}{p}}\\[.1cm]
&\hspace{0cm} \leqslant \tilde{\mathfrak{C}}\,\Vert\psi\Vert_{\textsc{bl}}\exp\left(\frac{1}{p}\int_s^{s+n\tau}\lambda(u)du\right),
\end{align}
where we applied H\"older's inequality and estimate (\ref{ExpPER}) in the last two lines respectively. 

Now, let $\varphi\in \mathcal{C}_\infty(\Rd)$ be given. Setting $\psi= \Pt_{s+n\tau,\,s+\tau+n\tau}\hspace{.02cm}\varphi = \Pt_{s,s+\tau}\hspace{.02cm}\varphi $ in (\ref{LipCH}), which holds due to (\ref{Periodic_FL}), and  using the invariance of $\mu_s$ under the transition evolution $(\Pt^*_{s,s+n\tau})_{n\in \N_0}$, we obtain by   Proposition \ref{StrongF} that 
\begin{align}\label{CH}
\bigg\vert \Pt_{s,s+\tau+n\tau}\varphi(x)-\int_{\Rd}\Pt_{s,s+\tau}\varphi(y)\mu_s(dy)\bigg\vert&=\bigg\vert \int_{\Rd}\Big(\Pt_{s,\,s+\tau+n\tau}\varphi(x)-\Pt_{s,\,s+\tau+n\tau}\varphi(y)\Big)\mu_s(dy)\bigg\vert \notag\\[.1cm]
&\hspace{.0cm}\leqslant \tilde{\mathfrak{C}}\,\Vert \Pt_{s+n\tau,\,s+\tau+n\tau}\varphi\Vert_{\textsc{bl}}\exp\left(\frac{1}{p}\int_s^{s+n\tau}\lambda(u)du\right)\notag\\[.1cm]
&=\tilde{\mathfrak{C}}\,\Vert \Pt_{s,s+\tau}\varphi\Vert_{\textsc{bl}}\exp\left(\frac{1}{p}\int_s^{s+n\tau}\lambda(u)du\right)\notag\\[.1cm]
&\hspace{0cm}\leqslant \mathfrak{C}_{\,\tau}\Vert \varphi\Vert_{\infty}\exp\left(\frac{1}{p}\int_s^{s+n\tau}\lambda(u)du\right),
\end{align}
where $\mathfrak{C}_{\,\tau} = C_\tau\tilde{\mathfrak{C}},$ and  $C_\tau$ is a constant appearing in Proposition \ref{StrongF}. 

\smallskip
\noindent \underline{\underline{Step III}}: Let $A\subset \Rd$ be a closed set, take $\varphi = \I_{A}$, and consider the sequence $(\varphi_m)_{m\in \N_1}$ of functions  defined by 
\vspace{-.2cm}\begin{align*}
\varphi_m(x) = \begin{cases} 1, \quad & \text{if}\;\; x\in A,\\
1-2^m \text{d}(x,A), & \text{if}\;\;\text{d}(x,A)\leqslant 2^{-m},\\
0, &\text{if}\;\; \text{d}(x,A)\geqslant 2^{-m},
\end{cases}
\end{align*}
where $\text{d}(x,A) = \inf\{\vert x-y\vert: y\in A\}, \; x\in \Rd.$ Then
$$\varphi_m(x)\rightarrow \varphi(x), \quad \text{as}\; m\rightarrow\infty\quad \forall \;x\in \Rd.$$
Next, for $s\in [0, \tau),$ we have 
\begin{align*}
\Pt_{s,s+n\tau}\,\varphi_m(x)\rightarrow \Pt_{s,s+n\tau}\,\varphi(x)= \Pt_{s,s+n\tau}\,\I_A(x),
\end{align*}
which  implies that $P(s,\ccdot\,;s+n\tau,A)= \Pt_{s,s+n\tau}\,\I_A\in \mathcal{C}_\infty(\Rd)$ and, since  $\mu_s$ is invariant under $(\Pt^*_{s,s+n\tau})_{n\in \N_0}$, (\ref{CH}) leads to 
\begin{align}\label{LIpTR}
 \big\vert P(s,x;s+n\tau,A)-\mu_s(A)\big\vert\leqslant \mathfrak{C}_\tau\exp\left(\frac{1}{p}\int_s^{s+n\tau}\lambda(u)du\right).
\end{align}
By the covering lemma (e.g.,~\cite{AMB00}), the inequality (\ref{LIpTR}) holds for any $A\in \Bb(\Rd)$, and  thus for $\mathcal{J}\subseteq \Sc,$ we have 
\begin{align*}
\int_{\mathcal{J}}\big\vert P(s,x;s+n\tau,A)-\mu_s(A)\big\vert ds &\leqslant \int_0^\tau \big\vert P(s,x;s+n\tau,A)-\mu_s(A)\big\vert ds\\[.1cm]
&\hspace{.0cm} \leqslant \mathfrak{C}_\tau\int_0^\tau \exp\left(\frac{1}{p}\int_s^{s+n\tau}\lambda(u)du\right) ds\\[.1cm]
&\hspace{0cm} = \mathfrak{C}_\tau\int_0^\tau \exp\left(\frac{1}{pn\tau}\int_s^{s+n\tau}\lambda(u)du\right)^{\!\!n\tau} ds.
\end{align*}
Now, we use the Chapmann--Kolmogorov equation (\ref{ChK}) for the transition probability to obtain
\begin{align*}
\bigg\vert \int_{\Jc}\Big[ P(s, x; t+n\tau,A)-\mu_t(A)\Big]dt\bigg\vert &= \bigg\vert \bigg[\int_{\Jc}\int_{\Rd}P(t,y;t+n\tau,A)-\mu_t(A)\bigg] P(s,x;t,dy)dt\bigg\vert\\[.1cm]
&\hspace{.0cm} \leqslant\int_0^\tau\int_{\Rd}\mathfrak{C}_\tau \exp\left(\frac{1}{pn\tau}\int_t^{t+n\tau}\lambda(u)du\right)^{\!\!n\tau}\!\!P(s,x;t,dy)dt\\[.2cm]
&\hspace{0cm} = \mathfrak{C}_\tau\int_0^\tau\exp\left(\frac{1}{pn\tau}\int_t^{t+n\tau}\lambda(u)du\right)^{\!\!n\tau} dt.
\end{align*}
By condition (\ref{Rate1}) of Assumption \ref{A2.1}, there exists $0<\beta<1, \; 0<K<\infty,$ such that 
\begin{align*}
\bigg\vert \int_{\Jc}\Big( P(s,x; t+n\tau,A)-\mu_t(A)\Big) dt\bigg\vert\leqslant \int_{\Jc}\big\vert P(s,x;t+n\tau,A) -\mu_t(A)\big\vert dt\leqslant K \beta^{n\tau}.
\end{align*}
It then follows that 
\begin{align}\label{Kr_conv}
\frac{1}{\tau}\int_0^\tau\int_{\Rd}\bigg\vert \int_{\Jc}\bigg\{\frac{1}{N}\sum_{n=0}^{N-1}P(s,x;t+n\tau,A)-\mu_t(A)\bigg\} dt \bigg\vert \mu_s(dx)ds\leqslant \frac{K}{N}\sum_{n=0}^{N-1}\beta^{n\tau}\underset{N\rightarrow \infty}{\longrightarrow} 0.
\end{align}

\vspace*{0.5cm}
\noindent \underline{\underline{Step IV}}:
In this final step, with the help of Step III, we show the convergence of Krylov-Bogolyubov scheme for the $\tau$-periodic skew-product measures $(\tilde{\mu}_t)_{t\in \Rp}$ on the cylinder $\Sc\times\Rd$. For any $\mathcal{J}\times A\in \Bb([0,\tau))\times\Bb(\Rd)$ we have 
\begin{align*}
&\hspace{-.2cm}\int_{[0,\,\tau)\times\Rd}\bigg\vert\int_0^\tau\bigg(\frac{1}{N}\sum_{n=0}^{N-1}\tilde{P}\big(\tilde x; t+n\tau, \Jc\times A\big)-\tilde{\mu}_t(\Jc\times A)\bigg)dt\bigg\vert\bar{\tilde{\mu}}(d\tilde x) \\
&\hspace{1cm} = \frac{1}{\tau}\int_0^\tau\int_{\Rd}\bigg\vert \int_0^\tau\bigg(\frac{1}{N}\sum_{n=0}^{N-1}\tilde{P}((s,x);t+n\tau,\Jc\times A)-\tilde{\mu}_{t}(\Jc\times A)\bigg)dt\bigg\vert \mu_s(dx)ds\\
&\hspace{1cm}=\frac{1}{\tau}\int_0^\tau \int_{\Rd}\bigg\vert \int_0^\tau\bigg(\frac{1}{N}\sum_{n=0}^{N-1}P(s,x;t+s+n\tau,A)-\mu_{t}(A)\bigg)\delta_{(t+s\!\!\!\!\mod\tau)}(\Jc)dt\bigg\vert \mu_s(dx)ds\\
&\hspace{1cm} = \frac{1}{\tau}\int_0^\tau \int_{\Rd}\bigg\vert \int_0^{\tau-s}\bigg(\frac{1}{N}\sum_{n=0}^{N-1}P(s,x;t+s+n\tau,A)-\mu_{t}(A)\bigg)\delta_{(t+s)}(\Jc)dt\\
&\hspace{2cm}+\int_{\tau-s}^\tau\bigg( \frac{1}{N} \sum_{n=0}^{N-1}P(s,x;t+s+n\tau,A)-\mu_t(A)\bigg) \delta_{(t+s-\tau)}(\Jc)dt\bigg\vert \mu_s(dx)ds\\
&\hspace{1.cm} = \frac{1}{\tau}\int_0^\tau \int_{\Rd}\bigg\vert \int_0^{\tau-s}\bigg(\frac{1}{N}\sum_{n=0}^{N-1}P(s,x;t+s+n\tau,A)-\mu_{t}(A)\bigg) \delta_{(t+s)}(\Jc)dt\\
&\hspace{2cm}+\int_{-s}^0\bigg( \frac{1}{N} \sum_{n=1}^{N}P(s,x;t+s+n\tau,A)-\mu_t(A)\bigg)\delta_{(t+s)}(\Jc) dt\bigg\vert \mu_s(dx)ds\\
&\hspace{1.cm} = \frac{1}{\tau}\int_0^\tau \int_{\Rd}\bigg\vert \int_{\Jc}\bigg(\frac{1}{N}\sum_{n=0}^{N-1}P(s,x;t+n\tau,A)-\mu_{t}(A)\bigg)dt\\
&\hspace{2cm} -\frac{1}{N}\int_{-s}^0\bigg( P(s,x;t+s,A)-P(s,x;t+s+N\tau,A)\bigg)\delta_{(t+s)}(\Jc) dt\bigg\vert \mu_s(dx)ds\\
&\hspace{1.cm}\leqslant \frac{1}{\tau}\int_0^\tau \int_{\Rd}\bigg\vert \int_{\Jc}\bigg(\frac{1}{N}\sum_{n=0}^{N-1}P(s,x;t+n\tau,A)-\mu_{t}(A)\bigg) dt\bigg\vert \mu_s(dx)ds\\[.2cm]
&\hspace{2cm} + \frac{1}{N\tau}\int_0^\tau \int_{\Rd}\bigg\vert \int_{-s}^0\bigg( P(s,x;t+s,A)\\[.2cm] 
&\hspace{5.3cm}-P(s,x;t+s+N\tau,A)\bigg)\delta_{(t+s)}(\Jc) dt \bigg\vert \mu_s(dx)ds\underset{N\rightarrow \infty}{\longrightarrow} 0. \hspace{.5cm}\qed
\end{align*}

\begin{rem}\label{tau_ergavg}\rm 
The invariance of the $\tau$-periodic probability measures  under the discrete evolution $(\tilde{\mathcal{P}}^*_{n\tau})_{n\in \mathbb{N}_0}$ on their respective Poincar\'e sections was pointed out in Lemma \ref{semi_P}. It can be shown, as a consequence of \cite[Theorem 4.11]{Feng18}, that such  $\tau$-periodic probability measures are ergodic w.r.t.~the discrete evolution $(\tilde{\mathcal{P}}^*_{n\tau})_{n\in \mathbb{N}_0}$ on their respective Poincar\'e sections; given that we require Assumption \ref{Homa} to be satisfied, these measures are supported on all of $\Rd$. This fact will be useful in \S\ref{Linear response_per} concerned with ergodic averages in the context of the linear response.
\end{rem}

\begin{exa}[Stochastic Lorenz model with periodic forcing]\label{Lor_ex}\rm
Consider a modified  Lorenz system (e.g.,~\cite{Keller}) given by 
\begin{align}\label{Lorenz1}
\begin{cases}
\dot{x} = -\bar\alpha \hspace{.04cm}x+\bar\alpha \hspace{.04cm} y,\\
\dot{y} = -\bar\alpha \hspace{.04cm} x -\bar\beta \hspace{.04cm} y-xz,\\
\dot{z} = -\bar\gamma \hspace{.04cm}z+xy-\bar\gamma \bar\beta^{-2}\hspace{.04cm}(\bar\varrho+\bar\alpha),
\end{cases}
\end{align}
with parameters $ \bar\alpha, \bar\beta,\bar\gamma,  \bar\varrho>0$. We set  $(v_1,v_2,v_3) := (x,y,z)\in\R^3$  and consider the periodically-in-time and stochastically  perturbed version of (\ref{Lorenz1})  for $t\in \Rp$ in the form
\begin{align}\label{Lorenz}
dv_t = b(t,v_t)dt+\sigma(v_t)dW_t = \big[-Av_t-G(v_t)+F(t)\big]dt+\sigma(v_t)dW_t, \qquad v_0\in \R^3,
\end{align}
where
\begin{align*} 
A&= \begin{bmatrix}
\bar\alpha& -\bar\alpha& 0\\
\bar\alpha & \bar\beta& 0\\
0 & 0 & \bar\gamma
\end{bmatrix}, \;\; G(v)= \begin{bmatrix}
\;\;\;0\\ \;\;\,v_1v_3\\ -v_1v_2
\end{bmatrix}, \;\; F(t) =\begin{bmatrix}
\bar f\big(1+\bar\delta\sin\big(\frac{2\pi}{\tau}\, t\big)\big)\\ 0\\ -\bar\gamma\bar\beta^{-2}(\bar\varrho+\bar\alpha)
\end{bmatrix},  \;\sigma(v) = \bar\sigma \begin{bmatrix}
v_1 & 0& 0\\
0& v_2&0\\ 0& 0& v_3
\end{bmatrix},
\end{align*}
with $|\bar\delta|\leqslant |\bar f|<\infty $ and  $\bar\sigma\in\R\setminus\{0\}$ finite, $0<\tau<\infty$, and $W_t =(W_t^1, W_t^2, W_t^3)$ an independent Wiener process in $\R^3$. Although the above  system is in the `toy category', considering the effects of time-periodic forcing and stochastic perturbations is relevant in many atmosphere-ocean applications to model, for example, seasonal and diurnal cycles in climate models  (e.g., \cite{ris89,Majda05, majda08, Pal01, Pal10, majdaqi19}).
It is well-known that for $\bar\sigma=0$ the system (\ref{Lorenz}) has an absorbing ball for all values of the parameters, since for $V(t,v) = |v|^2$ we have 
\begin{align*}
\frac{1}{2}\frac{d V}{dt}=\langle b(v),v\rangle&= -\bar\alpha\left(v_1-\frac{F_1}{2\bar\alpha}\right)^2- \bar\beta v_2^2-\bar\gamma\left(v_3 +\frac{\bar\varrho+\bar\alpha}{2\bar\beta^2}\right)^2 +\frac{\bar\alpha \bar\gamma\bar\beta^{-2}(\bar\varrho+\bar\alpha)^2+F_1^2}{4\bar\alpha},
\end{align*}
where we skip the explicit time dependence and $F_1(t) = \bar f\big(1+\bar\delta\sin\big(\frac{2\pi}{\tau} t)\big)$. Note that the drift and diffusion coefficients, $b$, $\sigma$, in (\ref{Lorenz}) are smooth and satisfy the growth conditions (\ref{diss_linear_growth}) outlined  in Remark \ref{Dep1}(a);  since  for $0< \varkappa_1, \varkappa_3<\infty$, and $\bar F_1 = \sup_{[0,\tau]} |F_1(t)|$ we have 
\begin{align*}
\langle b(v),v\rangle &\leqslant-\bar\alpha\left(1-\frac{\bar F_1}{4\bar\alpha\varkappa_1}\right)v_1^2 - \bar\beta v_2^2-\bar\gamma\left(1-\frac{\bar\varrho+\bar\alpha}{4\bar\beta^2\varkappa_3}\right)v_3^2 +\bar F_1\varkappa_1+\bar \gamma\bar \beta^{-2}(\bar\varrho+\bar\alpha)\varkappa_3,
\end{align*}
 where we used the fact that $|x|\leqslant \varkappa +\frac{1}{4\varkappa}|x|^2$ for $\varkappa>0$. Thus,  we have 
 \begin{align}\label{Lorenz_growth}
\langle b(v),v\rangle \leqslant L_{b_1}-L_{b_2}|v|^2, \qquad \|\sigma(v)\|_{\textsc{hs}}^2 \leqslant  L_\sigma\big(1+|v|^2\hspace{.03cm}\big), 
\end{align}
where $ L_\sigma = \bar\sigma$ and 
\begin{equation}\label{LLL}
L_{b_1} = \varkappa_1 \bar F_1+\varkappa_3\hspace{.02cm} \bar\gamma\bar\beta^{-2}(\bar\varrho+\bar\alpha), \quad L_{b_2} = \min\left(\bar\beta,  \bar\alpha\left(1-\frac{\bar F_1}{4\bar\alpha\varkappa_1}\right),  \bar\gamma\left(1-\frac{\bar\varrho+\bar\alpha}{4\bar\beta^2\varkappa_3}\right)\right).
\end{equation}
Thus, (\ref{Lorenz}) has global  solutions and it generates a stochastic flow of diffeomorphisms on~$\R^3$.

Next, note that the linear part in (\ref{Lorenz}) satisfies
\begin{align*}
\langle Av, v\rangle_{\R^3}\geqslant \mathfrak{C}_A \vert v\vert^2, \qquad \mathfrak{C}_A =\min\{\bar\alpha,\bar\beta, \bar\gamma\},
\end{align*}
and the nonlinear term  $G(v) = B(v,v)$ is given by a bilinear map $B(v,w) = (0, v_1w_3, -v_1w_2)$, $v,w\in \R^3$, which satisfies (see also~\cite{Keller})
\begin{align}\label{DLor}
\begin{cases}\langle B(v,w), w\rangle_{\R^3} = \langle (0, v_1w_3, -v_1w_2), (w_1,w_2,w_3)\rangle_{\R^3} = 0,\\
\langle B(v,w), u\rangle_{\R^3} = \langle (0, v_1w_3, -v_1w_2), (u_1, u_2,u_3)\rangle_{\R^3} 
=-\langle B(v,u), w\rangle_{\R^3}, \\
\vert B(v,w)\vert \leqslant \vert v\vert\vert w\vert.
\end{cases}
\end{align}
Consider $V(t,v) = \vert v\vert^p $ for some $1<p<\infty$, so that  
\begin{align*}
\partial_{v_i}V(t,v) = pv_i\vert v\vert^{p-2}, \;\; \partial_{v_iv_j}^2V(t,v) = p(p-2)v_iv_j\vert v\vert^{p-4}+\delta_{ij}p\vert v\vert^{p-2}.
\end{align*}
Next, we have 
\begin{align}\label{3.40}
\langle G(v)-G(w), v-w\rangle = \langle B(v-w,v), v-w\rangle\leqslant \vert v-w\vert^2\vert v\vert,
\end{align}
which follows from (\ref{DLor}) after some simple manipulations\footnote{\,This identity is obtained with the help of (\ref{DLor}) by noticing that one has \begin{align*}\langle B(v-w,v), v-w\rangle &= \langle B(v,v), v-w\rangle -\langle B(w,v), v-w\rangle \\&= \langle B(v,v), v-w\rangle -\langle B(w,w), v-w\rangle -\langle B(w,v-w), v-w\rangle \\ &=\langle G(v), v-w\rangle -\langle G(w), v-w\rangle,\end{align*} where the last term in the second line vanishes due to the fact that  $\langle B(u,w), w\rangle_{\R^3}=0$ $\forall u,w\in \R^3$. }, so that 
\begin{align}\label{Lorenz_2gen}
\hspace{.1cm}\LG^{(2)} V(t,v-w) &= p\big\langle -Av+Aw-G(v)+G(w), v-w\big\rangle_{\R^3}\vert v-w\vert^{p-2}\notag\\
& \hspace{.3cm} +\frac{1}{2}\sum_{ij=1}^3\bigg[(v_i-w_i)\big(\sigma_i(v)-\sigma_i(w)\big)\big(\sigma_j(v)-\sigma_j(w)\big)(v_i-w_j)p(p-2)\vert v-w\vert^{p-4}\notag\\[-.2cm] &\hspace{5.2cm}+ \delta_{ij}\,p\big(\sigma_i(v)-\sigma_i(w)\big)\big(\sigma_j(v)-\sigma_j(w)\big)\vert v-w\vert^{p-2}\bigg]\notag\\[.2cm]
&\leqslant p\vert v\vert \vert v-w\vert^p-p\beta\vert v-w\vert^p+{\textstyle \frac{1}{2}}\bar\sigma^2p(p-1)\vert v-w\vert^p,
\end{align}
where $\mathcal{L}^{(2)}$ is the two-point generator associated with   (\ref{Lorenz}). Next, choose $p$ such that for $v,w\in L^{p+1}(\Om,\mathcal{F}_{-\infty}^{\hspace{0.03cm}t},\mathbb{P})$ and $0<\E|v-w|^p$. Then, from Jensen's inequality we have 
\begin{equation}\label{Lorenz_jens}
0<\E|v-w|^p\leqslant (\E|v-w|^{p+1})^{p/(p+1)}<\infty,
\end{equation}
while the H\"older's inequality leads to (with $\|X\|_p:=(\E[|X|^p])^{1/p}$) 
\begin{align}\label{Lorenz_hold}
\E\big[|v||v-w|^p\big]=\big\||v||v-w|^p\big\|_1&\leqslant  \|v\|_{p+1}\big(\E|v-w|^{p+1}\big)^{p/{p+1}}.
\end{align}
The bounds (\ref{Lorenz_jens}) and (\ref{Lorenz_hold}) imply  that there exists a constant $1\leqslant \mathfrak{C}_p<\infty$ such that 
\begin{equation}\label{CCp}
\big(\E|v-w|^{p+1}\big)^{p/(p+1)}=\mathfrak{C}_p\,\E|v-w|^p.
\end{equation}
Combining (\ref{CCp}), (\ref{Lorenz_hold}), and (\ref{Lorenz_2gen}) leads to 
\begin{equation*}
\E\big[\LG^{(2)} V(t,v-w)\big] \leqslant -\lambda_p \,\E\big[V(t,v-w)\big], 
\end{equation*}
where $\lambda_p = p\big(\mathfrak{C}_A-\frot\bar\sigma (p-1)-\mathfrak{C}_p\|v\|_{p+1}\big)$. Now, for $v_t = \phi(t,s,\om,\xi), w_t = \phi(t,s,\om,\eta)$ solving (\ref{Lorenz}), we have from the above 
\begin{equation}\label{exx1}
\E\big[\LG^{(2)} V\big(t,\phi(t,s,\xi)-\phi(t,s,\eta)\big)\big] \leqslant -\lambda_p(t,s) \,\E\big[V\big(t,\phi(t,s,\xi)-\phi(t,s,\eta)\big)\big], 
\end{equation}
so that combining  It\^o's lemma
\begin{equation*}
d \E\big[V\big(t,\phi(t,s,\xi)-\phi(t,s,\eta)\big)\big] = \E\big[\LG^{(2)} V\big(t,\phi(t,s,\xi)-\phi(t,s,\eta)\big)\big],
\end{equation*}
with (\ref{exx1}) we obtain 
\begin{equation}\label{exx2}
\E\big[V\big(t,\phi(t,s,\xi)-\phi(t,s,\eta)\big)\big] \leqslant  \,\E\big[V(s,\xi-\eta)\big]\exp\left(-{ \int_s^t\lambda_p(r,s)dr}\right). 
\end{equation}

\begin{figure}[t]
\captionsetup{width=1\linewidth}

\vspace{-0.cm}
\includegraphics[width = 16cm]{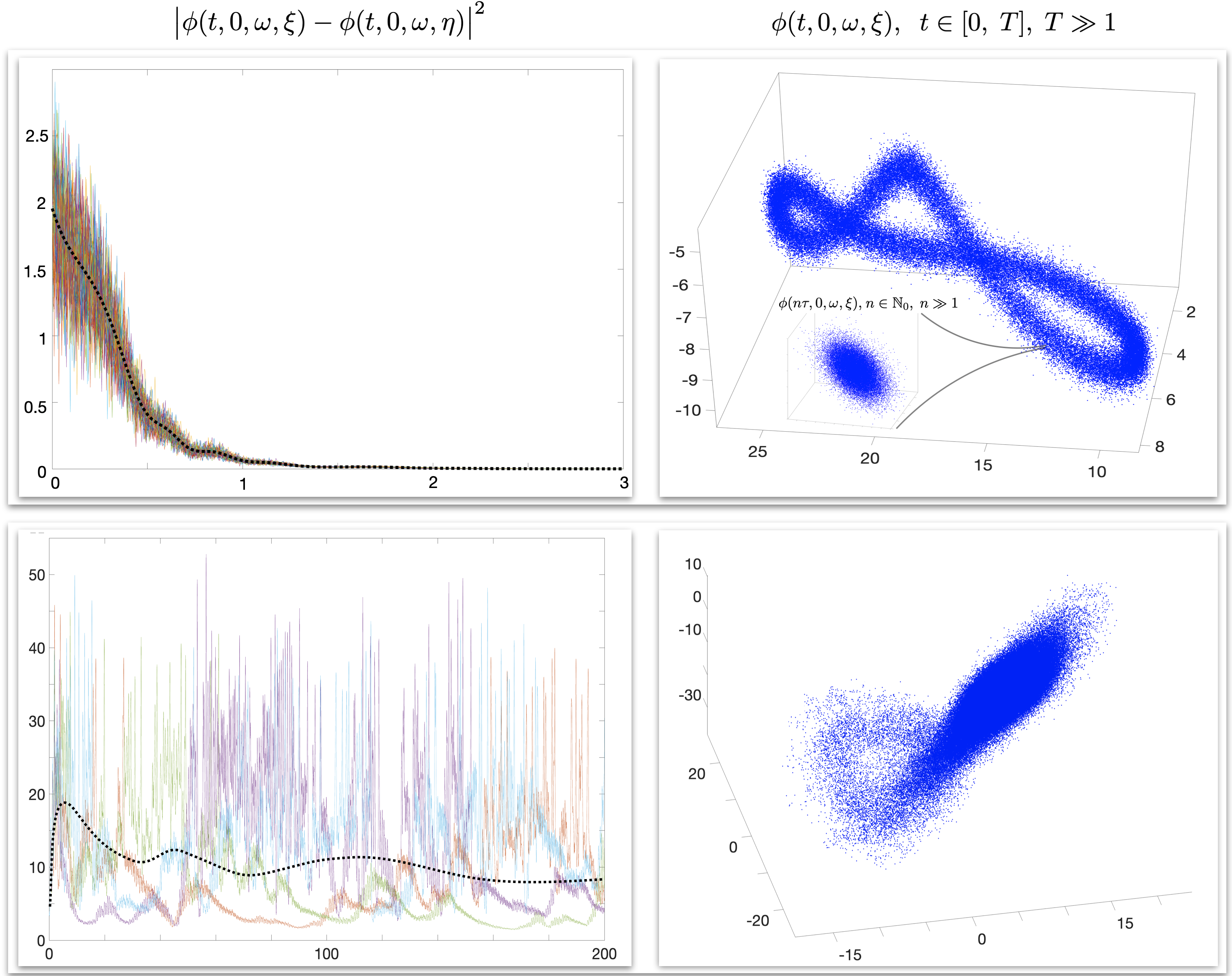}\\[.5cm]

\vspace{-.5cm}\caption{\footnotesize
Illustration of some aspects of the dynamics of the stochastic Lorenz model with time-periodic forcing (\ref{Lorenz}) with the flow of solutions $\{\phi(t,0,\om,\ccdot),\;t\geqslant 0\}$ in two different regimes. The top row corresponds to the regime in which the time-periodic measure exists and is supported on stable random periodic orbits of (\ref{Lorenz}); the inset of the top-right figure shows a finite sample from this measure on a Poincare section (i.e., on the subspace $\R^3$ of $\XX = \Sc\times\R^3$).  The top-left inset illustrates the relationship in (\ref{exx2}) in the case when $\lim_{t\rightarrow \infty} \lambda_p(t)<0$ and random periodic orbits exist (see text and Theorem \ref{Rand_SOL}); colours denote path-wise evolution of $|\phi(t,0,\om,\xi)-\phi(t,0,\om,\eta)|^2$ for fixed $\xi,\eta$, and the dotted black line denotes $\E|\phi(t,0,\ccdot,\xi)-\phi(t,0,\ccdot,\eta)|^2$. The bottom row illustrates a regime where  $\lim_{t\rightarrow \infty}|\phi(t,0,\om,\xi)-\phi(t,0,\om,\eta)|> 0$ and existence of random periodic orbits and periodic measures cannot be guaranteed.  Parameters in (\ref{Lorenz}) are: (top row) $\bar\alpha = 7.3, \bar\beta = 26, \bar\gamma= 7, \bar\varrho = 10$, $\bar f = 100, \bar\delta=0.9, \tau = 1, \bar\sigma=0.2$, and (bottom row) $\bar\alpha = 10, \bar\beta = 1,\bar\gamma= 8/3, \bar\varrho = 28$, $\bar f=23, \bar\delta=0.9, \tau = 1, \bar\sigma=0.2$. }\label{lor_fig}
\vspace*{-.5cm}
\end{figure}

Thus,  in order for Assumption \ref{A2.1}(ii) to hold, it is sufficient to require that 
$$ \mathfrak{C}_A-\frot\bar\sigma (p-1)- \limsup_{(t-s)\rightarrow\infty} \frac{1}{t-s}\int_s^t\mathfrak{C}_p(u,s)\|\phi(u,s,\xi)\|_{p+1}du>0.$$
Finally, we  choose $p=2$, so that  $V(t,v) = |v|^2$ and note that (see  Appendix \ref{app_hasminski})
$$\lim_{(t-s)\rightarrow\infty}\|\phi(t,s,\xi)\|_{3} =  \left(27\right)^{1/6}\left(\frac{L_{b_1}+L_\sigma}{L_{b_2}-L_\sigma}\right)^{1/2},$$ 
so that  simple but tedious algebraic manipulations lead to 
 \begin{equation}\label{constr}
 \hspace{1cm}\mathfrak{C}_A-\frot \bar\sigma- \bar{\mathfrak{C}}\left(\frac{L_{b_1}+L_\sigma}{L_{b_2}-L_\sigma}\right)^{1/2}>0, \quad \qquad \bar{\mathfrak{C}} = (27)^{1/6} \hspace{-0.3cm}\lim_{(t-s)\rightarrow\infty}\frac{1}{t-s}\int_s^t\mathfrak{C}_2(u,s)du, 
 \end{equation}
 where $L_{b_1}$, $L_{b_2}$, and $L_{\sigma}$ for the system (\ref{Lorenz}) are given in (\ref{LLL}).

\smallskip
Therefore, by Theorem \ref{Rand_SOL}, we conclude that the time-periodically forced stochastic Lorenz equation (\ref{Lorenz}) admits a family of periodic measures $\{\mu_{t}: t\in \Sc\}\subset \PP(\R^3)$ supported on stable periodic solutions of (\ref{Lorenz}) lifted to $\Sc\times \R^3$.  For $\bar \sigma\ne 0,$ Assumption~\ref{Homa} holds in addition to Assumption \ref{A2.1},  and Theorem \ref{Ps_erg} implies existence of ergodic $\tau$-periodic measures $\tilde{\mu}_t = \delta_{t\; \mmod \tau}\otimes\mu_t\in \PP\big(\Sc)\otimes\PP(\R^3\big)$ in the sense of Definition~\ref{PS-eR}. 
Numerical illustration of the convergence in (\ref{exx2}) is provided in Figure \ref{lor_fig} alongside a sample from the density of the ergodic measure supported on the attractor containing stable random periodic orbits. 

\end{exa}


\section{Linear response in the random time-periodic regime}\label{Linear response_per}
In this section, we derive a general formula  for the {\it linear response function} which characterises the change of  a statistical observable in response to  small perturbations of SDE dynamics  with a time-periodic ergodic probability measure. The results presented below build on and extend the derivations obtained  for time-dependent stochastic systems in \cite{Madja10}. First, we derive the  linear response formula associated with   perturbations of dynamics with time-periodic measures, and we represent it via formulas exploiting the  asymptotic statistical properties of the unperturbed dynamics; in line with terminology  from statistical physics, these are termed  {\it fluctuation-dissipation} formulas. In Theorem~\ref{FDTI}, we derive the fluctuation-dissipation formulas in the case when only the unperturbed dynamics has a time-periodic ergodic probability measure. In Theorem \ref{FDTII} we consider the linear response  associated with perturbations of  dynamics with a time-periodic ergodic probability measure under stronger conditions when perturbed dynamics also has a time-periodic ergodic measure. During the revision of the manuscript, we become aware of related results derived independently in \cite{X_chen} for non-autonomous SDE's;  those results are complementary to ours since they are confined to a finite time interval in the non-autonomous case with a restricted class of perturbations, and they do not deal with perturbations of time-periodic ergodic probability measures.  We conclude with some examples of the linear response for the  periodically forced stochastic Lorenz~model used earlier in Example~\ref{Lor_ex}. In principle, the results discussed below apply to a wider class of SDE's generating time-periodic measures under less stringent conditions than those in Assumption \ref{A2.1}; however, establishing the existence and ergodicity of such measures in a more general setting is not trivial and it is beyond the scope of this work.

\subsection{Setup and assumptions}
Consider the following SDE on $\Rd$ for $t\geqslant s$, $s\in \R$, 
\begin{align}\label{PSDE}
dX_t^{\alpha} = \hat b(\alpha(t),t, X_t^{\alpha})dt+ \hat\sigma(\alpha(t),t, X^{\alpha}_t)dW_{t-s}, \qquad 
 X^{\alpha}_{s}=x,  
\end{align}
 where  the maps $t\mapsto\hat b(0,t,\ccdot) = b(t,\ccdot)$, $t\mapsto\hat\sigma(0,t,\ccdot) = \sigma(t,\ccdot)$, are $\tau$-periodic and coincide with the coefficients in (\ref{NSDE});  $\alpha(\ccdot)\in \mathcal{C}^1_\infty(\R;\R)$ will be assumed sufficiently small in the sequel. Allowing for the explicit time dependence in $\alpha(t)$ enables one to consider time-dependent changes in the coefficients of (\ref{PSDE}) relative to those in the original dynamics for $\alpha=0$; for example, one can think of changes in the `climatological' forcing (e.g., \cite{ Abramov07,Abramov09, Majda05, Madja10, Madja10b,grit99,grit02,grit07,grit08,majdaqi19}) which is relevant for considerations in atmosphere-ocean science.   

Similar to~\S\ref{Erg_S},  we consider the Wiener probability space $(\Omega,\mathcal{F},\p)$, $\Omega:=\mathcal{C}_0(\R,\Rm)$, with $\mathcal{F}$ the Borel $\mathfrak{S}$-algebra on $\Om$, and  the probability measure $\p$ on $(\Om, \mathcal{F})$ induced by the $m$-dimensional Wiener process $W_t$. Furthermore,  we assume that there exists a proper  interval $\Ac\subseteq \R$  containing $\alpha=0$ such that, for all  $\alpha\in \Ac$, the coefficients $\hat b(\alpha,t,\ccdot)$, $\hat\sigma(\alpha,t,\ccdot)$ are sufficiently regular for  (\ref{PSDE}) to generate a forward stochastic flow  on $\Rd$   (see~Theorem \ref{SDEflow}); i.e.,  
$$ X^\alpha_t(\om) = \phi^\alpha(t,s,\om,x),\qquad s\leqslant t \quad \p\,\textrm{-\,a.s.}$$
The forward stochastic flow $\big\{\phi^\alpha(t, s, \ccdot, \ccdot): \;  s\leqslant t\big\}$ induced by~(\ref{PSDE})  has a one-point generator  
\begin{align}\label{La}
\mathcal{L}_{\alpha} = \sum_{i=1}\hat b_i(\alpha(t), t,x)\partial_{x_i}+\frot\sum_{i,j=1} \hat a_{ij}(\alpha(t),t,x)\partial^2_{x_ix_j}, \qquad \hat a:=\hat\sigma\hat\sigma^*.
\end{align}
As in the previous sections (see (\ref{NSDE1})), we lift the SDE (\ref{PSDE})  to  $\R\times\Rd$ to represent the dynamics~as 
\begin{align}\label{PSDEc}
d\tilde{X}^\alpha_t = \tilde{b}(\alpha(t),\tilde{X}^\alpha_t)dt +\tilde{\sigma}(\alpha(t), \tilde{X}_t^\alpha)d\tilde{W}_{t-s}, \qquad \tilde{X}^{\alpha}_s = \tilde x, \quad s\leqslant t,
\end{align}
where  $\tilde{x} = (s,x)\in \R\times\Rd$, and 
$$\tilde{b}\big(\alpha,\tilde x\big) = \big(1,\hat b(\alpha,s,x)\big)^{{*}},\quad \tilde{\sigma}\big(\alpha,\tilde x\big) = 
\begin{pmatrix}
0 &0\\
0 & \hat \sigma(\alpha,s,x)
\end{pmatrix}.$$

Finally, shifting the time $t\rightarrow t+s$ allows one to  represent the dynamics  in (\ref{PSDEc}) for $t\in \Rp$ as 
\begin{align}\label{PSDEc2}
d\tilde{X}^\alpha_{t+s} = \tilde{b}(\alpha(t+s),\tilde{X}^\alpha_{t+s})dt +\tilde{\sigma}(\alpha(t+s), \tilde{X}_{t+s}^\alpha)d\tilde{W}_{t+s}, \quad \tilde{X}^{\alpha}_s = \tilde x, 
\end{align}
where $\tilde{W}_{t+s} = \tilde{W}_{t+s}(\theta_{-s}\om)$ due to (\ref{tht_base2}). In what follows, we will always assume that the SDE with $\alpha=0$ satisfies the conditions of Theorem \ref{Rand_SOL}. 

The lifted process  $\tilde \Phi^\alpha:\, \Rp\times \Om\times \R\times \Rd \rightarrow\R\times\Rd$, is defined analogously to~(\ref{Lift}); namely 
\begin{align}\label{Lift_alph}
\tilde \Phi^\alpha\big(t, \om, \tilde x\big) := \big( t+s, \;\;\phi^\alpha(t+s,s,\theta_{-s}\om,x)\big), \quad \tilde x:=(s,x)\in \R\times\Rd, \;t\in \R^+. 
\end{align}
Note that if  the flow $\phi^\alpha$ is induced by the solutions of  (\ref{PSDE}) with the coefficients $\hat b\big(\alpha(t),t,x\big)$ and $\hat \sigma\big(\alpha(t),t,x\big)$ which are time periodic for all $\alpha\in \mathcal{A}$,  $\tilde \Phi^\alpha$ can be represented on a flat cylinder $[0,\,\hat\tau)\times \Rd$, $[0,\,\hat\tau)\eqsim \R\Mod\hat\tau$, $0<\hat\tau<\infty$, and the results of \S\ref{Erg_S} hold; one obvious case is for $\alpha =0$  when, by construction, the coefficients are time periodic with period $\tau$. If $\hat\tau=\tau$ for any $\alpha\in \mathcal{A}$, both the unperturbed ($\alpha=0$) and perturbed ($\alpha\ne0$) dynamics can be considered on the same cylinder. We will consider such a case in the last theorem of this section (Theorem \ref{FDTII}).

Similarly, the  transition evolutions $(\tilde{\mathcal{P}}^\alpha_t)_{t\in \Rp}$,  and their duals $(\tilde{\mathcal{P}}^{\alpha*}_t)_{t\in \Rp}$ can be defined through (\ref{tldP})-(\ref{tldP*}) with the transition kernel  $\tilde{P}^\alpha\big(\tilde{x}; t, \tilde{A}\big) := \p\big(\{\om: \tilde\Phi^\alpha(t,\om,\tilde{x})\in \tilde{A}\,\}\big)$; namely\footnote{\,As before,  we abuse the notation and  $\mathbb{M}_\infty\big(\R\times\Rd\big)$ denotes functions measurable w.r.t.~$\mathcal{B}(\R)\times\mathcal{B}(\Rd)$; see~\S\ref{gennot}.} 
\begin{align}
\tilde{\mathcal{P}}^\alpha_t\varphi(\tilde{x}) &:= \int_{\R\times\Rd}\varphi(\tilde{y})\tilde{P}^\alpha(\tilde{x}; t, d\tilde{y}) \hspace{.8cm} \forall \;\varphi\in \mathbb{M}_\infty\big(\R\times\Rd\big),\label{atldP}\\[.2cm]
\tilde\mu^\alpha_{t+r}(\tilde A)=\big(\tilde{\mathcal{P}}^{\alpha*}_t\hspace{.03cm}\tilde{\mu}^\alpha_r\big)(\tilde{A}) &:= \int_{\R\times\Rd}\tilde{P}^\alpha(\tilde{x}; t, \tilde{A})\tilde{\mu}^\alpha_r(d\tilde{x}) \quad \forall\;\tilde{\mu}^\alpha_r\in \PP\big(\R)\times\PP(\Rd\big),\;r\in \Rp,\label{atldP*}
\end{align}
with the short-hand notation  $\tilde\mu^\alpha_r(d\tilde x) = \delta_{r}(s)ds\otimes \mu^\alpha_{r}(dx)$ for the skew-product probability measures in the fibre bundle  $\PP(\R)\otimes \PP(\Rd)$, where $\tilde\mu^\alpha_r\in \PP(\R)\otimes \PP(\Rd)$ and $\mu^\alpha_r\in \PP(\Rd)$; see, e.g.,~\cite{sin} for more details concerning the structure of skew-product fibre bundles of probability measures. The definition in (\ref{atldP}) can be extend to $\varphi\in \mathbb{M}\big(\R\times\Rd\big)$ in a standard fashion.

The generator of the lifted one-point motion is given by  $\tilde{\LG}_{\alpha} = \partial_s+\LG_{\alpha}$ (see Definition \ref{inf_gnrt}). By construction (see, e.g., \cite{Madja10}), one can check that if $\mu^\alpha_t\in \PP(\Rd)$ is a solution of  forward Kolmogorov equation with the operator $\mathcal{L}^*_\alpha$, then $\tilde\mu^\alpha_t = \delta_{t}\otimes\mu^\alpha_t $ solves the lifted forward Kolmogorov equation with $\tilde{\mathcal{L}}^*_\alpha$ in the skew-product fibre bundle   $\PP(\R)\otimes\PP(\Rd)$.
 
In the sequel,  we derive fluctuation-dissipation formulas associated with  the  linear response for time-asymptotic  SDE's dynamics  in the random time-periodic regime with the ergodic measure $\bar{\tilde \mu}\in \PP\big(\Sc\big)\otimes\PP\big(\Rd\big)$, as in Theorem \ref{Ps_erg}.  We start with the definition of a {\it linear response function} which  approximates changes in the statistical observables due to sufficiently small perturbations of the unperturbed dynamics with a time-periodic ergodic probability measure.

\begin{definition}[{\bf Linear Response Function}]\rm
Assuming that $\E[\varphi(\tilde\Phi^\alpha)]\in L^1(\tilde\mu_0)$,  consider a family of statistical observables  
\begin{align}\label{Fobs}
{\mathbb{F}}_{\varphi}^{\tilde\mu_0}(t,\alpha) &=  \int_{\Sc\times\Rd}\E\big[\varphi\big(\tilde\Phi^{\alpha}(t,\ccdot,\tilde{x})\big)\big]\tilde{\mu}_0(d\tilde{x}) = \big\langle \tilde{\mathcal{P}}_t^\alpha \varphi, \tilde{\mu}_0\big\rangle,  \quad\varphi \in \mathcal{C}_\tau\big(\R\times\Rd\big),
\end{align}
where 
\begin{equation}\label{Ctau}
\mathcal{C}^2_\tau(\R\times\Rd):=\big\{\varphi\in \mathcal{C}^2(\R\times\Rd): \;\varphi(t+\tau,\ccdot)=\varphi(t,\ccdot), \;t\in \R\big\},
\end{equation}
the transition evolution $\{\tilde{\mathcal{P}}^\alpha_t\}_{t\in \Rp}$ is induced by $\tilde\Phi^\alpha$ in (\ref{Lift_alph}), and $\tilde{\mu}_0\in \PP\big(\Sc\big)\otimes\PP\big(\Rd\big)$ is  the  probability measure on the initial condition in (\ref{PSDEc2}), which is assumed throughout to be given by the $\tau$-periodic ergodic (skew-product) probability measure associated with the `unperturbed' dynamics with $\alpha=0$. If there exists a locally  integrable function $\mathcal{R}_{\varphi}^{\tilde\mu_0}$  such that the Gateaux derivative of $\tilde{\mathbb{F}}_{\varphi}^{\tilde\mu_0}(\ccdot,\alpha)$ at $\alpha=0$ satisfies 
\begin{align}\label{F_prt}
\Delta \mathbb{F}_{\varphi,\vartheta}^{\tilde\mu_0} (t):= \frac{d}{d \varepsilon}{\mathbb{F}}_\varphi^{\tilde\mu_0}(t, \varepsilon\vartheta)\Big|_{\varepsilon =0} =\int_0^t\mathcal{R}_{\varphi}^{\tilde\mu_0}(t-r, r)  \vartheta(r)dr, \qquad \vartheta\in \mathcal{C}^1_\infty(\Rp,\R), \quad \vartheta(0)=0,
\end{align}
we say that $\mathcal{R}_{\varphi}^{\tilde\mu_0}$ is a {\it linear response function} due to  perturbations of  the statistical observable~${\mathbb{F}}_{\varphi}^{\tilde\mu_0}$. 
\end{definition}
In other words, $\mathcal{R}_{\varphi}^{\tilde\mu_0}$ can be defined if the functional  ${\mathbb{F}}_{\varphi}^{\tilde\mu_0}(\,\cdot\,,\alpha)$ is Gateaux differentiable at $\alpha =0$ in the direction of $\vartheta$, and the Gateaux derivative is linear and continuous in the neighbourhood of  $\alpha=0$.   The formula (\ref{F_prt}) can be interpreted as an $\mathcal{O}(\varepsilon)$ approximation of the  change of the statistical observable ${\mathbb{F}}_{\varphi}^{\tilde\mu_0}$ in response to a sufficiently small perturbation $ \varepsilon\vartheta(t)$ around $\alpha=0$. The explicit time dependence in the perturbation $\vartheta(t)$ enables one to consider the linear response to small time-dependent changes in the coefficients of (\ref{PSDEc2}) relative to those in the original dynamics for $\alpha\,\,{=}\,\,0$; for example, one can consider  changes in the ``climatological'' forcing  (e.g., \cite{ Abramov07,Abramov09, Majda05, Madja10, Madja10b,grit99,grit02,grit07,grit08,majdaqi19}).

 \smallskip
 Throughout the remainder of this section, we impose the following regularity conditions which reduce to Assumption~\ref{A2.1} when $\alpha= 0$, and which imply the smoothing property (e.g., \cite{Da Prato2}) of the transition evolutions  $(\tilde{\mathcal{P}}^\alpha_t)_{t\in \Rp}$ (in the $x$-component of the extended state space $\R\times\Rd$):

\begin{Assum}\rm\label{A2.1b} 
 Assume that there exists a proper interval $\Ac\subseteq\R$ containing $\alpha=0$, and that  the following conditions are satisfied for all $s\leqslant t$, $s\in \R$,  and for all $\alpha\in \Ac$: 
\begin{itemize} [leftmargin = .8cm]
\item[(i)]   $D_\alpha^n\hspace{.03cm}\hat b(\alpha,t,\ccdot)\in \tilde{\mathcal{C}}^\infty(\Rd)$, $|D_\alpha^n\hspace{.03cm}\hat b(\alpha,t,x)|\leqslant \mathfrak{C}_\ell(1+|x|^\ell)$, $0\leqslant \ell<\infty,\,0\leqslant \mathfrak{C}_\ell<\infty$, $n\in\mathbb{N}_0$, and \\[.05cm] 
$t\mapsto \hat b(\ccdot,t,\ccdot)$ is differentiable on $\Ac\times\R\times \Rd$. 

\vspace{.15cm}\item[(ii)] $D_\alpha^n\hspace{.03cm}\hat \sigma_k(\alpha,t,\ccdot)\in \bar{\bar{\mathcal{C}}}^\infty(\Rd)$, $n\in\mathbb{N}_0$, and  $t\mapsto \hat \sigma_k(t,\ccdot,\ccdot)$ is differentiable on $\Ac\times\R\times \Rd$,
and 
\begin{align*}
\vert \partial_t\partial_{\alpha}^n D_x^\beta\hat \sigma_k(\alpha,t,x)\vert< \mathfrak{C}<\infty, \quad (\alpha,t,x)\in \Ac\times\R\times\Rd, \quad 1\leqslant k\leqslant m,
\end{align*}
for any multi-index $\beta$, and $\hat\sigma_k$, $1\leqslant k\leqslant m$ the columns of $\hat\sigma$. 

\vspace{.15cm}\item[(iii)] $\text{Lie}\big(\hat \sigma_{1}(\alpha,t,x),\cdots,\hat\sigma_m(\alpha,t,x)\big) = \Rd$, for all $s\leqslant t$,  where
\begin{align*}
\text{Lie}\big( \hat\sigma_1,\cdots,\hat \sigma_m\big):= \text{span}\big\{ \hat \sigma_i, [\hat \sigma_i, \hat\sigma_j], [\hat \sigma_i, [\hat\sigma_j, \hat\sigma_k]], \cdots,\; 1\leqslant i,j,k\leqslant m\big\},
\end{align*}
and $[F, G\hspace{.04cm}](\alpha,t,x)$ is the Lie bracket between the vector fields $F$ and $G$ defined by 
\begin{align*}
 [F,G\hspace{.04cm}](\alpha,t,x):= D_xG(\alpha,t,x)F(\alpha,t,x)-D_xF(\alpha,t,x)G(\alpha,t,x).
\end{align*}

\end{itemize}
\end{Assum}

\begin{rem}\rm Assumption \ref{A2.1b}, which is a version of the H\"ormander condition, implies  the existence of a smooth density of the time-marginal probability measure on $(\Rd, \mathcal{B}(\Rd))$ induced by the law of the solutions of (\ref{PSDE}); i.e.,  $\mu_t^\alpha(dx) = \rho_t^\alpha(x)dx$, $\rho_t^\alpha\in \mathcal{C}^\infty_\infty(\Rd)\cap L^1_+(\Rd)$ for all $s\leqslant t$.  In order to simplify the subsequent derivations, we will abuse notation and set the following  
\begin{equation}\label{abus}
\tilde \mu_t^\alpha(d\tilde x)\equiv \tilde\rho_t^\alpha(\tilde x)d\tilde x\equiv\delta_t(s)ds\otimes\rho_t^\alpha(x)dx,
\end{equation}  
when dealing with the skew-product  probability measures $\tilde \mu_t^\alpha\in \PP\big(\R)\otimes \PP(\Rd\big)$, $\tilde \mu_t^\alpha= \delta_t\otimes\mu_t^\alpha$, of the lifted process on $\R\times\Rd$ associated with the SDE (\ref{PSDEc2}). This intuitive convention is consistent with the convention introduced in (\ref{PPhi}) for  transition evolutions.
\end{rem}

\smallskip
In addition to Assumption \ref{A2.1b}, the last theorem of this section, which is  concerned with the linear response when both the unperturbed and perturbed measures are $\tau$-periodic and ergodic, will require the following assumption (which coincides with Assumption \ref{A2.1} for $\alpha=0$).
 
\begin{Assum}\label{A2.1A}\rm\;
Let $V\in \mathcal{C}^{1,2}(\R\times\Rd, \Rp)$ such that $V(t,0) = 0$ for all $t\in \R$, and the coefficients $\hat b(\alpha,t,x), \,\hat\sigma(\alpha,t,x)$, in (\ref{PSDE}) be such that for $\alpha\in \mathcal{A}$, where $\mathcal{A}\subseteq \R$ is  a proper  interval containing $\alpha=0$, the following hold: 
\begin{itemize}[leftmargin=0.9cm]
\item[(i)] There exist functions $\lambda_{\alpha}\in L^1(\R)$ with  $\mathcal{A}\ni\alpha\mapsto\lambda_{\alpha}(t)$ bounded for $t\in \R$, and a constant $1\leqslant \mathfrak{C}<\infty$ such that for all $\xi, \eta\in L^p(\Om, \F_{-\infty}^s,\p)$ and some $1<p<\infty$ we have 
 \begin{align}\label{Lyap_f1}
\begin{cases}
\E\vert \xi\vert^p\leqslant \E[V(t,\xi)]\leqslant \mathfrak{C}\,\vert \xi\vert^p,\\[.2cm]
\displaystyle \E\big[\mathcal{L}_{\alpha}^{(2)}V(t,\xi-\eta)\big]\leqslant\lambda_{\alpha}(t) \E\big[V(t,\xi-\eta)\big],
\end{cases}
\end{align} 
 where  $\mathcal{L}_{\alpha}^{(2)}$ is defined  analogously to the two-point generator $\mathcal{L}^{(2)}$ in (\ref{2pp}) but based on the coefficients $\hat b(\alpha, t,\ccdot)$, $\hat \sigma(\alpha, t,\ccdot)$ of (\ref{PSDE}); moreover, $\mathcal{L}_{\alpha}^{(2)}\equiv \mathcal{L}^{(2)}$ for $\alpha=0$.

\item[(ii)] There exists $\bar{\bar{\lambda}}>0$ such that 
\begin{align}\label{Dissp32}
\sup_{\alpha\in \Ac }\bigg\{\limsup_{(t-s)\rightarrow\infty}\frac{1}{t-s}\int_{s}^t\lambda_\alpha(u)du\bigg\}<-\bar{\bar{\lambda}}<0.
\end{align}

\item[(iii)]  Given the one-point motion  $t\mapsto\phi^\alpha(t,s,\om,\xi)$ induced by (\ref{PSDE}) for  $\om\in \Om$, $\xi\in \Rd$, and $s\leqslant t$, there exists $0<\mathfrak{D}_\alpha<\infty$ independent of $s,t\in \R$ such that\footnote{\,This condition can be replaced by a stronger but more concrete constraint on the global existence of the $p$-th absolute moment of $\phi^\alpha$; see Lemma \ref{lem_app1} in Appendix \ref{app_hasminski}.}  for all $\xi\in L^{p}(\Om,\F_{-\infty}^s,\p)$
\begin{align} \label{Temp12}
\limsup_{(t-s)\rightarrow \infty}\E\left[V\big(t,\phi^\alpha(t,s,\xi) -\xi\big)\right]<\mathfrak{D}_\alpha,
\end{align}
where $\E\big[V\big(\phi^\alpha(t,s,\xi)\big)\big]:=\int_\Om V\big(\phi^\alpha(t,s,\om,\xi)\big)\p(d\om)$. 

\end{itemize}
\end{Assum}

\smallskip
\subsection{Preparatory lemmas}
\noindent We start with the following standard and preparatory results which utilise relatively well-known results from \cite{Da Prato2,Stroock08}, and are aimed at representing $\tilde{\mathcal{P}}^\alpha_{t}\varphi(\tilde{x}) -\tilde{\mathcal{P}}_{t}\varphi(\tilde{x})$ in the form amenable to further analysis in the context of the linear response. The main results are derived in \S\ref{LRsp}.

\vspace*{-.1cm}\begin{lem}\label{PertPropA}
Suppose that the coefficients $\hat b(\alpha,t,x)$, $\hat\sigma(\alpha,t,x)$ in the SDE (\ref{PSDE}) satisfy Assumption~\ref{A2.1b}  so that global solutions of (\ref{PSDE}) exist for all time, and  $\varphi\in \mathcal{C}^2(\R\times\Rd)$ is such that for any fixed $\tilde x=(s,x)\in \R\times\Rd$ and for all $\alpha\in \mathcal{A}$, \,$\E[D_x^\beta\varphi(\tilde\Phi^\alpha(t,\tilde x))]<\infty$, $ |\beta| \leqslant 2$, where $\{\tilde \Phi^\alpha(t,\ccdot,\ccdot)\!: \;t\in \Rp\}$ in (\ref{Lift_alph}) is generated  by  the lifted SDE~(\ref{PSDEc2}) on $\R\times\Rd$.

Then, the function $v(r,\tilde{x}) := \tilde{\mathcal{P}}^\alpha_{t-r}\varphi(\tilde{x})$ with $\tilde{x}:=(s,x)\in \R\times\Rd$, and $\tilde{\mathcal{P}}^\alpha_{t}$ defined in (\ref{atldP}),
is the unique solution of the {\it backward Kolmogorov equation} with the terminal condition
\begin{align}\label{bkeq}
\begin{cases}
\partial_rv(r,\tilde{x}) = -\tilde{\LG}_\alpha v(r,\tilde{x}), \quad 0\leqslant r\leqslant  t,\\
v(t,\tilde{x}) = \varphi(\tilde{x}).
\end{cases}
\end{align}
Moreover, for any $\varphi\in \mathcal{C}_\infty^{2}(\R\times\Rd)$, there exists a constant $\mathfrak{C} >0$ such that 
\begin{align}\label{bkeBB}
\Vert \tilde{\mathcal{P}}^\alpha_{t-r}\varphi\Vert_{2,\infty}\leqslant \mathfrak{C}\Vert \varphi\Vert_{2,\infty}, \quad 0\leqslant r\leqslant t,
\end{align}
where $ \displaystyle \Vert \varphi\Vert_{2,\infty} =\Vert \varphi\Vert_{\infty}+ \textstyle{\sum_{1\leqslant\vert\beta\vert\leqslant 2}}\Vert D_x^\beta\varphi\Vert_{\infty}.$
\end{lem}
\noindent {\it Proof.} See, e.g., \cite[Thm 4.8.11]{Kunita} or \cite{Da Prato2, Stroock08,Bertoldi07}. Sufficient conditions for $\E[D_x^\beta\varphi(\tilde\Phi^\alpha)]<\infty$,  $ |\beta| \leqslant 2$,   $\varphi\in \mathcal{C}^2(\R\times\Rd)$, $\alpha\in \mathcal{A}$, are given in Proposition \ref{PertProp}.  It can also be shown that, under Assumption \ref{A2.1b}, the problem (\ref{bkeq}) has unique classical solutions for $\varphi\in \mathcal{C}(\R\times\Rd)$  due to the smoothing property of~$(\tilde{\mathcal{P}}^\alpha_t)_{t\in \Rp}$  (e.g., \cite{Da Prato2}; the general case can be obtained by approximating   $\varphi\in \mathcal{C}(\R\times\Rd)$ by $\varphi_n\in \mathcal{C}_\infty(\R\times\Rd)$ converging uniformly to $\varphi$ on compact subsets of $\R\times\Rd$.

\begin{lem}\label{PertProp2}
Suppose that the conditions of Lemma \ref{PertPropA} hold. Then 
\begin{align}\label{PFormu}
\tilde{\mathcal{P}}^{\alpha}_{t}\varphi(\tilde{x})-\tilde{\mathcal{P}}_{t}\varphi(\tilde{x}) = \int_0^t\tilde{\mathcal{P}}_{r}\big( \tilde{\LG}_{\alpha}-\tilde{\LG}\big)\tilde{\mathcal{P}}^\alpha_{t-r}\hspace{.04cm}\varphi(\tilde{x})dr, \qquad \,0\leqslant r\leqslant t, \;\;\tilde{x}\in \R\times\Rd.
\end{align}
\end{lem}
\noindent {\it Proof.}
It can be obtained from  Lemma \ref{PertPropA} that the function $u(r,\tilde{x}) = \tilde{\mathcal{P}}^\alpha_{t-r}\varphi(\tilde{x}) -\tilde{\mathcal{P}}_{t-r}\varphi(\tilde{x})$,  $0\leqslant r\leqslant t$, $\tilde{x}\in \R\times\Rd$, uniquely  solves the  inhomogeneous Cauchy problem
(see, e.g., \cite{Stroock08,Da Prato2, Bertoldi07}) 
\begin{align}\label{inh_pde}
\begin{cases}
\partial_ru(r,\tilde{x}) = -\tilde{\LG} u(r,\tilde{x}) -\alpha(r)\,\hat{\mathfrak{f}}_{t,\alpha}(r,\tilde{x}),\\
u(t,\tilde{x}) =0,
\end{cases}
\end{align} 
where $\alpha\in \mathcal{C}^1_\infty(\R,\R)$,  $\alpha(r)\in\mathcal{A}$, and 
\begin{equation}\label{alph_g}
\alpha(r)\hat{\mathfrak{f}}_{t,\alpha}(r,\tilde{x}):=\big( \tilde{\LG}_{\alpha}- \tilde{\LG}\big) \tilde{\mathcal{P}}^\alpha_{t-r}\varphi(\tilde{x}).
\end{equation}

\smallskip
Next, consider  the solutions of  (\ref{PSDEc2}) with $\alpha=0$ represented through (\ref{Lift_alph}) as  $\tilde\Phi^0=\tilde\Phi(r,\om,\tilde x)$  where $\tilde\Phi$ is defined in  (\ref{Lift}) and solves the SDE (\ref{NSDE11}). Then,  by It\^o's formula 
\begin{align}\label{duP}
du(r,\tilde{\Phi}_r) &= \Big[ \partial_ru(r,\tilde{\Phi}_r)+\tilde{\LG}u(r,\tilde{\Phi}_r)\Big]dr+D_xu(r,\tilde{\Phi}_r)^{T}\tilde{\sigma}(r,\tilde{\Phi}_r)d\tilde{W}_r\notag\\
&=-\alpha(r)\,\hat{\mathfrak{f}}_{t,\alpha}(r,\tilde{\Phi}_r) dr+D_xu(r,\tilde{\Phi}_r)^{T}\tilde{\sigma}(r,\tilde{\Phi}_r)d\tilde{W}_r, 
\end{align}
where $\tilde \Phi_r\equiv \tilde\Phi(r,\om,\tilde x)$ to simplify notation. Combining (\ref{duP}) with (\ref{inh_pde}), and using the explicit form of $u(t,\tilde x)$ leads to 
\begin{equation*}
\hspace{1cm}\tilde{\mathcal{P}}^\alpha_{t}\varphi(\tilde{x}) -\tilde{\mathcal{P}}_{t}\varphi(\tilde{x}) = \int_0^t\alpha(r)\E\big[ \,\hat{\mathfrak{f}}_{t,\alpha}(r,\tilde{\Phi}(r,\tilde x)\big]dr=\int_0^t\tilde{\mathcal{P}}_{r}\big( \tilde{\LG}_{\alpha}- \tilde{\LG}\big) \tilde{\mathcal{P}}^\alpha_{t-r}\varphi(\tilde{x}) dr. \hspace{1.3cm}
\end{equation*}
The above identity is well-defined due to the underlying assumptions, and it  is discussed further in Proposition \ref{PertProp}. \qed

\begin{definition}[{$\alpha$-linearised generator}] \rm Given the infinitesimal generator $\tilde{\mathcal{L}}_\alpha = \partial_s+\mathcal{L}_\alpha$ with $\mathcal{L}_\alpha$ defined in (\ref{La}) and $(\hat b,\,\hat \sigma)$ satisfying Assumption \ref{A2.1b}, the {\it $\alpha$-linearised generator} is defined~by 
\begin{align}\label{Va}
\tilde{\mathcal{V}}\varphi(\tilde x) =  \mathfrak{b}(\tilde x)D_x\varphi(\tilde x)+\frac{1}{2}\text{Tr}\Big(\mathfrak{a}(\tilde x)D_x^2\varphi(\tilde x)\Big), \qquad \tilde x = (s,x),\;\;\; \varphi\in \mathcal{C}^2\big(\R\times\Rd\big),
\end{align}

\vspace{-.2cm}
\noindent where 
 $$\mathfrak{b}(\tilde x) = \partial_\alpha \hat{b}(\alpha,s,x)\big|_{\alpha=0},  \qquad \mathfrak{a}(\tilde x):= \sigma(s,x)H^*(s,x)+\sigma^*(s,x)H(s,x),$$
with   $H_{ik}(s,x)=\partial_\alpha\hat{\sigma}_{ik}(\alpha,s,x)\big|_{\alpha=0}$, $1\leqslant k \leqslant m, \; 1\leqslant i\leqslant d$, and $\hat\sigma(0,t,\ccdot) = \sigma(t,\ccdot)$ as in (\ref{PSDE}).  The $L^2(\tilde\mu_t)$ dual of $\tilde{\mathcal{V}}$ is given by  
\begin{align}\label{DUal}
\tilde{\mathcal{V}}^*\tilde \rho_t(\tilde x) = -D_x\big( \mathfrak{b}(\tilde{x})\tilde\rho_t(\tilde x)\big)+\frac{1}{2}\text{Tr}\Big(D_x^2\big(\mathfrak{a}(\tilde x)\tilde\rho_t(\tilde x)\big)\Big), 
\end{align}
where $\tilde\rho_t(\tilde x)$ is understood in the sense of (\ref{abus}).
\end{definition}

Note that the properties of $\mathfrak{a}(\tilde{x})$ and $\mathfrak{b}(\tilde{x})$  are fully controlled through Assumption \ref{A2.1b}.

\begin{prop}\label{PertProp}
Suppose that the conditions of Lemma \ref{PertPropA} are satisfied, and  consider a  function $\mathfrak{f}_{t,\alpha}: [0, \,T)\times\R\times\Rd\rightarrow \R$,  $\alpha\in \Ac$, defined~by
\begin{align}\label{CrucialPer}
\mathfrak{f}_{t,\alpha}(r,\tilde{x}) = \tilde{\mathcal{V}}\,\tilde{\mathcal{P}}^{\alpha}_{t-r}\varphi(\tilde{x}), \qquad 0\leqslant r\leqslant t\leqslant T.
\end{align}
Then, $\mathfrak{f}_{t,\alpha}\in \mathcal{C}(\R\times\Rd)$ and $\mathfrak{f}_{t,\alpha}<\infty$ for any fixed $\tilde{x}\in \R\times\Rd$.

If  Assumption \ref{A2.1b} holds for $0\leqslant \ell<\infty$,  and the initial condition in  the lifted SDE~(\ref{PSDEc2}) has  $p\geqslant \max(2,\ell)$ finite moments, then there exists a constant $\mathfrak{C} = \mathfrak{C}(T, k,\varphi)>0$ such that for any fixed $\tilde{x}\in \R\times\Rd$ one has 
\begin{align}\label{CrucialPer2}
\sup_{0\leqslant  r\leqslant t\leqslant T} \E \big|\mathfrak{f}_{t,\alpha}(r,\tilde{\Phi}(r,\tilde x))\big|\leqslant   \mathfrak{C}<\infty,
\end{align}
where $\{\tilde \Phi(r,\ccdot,\ccdot)\!: \;r\in \Rp\}$ is the RDS (\ref{Lift}) generated by  the SDE (\ref{PSDEc2}) with~$\alpha=0$ (or (\ref{NSDE11})).

The sufficient condition for (\ref{CrucialPer2}) to hold as $T\rightarrow \infty$ is  that $(\hat b,\,\hat\sigma)$ satisfy  conditions (\ref{diss_linear_growth})-(\ref{diss_linear_growth_mom}) with $p=\max(2,\ell)$. If  $\varphi(\ccdot,x)\leqslant \mathfrak{C}_l(1+|x|^l)$, $0\leqslant l<\infty,\,0\leqslant\mathfrak{C}_l<\infty$,  then (\ref{CrucialPer2}) holds for $T\rightarrow \infty$ when $(\hat b,\,\hat\sigma)$ satisfy  (\ref{diss_linear_growth})-(\ref{diss_linear_growth_mom}) with $p=\max(2,\ell,l)$; i.e., the conditions $\E[D_x^\beta\varphi(\tilde\Phi^\alpha)]<\infty$, $|\beta| \leqslant 2$,  in Lemma \ref{PertPropA} can be replaced by assuming a polynomial growth of~$\varphi$.
\end{prop}

\begin{rem}\rm
Note that for a dissipative dynamics satisfying (\ref{diss_linear_growth})  the dissipation coefficient $L_{b_2}$ might not be large enough to satisfy (\ref{diss_linear_growth_mom}) with a given $p\geqslant 2$. Thus, not all dissipative dynamics automatically satisfy Proposition \ref{PertProp} for all time.
\end{rem}

\noindent {\it Proof.} 
For  $(s,x)\in \R\times\Rd,$ one can obtain directly from  (\ref{Va}) that 
\begin{align}\label{gtalph}
\mathfrak{f}_{t,\alpha}(r,\tilde x) &= \tilde{\mathcal{V}}\tilde{\mathcal{P}}_{t-r}^{\alpha}\varphi(\tilde x)= \mathfrak{b}(\tilde x)D_x\tilde{\mathcal{P}}^{\alpha}_{t-r}\varphi(\tilde x)+\frac{1}{2}\text{Tr}\Big( \mathfrak{a}(\tilde x)D_x^2\tilde{\mathcal{P}}^{\alpha}_{t-r}\varphi(\tilde x)\Big).
\end{align}
The regularity and growth conditions of the coefficients $(\hat b, \,\hat \sigma)$ imposed in  Assumption \ref{A2.1b} ensure the existence of global solutions to (\ref{PSDEc2})  which are represented via $ \tilde\Phi^\alpha(t-r,\ccdot,\tilde x)$, $t-r\in \Rp$, and generate $\tilde{\mathcal{P}}^{\alpha}_{t-r}$ in~(\ref{atldP}).  If $p\geqslant 2$ moments of the initial condition of (\ref{PSDEc2}) are finite and $\E[D_x^\beta\varphi(\tilde\Phi^\alpha(t,\tilde x))]<\infty$, $ |\beta| \leqslant 2$, $t\leqslant T$, then 
by the assumption on $(\hat b, \,\hat \sigma)$ and Lemma~\ref{PertPropA}, we have $\mathfrak{f}_{t,\alpha}\in \mathcal{C}(\R\times\Rd)$ and $\mathfrak{f}_{t,\alpha}<\infty$ for any fixed $\tilde{x}\in \R\times\Rd$. 

Regarding (\ref{CrucialPer2}),  the polynomial growth of $\hat b$ combined with  the  standard calculation utilising It\^o's formula guarantees the existence of $\max(2,\ell)$ finite moments of the solution  for $T<\infty$ (Theorem~\ref{SDEflow}). Thus, (\ref{CrucialPer2}) follows by the Cauchy-Schwarz inequality applied to $\E|\mathfrak{f}_{t,\alpha}|$ and the finiteness of the moments of $\mathfrak{a}, \mathfrak{b}$ for $T<\infty$.  

Considering (\ref{CrucialPer2}) for  $T\rightarrow \infty$, may require additional dissipative constraints on  the drift $\hat b$, as outlined  below.  Moreover, we specify two explicit classes of $\varphi\in \mathcal{C}_\infty^{2}(\R\times\Rd)$ for which $\E[D_x^\beta\varphi(\tilde\Phi^\alpha(t,\tilde x))]<\infty$, $ |\beta| \leqslant 2$,  (\ref{CrucialPer2}) holds for all time.  First, if $\varphi\in \mathcal{C}_\infty^{2}(\R\times\Rd)$, then 
\begin{align}\label{bb1}
\sup_{0\leqslant r\leqslant t\leqslant T}\E| \mathfrak{f}_{t,\alpha}(r,\tilde \Phi(r,\tilde x)|  \leqslant {C}_T\sup_{0\leqslant r\leqslant t\leqslant T}\Vert \tilde{\mathcal{P}}^\alpha_{t-r}\varphi\Vert_{2,\infty}.
\end{align}
The second part of Lemma \ref{PertPropA}, the term on the right of (\ref{bb1}) is bounded by  ${C}\Vert \varphi\Vert_{2,\infty}$ so that one can set explicitly $\mathfrak{C} = {C}_T\,{C}\Vert \varphi\Vert_{2,\infty}$ in (\ref{CrucialPer2}). Given, the existence of global solutions of (\ref{PSDEc2}) for all time (Assumption \ref{A2.1b} and Theorem \ref{SDEflow}), the bound (\ref{CrucialPer2}) can be extended to $T\rightarrow \infty$  provided that the dissipative conditions  (\ref{diss_linear_growth})-(\ref{diss_linear_growth_mom}) hold for $(\hat b, \hat\sigma)$ with $p=\max(2,\ell)$; so that $\E|\mathfrak{b}(\tilde\Phi(r,\tilde x))|<\infty$, $\E|\mathfrak{a}(\tilde\Phi(r,\tilde x))|<\infty$, $r\in \Rp$, in (\ref{gtalph}). More generally, if $\varphi\in \mathcal{C}^{2}(\R\times\Rd)$ and  $\varphi(\ccdot,x)\leqslant \mathfrak{C}_l(1+|x|^l)$ with $0\leqslant l<\infty, \,0\leqslant \mathfrak{C}_l<\infty$, the bound (\ref{CrucialPer2}) can be extended to $T\rightarrow \infty$  provided that the dissipative conditions  (\ref{diss_linear_growth})-(\ref{diss_linear_growth_mom}) hold for $(\hat b, \hat\sigma)$ with $p=\max(2,\ell,l)$. Both assertions can be obtained through derivations analogous to Lemma \ref{lem_app1} in  Appendix \ref{app_hasminski}. \qed

\begin{cor}\label{PertProp_cor}
Let $\hat{\mathfrak{f}}_{t,\alpha}: [0, \,T)\times\R\times\Rd\rightarrow \R$,  $\alpha\in \Ac$, be defined~by  
\begin{equation}\label{alph_gg}
\alpha(r)\hat{\mathfrak{f}}_{t,\alpha}(r,\tilde{x}):=\big( \tilde{\LG}_{\alpha}- \tilde{\LG}\big) \tilde{\mathcal{P}}^\alpha_{t-r}\varphi(\tilde{x}), \qquad   \tilde x\in \R\times\Rd,  \;\;0\leqslant r\leqslant t\leqslant T,
\end{equation}
as in  (\ref{alph_g}) of  Lemma~\ref{PertProp2}; so that $\hat{\mathfrak{f}}_{t,\alpha}(r,\tilde{x}) = {\mathfrak{f}}_{t,\alpha}(r,\tilde{x})+\mathcal{O}(\alpha)$ with ${\mathfrak{f}}_{t,\alpha}$ defined in Proposition~\ref{PertProp}. Then, under the same conditions as those in Proposition \ref{PertProp},  $\hat{\mathfrak{f}}_{t,\alpha}\in \mathcal{C}(\R\times\Rd)$ and $\hat{\mathfrak{f}}_{t,\alpha}<\infty$ for any fixed $\tilde{x}\in \R\times\Rd$. Furthermore, for $p\geqslant 2$ chosen as in Proposition \ref{PertProp} one has 
\begin{align}\label{hatCrucialPer2}
\sup_{0\leqslant  r\leqslant t\leqslant T} \E \big|\hat{\mathfrak{f}}_{t,\alpha}(r,\tilde{\Phi}(r,\tilde x))\big|<\infty, \qquad \tilde{x}\in \R\times\Rd,
\end{align}
which can be extended to $T\rightarrow \infty$ in a way analogous to that  in Proposition \ref{PertProp}.
\end{cor}
\noindent {\it Proof.} This is a direct consequence of Proposition \ref{PertProp} and the fact that the $\mathcal{O}(\alpha)$ terms involve $D_x\tilde{\mathcal{P}}^{\alpha}_{t-r}\varphi(\tilde x)$, $D_x^2\tilde{\mathcal{P}}^{\alpha}_{t-r}\varphi(\tilde x)$,  with coefficients given by $D^n_\alpha \hat b$ and  $D^n_\alpha \hat \sigma$ which are controlled through Assumption \ref{A2.1b} and the polynomial bound on the growth of the derivatives of $\hat b$ w.r.t.~$\alpha$.

\vspace{.1cm}
\subsection{Linear response and fluctuation-dissipation formulas for time-periodic measures}\label{LRsp}

Here, in \S\ref{linresp_ss} we derive a general expression for the linear response function characterising the change in the statistical observable (\ref{Fobs}) to small perturbations of dynamics of an SDE  whose time-asymptotic dynamics is characterised by time-periodic ergodic probability measures~(see \S\ref{Erg_S}). This is followed in \S\ref{fdt_ss} by deriving a more tractable representation of the response function in terms of fluctuation-dissipation type formulas which allow one to express the change  in the statistical observables through statistical characteristics of the unperturbed~dynamics. 

\subsubsection{The linear response}\label{linresp_ss} First, we derive a general formula for the linear response function associated with perturbations of the time-asymptotic dynamics of the SDE (\ref{NSDE11}). The derived formula is equivalent to the one obtained  formally in \cite{Madja10}.  

\begin{theorem}[\textit{\textbf{Linear response}}]\label{Linear_resp_formula}
Suppose that Assumption~\ref{A2.1b} holds,  Assumption \ref{A2.1A} holds for $\alpha=0$, and that Proposition \ref{PertProp} is satisfied. Consider the family of transition semigroups $\{\tilde{\mathcal{P}}_t^\alpha: (t,\alpha)\in \Rp\times\Ac\}$ induced by the SDE (\ref{PSDEc2}) which for $\alpha =0$ admit $\tau$-periodic ergodic (skew-product) measures $(\tilde{\mu}_t)_{t\geqslant 0}$, $\tilde\mu_t\in \PP\big(\Sc\big)\otimes\PP\big(\Rd\big)$,  where  $\tilde\mu_t(d\tilde x) = \delta_{(t+s\!\!\mod \tau)}(s)\otimes\rho_t(x)dsdx$.

Then, given the observable ${\mathbb{F}}_{\varphi}^{\tilde\mu_0}(t,\alpha) = \big\langle \tilde{\mathcal{P}}_t^\alpha \varphi, \tilde{\mu}_0\big\rangle$ in (\ref{Fobs}),  for any $\varphi\in \mathcal{C}^2_\tau(\R\times\Rd)$ such that $\E[D_x^\beta\varphi(\tilde\Phi^\alpha(t,\ccdot))]\in L^1(\tilde\mu_0)$, $ |\beta| \leqslant 2$,   and the perturbation $\alpha(\cdot) = \varepsilon\hspace{.03cm}\vartheta(\cdot)\in \mathcal{C}_{\infty}^{1}(\Rp,\R)$, $\vartheta(0)=0$,  such that $\varepsilon\hspace{.03cm}\vartheta\in \mathcal{A}$, the following holds
\begin{equation}\label{LinR}
\Delta \mathbb{F}_{\varphi,\vartheta}^{\tilde\mu_0} (t) = \int_0^t\mathcal{R}^{\tilde{\mu}_0}_\varphi(t-r,r)\vartheta(r) dr,
\end{equation}
with the linear response function given by 
\begin{equation}\label{Resp_fun}
\mathcal{R}^{\tilde{\mu}_0}_\varphi(t-r,r) =   \int_{\Sc\times\Rd}\tilde{\mathcal{P}}_{t-r}\varphi(\tilde x)\big( \tilde{\mathcal{V}}^*{\tilde \rho}_r\big)({\tilde x}) d{\tilde x}, 
\end{equation}
where $\tilde{\mathcal{V}}^*$ is defined in (\ref{DUal}), and $\tilde \rho_r(\tilde x)d\tilde x$ is understood in the sense of (\ref{abus}).
\end{theorem}

\begin{rem}\label{pert_fact_rem}\rm 
$\Delta \mathbb{F}_{\varphi,\vartheta}^{\tilde\mu_0} (t)$ can be interpreted as the approximate $\mathcal{O}(\varepsilon)$ response of the observable $ {\mathbb{F}}_{\varphi}^{\tilde\mu_0}(t,0)$  in (\ref{PSDEc2})  to a sufficiently small perturbation $\alpha(t)=\varepsilon\vartheta(t)$, $\vartheta(0)=0$. Note that the linear response formula is solely based on  quantities defined for the `unperturbed' dynamics~($\alpha=0$).  Moreover,  the perturbation  $\alpha(t)=\varepsilon\vartheta(t)$ does not necessarily have to factorise the coefficients of (\ref{PSDE}) as, e.g.,  $\hat b(\alpha(t),t,x) = \hat b(0,t,x)+\varepsilon\hspace{.03cm}\mathfrak{F}({x})\vartheta(t)$. For example, consider the following 
$\hat b(\alpha(t),t,x) = -(2-\sin^2\big(x(1+\alpha(t))\big)x\cos^2(t)$ 
in (\ref{PSDE}) and the derivations below.  
\end{rem}

\noindent {\it Proof.} First, we set $\alpha = \varepsilon\vartheta\in \mathcal{A}$ and  show that under the assumptions of the proposition the following holds for any $\varphi\in \mathcal{C}^2_\tau(\R\times\Rd)$ such that $\E[D_x^\beta\varphi(\tilde\Phi^\alpha(t,\ccdot))]\in L^1(\tilde\mu_0)$, $ |\beta| \leqslant 2$ (see Proposition \ref{PertProp} for a sufficient condition for this to hold for all time),  
\begin{align*}
\lim_{\varepsilon\rightarrow 0}\frac{1}{\varepsilon}\big(\tilde{\mathcal{P}}_t^{\varepsilon\vartheta(t)}\varphi(\tilde{x}) -\tilde{\mathcal{P}}_t\varphi(\tilde{x})\big)  = \int_0^t\vartheta(r)\tilde{\mathcal{P}}_r(\tilde{\mathcal{V}}\tilde{\mathcal{P}}_{t-r}\varphi)(\tilde{x})dr, \qquad 0\leqslant r\leqslant t, \;\vartheta\in \mathcal{C}^1_\infty(\Rp,\R),
\end{align*}
with $\tilde{\mathcal{V}}$ defined in (\ref{Va}). To this end, it follows from Lemma \ref{PertProp2} that for $\varepsilon>0$ sufficiently small
\begin{align*}
\frac{1}{\varepsilon}\Big( \tilde{\mathcal{P}}_{t}^{\varepsilon\vartheta(t)}\varphi(\tilde{x})-\tilde{\mathcal{P}}_t\varphi(\tilde{x})\Big) &= \frac{1}{\varepsilon}\int_0^t\tilde{\mathcal{P}}_r\big( \tilde{\LG}_{\varepsilon\vartheta(r)}-\tilde{\LG}\big)\tilde{\mathcal{P}}_{t-r}^{\varepsilon\vartheta(r)} \varphi(\tilde{x}) dr=\int_0^t\vartheta(r)\E\big[\,\hat{\mathfrak{f}}_{t,\varepsilon}(r,\tilde{\Phi}(r,\tilde x))\big]dr,
\end{align*}
where $\tilde \Phi(r,\om,\tilde x)$ represents the solution of~(\ref{PSDEc2}) with $\alpha=0$ (or the solution of (\ref{NSDE11})).

 By Proposition \ref{PertProp}, Corollary \ref{PertProp_cor},  and the dominated convergence theorem, we have 
\vspace{.1cm}
\begin{align*}
\lim_{\varepsilon\rightarrow 0}\frac{1}{\varepsilon}\Big( \tilde{\mathcal{P}}_{t}^{\varepsilon\vartheta(t)}\varphi(\tilde{x})-\tilde{\mathcal{P}}_t\varphi(\tilde{x})\Big)&=\int_0^t\vartheta(r)\E\Big[\lim_{\varepsilon\rightarrow 0} \hat{\mathfrak{f}}_{t,\varepsilon}(r,\tilde{\Phi}(r,\tilde x)\Big] dr
=\int_0^t\vartheta(r)\E\big[\tilde{\mathcal{V}}\tilde{\mathcal{P}}_{t-r}\varphi(\tilde{\Phi}(r,\tilde x)\big]\\[.0cm]
&=\int_0^t \vartheta(r)\tilde{\mathcal{P}}_r\big(\tilde{\mathcal{V}}\tilde{\mathcal{P}}_{t-r}\varphi(\tilde{x})\big) dr.
\end{align*}

Using Fubini's theorem and  Proposition \ref{PertProp}, we have 
\begin{align*}
\lim_{\varepsilon\rightarrow 0}\frac{1}{\varepsilon}\big\langle \tilde{\mathcal{P}}_t^{\varepsilon\vartheta(t)}\varphi-\tilde{\mathcal{P}}_t\varphi, \tilde{\mu}_0\big\rangle&=\int_0^t\vartheta(r)\int_{\Sc\times\Rd}\tilde{\mathcal{P}}_r\big(\tilde{\mathcal{V}}\tilde{\mathcal{P}}_{t-r}\varphi(\tilde{x})\big) \tilde{\mu}_0(d\tilde{x}) dr\\[.2cm]
&=\int_0^t\vartheta(r) \int_{\Sc\times\Rd}\tilde{\mathcal{V}}\tilde{\mathcal{P}}_{t-r}\varphi(\tilde{x})(\tilde{\mathcal{P}}_r^*\tilde{\mu}_0)(d\tilde{x})dr\\[.2cm]
&=\int_0^t\vartheta(r) \int_{\Sc\times\Rd}\tilde{\mathcal{V}}\tilde{\mathcal{P}}_{t-r}\varphi(\tilde{x})\tilde\mu_r(d\tilde{x})dr.
\end{align*}
 By the H\"ormander Lie bracket condition in Assumption \ref{A2.1b}, and ergodicity of the time-periodic measures $\tilde\mu_r$, there exists $0<{\rho}_r\in \mathcal{C}^\infty(\Rd)\cap L^1_+(\Rd)$ such that 
$(\tilde{\mathcal{P}}_r^*\tilde{\mu}_0)(d\tilde{x}) =\tilde\mu_r(d\tilde x)= \tilde{\rho}_r(\tilde{x})d\tilde{x}$, where $\tilde{\rho}_r(\tilde{x})d\tilde{x}$ is understood in the sense of (\ref{abus}).
Thus, for any $\varphi\in \mathcal{C}^2_\tau(\R\times\Rd)$, such that $\E[D_x^\beta\varphi(\tilde\Phi^\alpha(t,\ccdot))]\in L^1(\tilde\mu_0)$, $ |\beta| \leqslant 2$, we have 
\begin{align*}
\lim_{\varepsilon\rightarrow 0}\frac{1}{\varepsilon}\big\langle \tilde{\mathcal{P}}_t^{\varepsilon\vartheta(t)}\varphi-\tilde{\mathcal{P}}_t\varphi, \tilde{\mu}_0\big\rangle&= \int_0^t\vartheta(r) \int_{\Sc\times\Rd}\tilde{\mathcal{V}}\tilde{\mathcal{P}}_{t-r}\varphi(\tilde{x})\tilde{\rho}_r(\tilde{x})d\tilde{x}dr\\[.2cm]
&=\int_0^t\Big(\int_{\Sc\times\Rd}\tilde{\mathcal{P}}_{t-r}\varphi(\tilde{x})\big(\tilde{\mathcal{V}}^*\tilde{\rho}_r\big)(\tilde{x})d\tilde{x}\Big) \vartheta(r)dr.
\end{align*}
Finally, by the definition of the response functional $\mathcal{R}^{\tilde{\mu}_0}_\varphi$  we have for $\vartheta\in \mathcal{C}_{\infty}^1(\Rp,\R)$, that 
\begin{align*}
\hspace{1.5cm} \int_0^t \mathcal{R}^{\tilde{\mu}_0}_\varphi(t-r,r)\vartheta(r) dr& = \lim_{\varepsilon\rightarrow 0}\frac{1}{\varepsilon}\Big[ \mathbb{F}_\varphi^{\tilde{\mu}_0}(t,\varepsilon\vartheta(t))-\mathbb{F}_\varphi^{\tilde{\mu}_0}(t,0)\Big]\\[.2cm]
&=\lim_{\varepsilon\rightarrow 0}\frac{1}{\varepsilon}\big\langle \tilde{\mathcal{P}}_t^{\varepsilon\vartheta(t)}\varphi-\tilde{\mathcal{P}}_t\varphi, \tilde{\mu}_0\big\rangle\\[.2cm]
&=\int_0^t\Big(\int_{\Sc\times\Rd}\tilde{\mathcal{P}}_{t-r}\varphi(\tilde{x})\big(\tilde{\mathcal{V}}^*\tilde{\rho}_r\big)(\tilde{x})d\tilde{x}\Big) \vartheta(r)dr, \hspace{2.3cm}
\end{align*}
where $\tilde{\rho}_r(\tilde{x})d\tilde{x}$ is understood in the sense of (\ref{abus}). \qed

\subsubsection{{\bf Fluctuation-dissipation formulas}}\label{fdt_ss}
Given the general framework for the linear response in the time-periodic regime, we now derive a set of more tractable representations of the response function (\ref{Resp_fun}) via formulas exploiting the  time-asymptotic statistical properties of the unperturbed dynamics (\ref{NSDE11}), or (\ref{PSDE}) with $\alpha=0$; in line with the terminology inherited from statistical physics, these are termed  {\it `fluctuation-dissipation'} \,formulas. The first set of results in Theorem~\ref{FDTI} shadows and formalises formulas derived in \cite{Madja10}, while the results in a more restrictive Theorem~\ref{FDTII} concern the linear response in situations when the `$\alpha$-perturbations' do not destroy the time periodicity and ergodicity of the dynamics in the sense that the coefficients in (\ref{PSDE}) remain $\tau$-periodic for all $\alpha\in \mathcal{A}$. It turns out the two results are related in a specific way. 

\begin{definition}[{\bf Correlation function}] Given the skew-product RDS $\big\{\tilde\Phi(t,\ccdot,\ccdot):\;t\in \Rp\big\}$ on $\Sc\times\Rd$ gicen in (\ref{Lift}), and $\varphi, \psi\in \mathcal{C}_\tau^2(\R\times\Rd)$, such that  for any fixed $\tilde x\in \Sc\times\Rd$, $\E[\varphi\big(\tilde\Phi(t,\tilde x)\big)], \,\E[\psi\big(\tilde\Phi(t,\tilde x)\big)]<\infty$,  the correlation of the random variables $\varphi\big(\tilde\Phi(t,\ccdot,\tilde x)\big)$ and $\psi\big(\tilde\Phi(r,\ccdot,\tilde x)\big)$ is given by  
\begin{align}\label{corr}
 \E\big[\varphi\big(\tilde\Phi(t,\tilde x)\big)\psi\big(\tilde\Phi(r,\tilde x)\big)\big]= \tilde{\mathcal{P}}_{r}\big(\psi(\tilde x)\tilde{\mathcal{P}}_{t-r}\varphi(\tilde x)\big), \quad 0\leqslant  r\leqslant t,
 \end{align}

\noindent where  $(\tilde{\mathcal{P}}_{t})_{t\in \Rp}$ is  defined in (\ref{tldP}), and  (\ref{corr}) follows from the Markov property of the RDS~$\tilde\Phi$. The correlation function based on $\varphi, \psi\in \mathcal{C}^2_\tau(\R\times\Rd)$ and $\tilde\mu\in \PP\big(\Sc\big)\otimes\PP\big(\Rd\big)$ is defined~as
\begin{align}\label{corr_fun}
\mathcal{K}^{\tilde{\mu}}_{\varphi,\psi}(t-r, r):= \int_{\Sc\times\Rd}\tilde{\mathcal{P}}_{r}\big(\psi\tilde{\mathcal{P}}_{t-r}\varphi\big)d\tilde{\mu} = \int_{\Sc\times\Rd}\psi\tilde{\mathcal{P}}_{t-r}\varphi \hspace{.04cm} d(\tilde{\mathcal{P}}^*_{r}\tilde{\mu}), \qquad 0\leqslant r\leqslant t.
\end{align}
\end{definition}

\begin{theorem}[\textbf{FDT I}]\label{FDTI} Suppose that  Assumption \ref{A2.1b} holds,  Assumption \ref{A2.1A} holds for $\alpha=0$, and Proposition \ref{PertProp} is satisfied. Then, for any  $\varphi\in \mathcal{C}^2_\tau\big(\R\times\Rd\big)$ such that $\E[D_x^\beta\varphi(\tilde\Phi^\alpha(t,\!\cdot))]\!\in \!L^1(\tilde\mu_0)$, $ |\beta| \leqslant 2$, the following holds:
\begin{itemize}[leftmargin = 0.8cm]
\item[(i)] There exists a family $(\tilde\mu_t)_{t\geqslant 0}$, ${\tilde\mu}_t\in \PP\big(\Sc\big)\otimes\PP\big(\Rd\big)$, of $\tau$-periodic, skew-product probability measures  and a unique $\tilde{\mathcal{P}}^*_t$-ergodic probability measure,  $ \bar{\tilde\mu}\in \PP\big(\Sc\big)\otimes\PP\big(\Rd\big)$, which is associated with the skew-product RDS $\{\tilde \Phi(t,\ccdot,\ccdot): t\in \Rp\}$ on $\Sc\times\Rd$ and generated by the SDE (\ref{PSDEc2}) with $\alpha=0$.

\vspace{.3cm}\item[(ii)]  The linear response function in (\ref{Resp_fun}) is given~by   
\begin{equation}\label{linres_tper} 
\hspace{.6cm}\mathcal{R}^{\tilde{\mu}_0}_{\varphi}(t-r,r) = \mathcal{K}^{\tilde{\mu}_0}_{\varphi, \mathbb{B}_r}(t-r, r) = \!\!\int_{\Sc\times\Rd} \hspace{-.4cm}\mathbb{B}_r(\tilde x)\tilde{\mathcal{P}}_{t-r}\varphi(\tilde x)\tilde{\rho}_r(\tilde x)d\tilde x, \quad \;\;\mathbb{B}_r(\tilde x) = \frac{\tilde{\mathcal{V}}^*\tilde{\rho}_r(\tilde{x})}{\tilde{\rho}_r(\tilde{x})},
\end{equation}
where $0\leqslant r\leqslant t$,  $\tilde\mu_r(d\tilde x) = \tilde\rho_r(\tilde x)d\tilde x$ is understood in the sense of (\ref{abus}),  the correlation function  $\mathcal{K}^{\tilde{\mu}_0}_{\varphi, \,\scaleobj{.7}{(\boldsymbol{\ccdot})}}$ is defined in~(\ref{corr_fun}), and  the operator $\tilde{\mathcal{V}}^*$ is defined in~(\ref{DUal}).

\vspace{.3cm}\item[(iii)] The linear response for  perturbations of observables based on the ergodic measure $\bar{\tilde \mu}$  is  
\begin{align}\label{linres_inv}
\bar{\mathcal{R}}_{\varphi}(t-r)= \bar{\mathcal{K}}_{\varphi, \mathbb{B}}(t-r) =  \int_{\Sc\times\Rd} \!\!\!\mathbb{B}(\tilde x)\tilde{\mathcal{P}}_{t-r}\varphi(\tilde x)\bar{\tilde{\rho}}(\tilde x)d\tilde x, \qquad \mathbb{B}(\tilde x) = \frac{\tilde{\mathcal{V}}^*\bar{\tilde{\rho}}(\tilde{x})}{\bar{\tilde{\rho}}(\tilde{x})},
\end{align}
where $0\leqslant r\leqslant t$, $ \bar{\mathcal{R}}_{\varphi}(t-r):=\mathcal{R}^{\bar{\tilde{\mu}}}_\varphi(t-r,0)$, and $\bar{\mathcal{K}}_{\varphi, \mathbb{B}}(t-r) := \mathcal{K}^{\bar{\tilde{\mu}}}_{\varphi, \mathbb{B}}(t-r, 0)$.
\end{itemize}
The above hold for all time if Proposition \ref{PertProp} is satisfied for all time.

\end{theorem}
\noindent {\it Proof.} 
Part (i) is a direct consequence of Theorem \ref{Rand_SOL} and Theorem \ref{Ps_erg}.
For Part (ii), we have from the representation  of the response functional $\mathcal{R}^{\tilde{\mu}_0}_{\varphi}$ in (\ref{Resp_fun}) together with the operator $\tilde{\mathcal{V}}^*$ in (\ref{DUal}) that the following holds  for $0\leqslant r\leqslant t$
\begin{align*}
&\mathcal{R}^{\tilde{\mu}_0}_{\varphi}(t-r,r) \\[.1cm] 
&= \int_{\Sc\times\Rd}\tilde{\mathcal{P}}_r\tilde{\mathcal{V}}\big(\tilde{\mathcal{P}}_{t-r}\varphi(\tilde{x})\big)\tilde{\mu}_0(d\tilde{x})= \int_{\Sc\times\Rd}\tilde{\mathcal{V}}\big(\tilde{\mathcal{P}}_{t-r}\varphi(\tilde{x})\big)\tilde{\rho}_r(\tilde{x})d\tilde{x}\\[.2cm]
& = \int_{\Sc\times\Rd}\bigg\{G_i( \tilde{x})\partial_{x_i}\tilde{\mathcal{P}}_{t-r}\varphi(\tilde{x})+\frac{1}{2}\left[\sigma_{ik}(\tilde{x})H_{jk}(\tilde{x})  + \sigma_{jk}(\tilde{x})H_{ik}( \tilde{x})\right]\partial^2_{x_ix_j}\tilde{\mathcal{P}}_{t-r} \varphi(\tilde{x})\bigg\}\tilde{\rho}_{r}(\tilde{x})d\tilde{x}\\[.2cm]
& =\int_{\Sc\times\Rd} \bigg\{-\partial_{x_i}[ G_i(\tilde{x})\tilde{\rho}_r(\tilde{x})]+ \frac{1}{2}\partial^2_{x_ix_j}\Big(\left[ \sigma_{ik}(\tilde{x})H_{jk}(\tilde{x})  + \sigma_{jk}(\tilde{x})H_{ik}(\tilde{x}\right]\tilde{\rho}_r(\tilde{x})\Big)\bigg\}\tilde{\mathcal{P}}_{t-r}\varphi(\tilde{x}) d\tilde{x}\\[.2cm] 
& =\int_{\Sc\times\Rd} \frac{\tilde{\mathcal{V}}^*\tilde{\rho}_r(\tilde{x})}{\tilde{\rho}_r(\tilde{x})}\tilde{\mathcal{P}}_{t-r}\varphi(\tilde{x}) \tilde{\rho}_r(\tilde{x})d\tilde{x} = \int_{\Sc\times\Rd} \mathbb{B}_r(\tilde{x})\tilde{\mathcal{P}}_{t-r}\varphi(\tilde{x})\tilde{\rho}_r(\tilde{x})d\tilde{x} = \mathcal{K}^{\tilde{\mu}_0}_{\varphi, \mathbb{B}_r}(t-r,r),
\end{align*}
where $\varphi\in\mathcal{C}^2_\tau \big(\R\times\Rd\big)$, $\E[D_x^\beta\varphi(\tilde\Phi^\alpha(t,\ccdot))]\in L^1(\tilde\mu_0)$, $ |\beta| \leqslant 2$, as in Theorem~\ref{Linear_resp_formula} (see also Proposition \ref{PertProp} for a sufficient condition), and  $\tilde{\rho}_r(\tilde{x})d\tilde{x}$ is understood in the sense of (\ref{abus}); i.e., 
\begin{align}\label{cresp}
\mathcal{R}^{\tilde{\mu}_0}_{\varphi}(t-r,r) &= \int_{\Sc\times\Rd} \big(\tilde{\mathcal{V}}^*\rho_r(x)\big)\big(\mathcal{P}_{r,t}\varphi(t,x)\big)d{x}\notag\\
&= \int_{\Sc\times\Rd} \mathbb{B}_r(x)\big(\mathcal{P}_{r,t}\varphi(t,x)\big)\rho_r(x)d{x}, \quad 0\leqslant r\leqslant t.
\end{align}
As regards Part (iii), notice  that due to the fact that $\tilde{\mathcal{P}}^*_r\bar{\tilde\mu} = \bar{\tilde\mu}$ for any $r\in \Sc$, we have 
\begin{align*}
\hspace{2cm}\mathcal{R}^{\bar{\tilde{\mu}}}_{\varphi}(t-r,r) &=\int_{\BT}\tilde{\mathcal{P}}_{r}\tilde{\mathcal{V}}\big(\tilde{\mathcal{P}}_{t-r}\varphi\big) \bar{\tilde{\mu}}(d\tilde x) =\int_{\BT}\tilde{\mathcal{V}}\big(\tilde{\mathcal{P}}_{t-r}\varphi\big) \tilde{\mathcal{P}}^*_{r}\bar{\tilde{\mu}}(d\tilde x) \\[.2cm] 
&=\int_{\BT}\tilde{\mathcal{V}}\big(\tilde{\mathcal{P}}_{t-r}\varphi\big) \bar{\tilde{\mu}}(d\tilde x)=\mathcal{R}^{\bar{\tilde{\mu}}}_{\varphi}(t-r,0),
\end{align*}
and the desired result can be derived by following analogous derivations to those above. \hspace{1.cm} \qed

\begin{rem}\label{remark_FDTI}\rm\mbox{}
\begin{itemize}[leftmargin=0.8cm]
\item[(i)] Theorem \ref{FDTI} implies that for the skew-product RDS  $\big\{\tilde\Phi(t,\ccdot, \ccdot)\!:\; t\in \Rp\big\}$ in (\ref{Lift}) induced by the lifted SDE  (\ref{NSDE11}) on $\Sc\times\Rd,$ the change in the value of an observable ${\mathbb{F}}_{\varphi}^{\tilde\mu_0}(t,\alpha)=\langle \tilde{\mathcal{P}}_t \varphi, \tilde{\mu}_0\rangle$ in (\ref{Fobs}) in response to a sufficiently small and regular perturbation  can be represented by the correlation function utilising the unperturbed dynamics/fluctuations. The operator $\tilde{\mathcal{V}}$ defined in (\ref{Va}) does not depend on time due to the $\tau$-periodicity of the coefficients of (\ref{PSDE}) at $\alpha=0$ and the skew-product formulation on $\Sc\times\Rd.$ 
\item[(ii)] The response functions (\ref{linres_tper}) and (\ref{linres_inv})   evaluated on the unperturbed dynamics are amenable to practical approximations via the  appropriate long-time averages. 

By $\tilde{\mathcal{P}}^*_{n\tau}$\,-\,ergodicity of $\tilde \mu_r \in \PP\big(\Sc\big)\otimes\PP\big(\Rd\big)$, $r\in \Sc$,  $r\leqslant t$   we have (see Remark \ref{tau_ergavg})
\begin{align*}
\hspace{.6cm}\mathcal{R}^{\tilde{\mu}_0}_{\varphi}(t-r,r) &= \int_{\Sc\times\Rd} \hspace{-.0cm}\mathbb{B}_r(\tilde x)\tilde{\mathcal{P}}_{t-r}\varphi(\tilde x)\tilde{\rho}_r(\tilde x)d\tilde x \notag\\ 
&=\lim_{n\rightarrow \infty}\frac{1}{N}\sum_{n=0}^N \,\mathbb{B}_r\big(\tilde{\mathcal{P}}_{n\tau}(\tilde x)\big)\tilde{\mathcal{P}}_{t-r}\varphi\big(\tilde{\mathcal{P}}_{n\tau}(\tilde x)\big)\\
&= \lim_{n\rightarrow \infty}\frac{1}{N}\sum_{n=0}^N\,\mathbb{B}_r(\mathcal{P}_{t,t+n\tau}(x))\big(\mathcal{P}_{r,t}\varphi(t,\mathcal{P}_{t,t+n\tau}(x))\big),  
\end{align*}
where $\varphi\in \mathcal{C}^2_\tau\big(\R\times\Rd\big)$, $\E[D_x^\beta\varphi(\tilde\Phi^\alpha(t,\ccdot))]\in L^1(\tilde\mu_0)$, $ |\beta| \leqslant 2$, and in the last line above, we used explicitly the skew-product formulation (\ref{abus}) and the representation (\ref{cresp}); $(\mathcal{P}_{r,t})_{t\geqslant r}$ is the family of transition evolutions defined in (\ref{calP}).

Similarly, by $\tilde{\mathcal{P}}^*_{t}$-\,ergodicity of $\bar{\tilde{\mu}}$ (see Theorem \ref{Ps_erg}) we have for any $u\in \Rp$
\begin{align*}
\hspace{.7cm}\bar{\mathcal{R}}_{\varphi}(u) &= \int_{\Sc\times\Rd} \hspace{-.0cm}\mathbb{B}(\tilde x)\tilde{\mathcal{P}}_{t-r}\varphi(\tilde x)\bar{\tilde{\rho}}(\tilde x)d\tilde x \notag\\
&=\lim_{n\rightarrow \infty}\frac{1}{N}\int_{0}^N \,\mathbb{B}\big(\tilde{\mathcal{P}}_{\zeta}(\tilde x)\big)\tilde{\mathcal{P}}_{u}\varphi\big(\tilde{\mathcal{P}}_{\zeta}(\tilde x)\notag\\
&=\lim_{n\rightarrow \infty}\frac{1}{N}\int_{0}^N \,\mathbb{B}\big({\mathcal{P}}_{0,\zeta}(x)\big){\mathcal{P}}_{0,u}\varphi\big(\zeta,{\mathcal{P}}_{0,\zeta}(x)\big)d\zeta.
\end{align*}
 
\item[(iii)]  Note that the function $\mathbb{B}_r$ in (\ref{Resp_fun}) of Theorem \ref{FDTI} is unique almost everywhere. To see this,  suppose that there exists $\tilde{\mathbb{B}}_r\in L^1(\tilde{\mu}_r)$ such that 
\begin{align*}
\mathcal{K}^{\tilde\mu}_{\varphi, \mathbb{B}_r}(t-r,r) = \mathcal{K}^{\tilde\mu}_{\varphi, \tilde{\mathbb{B}}_r}(t-r,r), \quad \varphi\in \mathcal{C}^2_c(\Sc\times\Rd), \quad 0\leqslant r\leqslant  t.
\end{align*} 
This implies that 
\begin{align*}
\big\langle \tilde{\mathcal{P}}_{t-r}\varphi, \mathbb{B}_r-\tilde{\mathbb{B}}_r\big\rangle_{\tilde\mu_r}=0, \quad \varphi \in \mathcal{C}^2_c(\Sc\times\Rd),\quad 0\leqslant r \leqslant t.
\end{align*}
Given that  $\varphi\in \mathcal{C}^2_c(\Sc\times\Rd)$ is bounded,  taking limit as $r\rightarrow t$ in the above and applying the dominated convergence theorem, we obtain $\mathbb{B}_t=\tilde{\mathbb{B}}_t$\; $\tilde\mu_t$\,-\,a.e.~by arbitrariness of $\varphi$.

\end{itemize}
\end{rem}

Throughout the reminder of this section, we shall assume that $t\mapsto \hat b(\alpha,t,x), \,\hat \sigma(\alpha,t,x)$, are $\tau$-periodic for all $(\alpha,x)\in \Ac\times\Rd$.  
Thus, under Assumptions \ref{A2.1A} and \ref{A2.1b}, the family of (skew-product) measures $\{\tilde{\mu}^\alpha_t: t\in \Rp, \; \alpha\in \Ac\}$ is $\tau$-periodic and ergodic (as in  the case of \S\ref{permeas}) and, as a consequence, the time-averaged ($\tilde{\mathcal{P}}_t^{\alpha*}$-\,ergodic) measures  $\bar{\tilde{\mu}}^\alpha$ (see (\ref{ERm})) satisfy 
\begin{align}\label{StaFKP}
\int_{\Sc\times\Rd}\tilde{\LG}_\alpha\varphi(\tilde x) \bar{\tilde{\mu}}^\alpha(d\tilde x) = 0, \qquad \alpha\in \Ac, \;\;\varphi \in \D(\tilde{\LG}_\alpha)\cap\mathcal{E}_\alpha,
\end{align}
where $\tilde{\LG}_\alpha$ is the generator of the RDS  $\{\tilde \Phi^\alpha(t,\ccdot,\ccdot)\!: \,t\in \Rp\}$ on $\Sc\times\Rd$, for $\alpha\in \Ac$ (in the same form as $\tilde{\mathcal{L}}$ in Definition \ref{inf_gnrt}),  and the domain $\D(\tilde{\LG}_\alpha)$ is defined by
  \begin{align}\label{Dalph}
\mathcal{D}(\tilde{\LG}_\alpha):=&\mathbb{H}^1\big(\Sc\times\Rd; \bar{\tilde{\mu}}^\alpha\big)\notag\\[.2cm] 
=&\Big\{\varphi: \Sc\times\Rd\rightarrow\R:\; \varphi(0,\ccdot) = \varphi(\tau, \ccdot)\; \textrm{and}\; \notag\\[.2cm] 
&\quad \int_{\Sc\times\Rd} \vert \varphi(s,x)\vert^2\bar{\tilde{\mu}}^\alpha(dsdx) + \int_{\Sc\times\Rd} \vert D\varphi(s,x)\vert^2\bar{\tilde{\mu}}^\alpha(dsdx)<\infty, \;\alpha\in \mathcal{A}\Big\} .
 \end{align}
 and 
 \begin{equation}\label{Etau}
 \mathcal{E}_\alpha:=\Big\{\varphi\in \mathcal{C}^2(\Sc\times\Rd):\; \varphi(0,\ccdot) = \varphi(\tau, \ccdot),\; \E[D_x^\beta\varphi(\tilde\Phi^\alpha(t,\ccdot))]\in L^1(\bar{\tilde{\mu}}^\alpha), |\beta| \leqslant 2,\;\alpha\in \mathcal{A}\Big\}.
 \end{equation} 
 Sufficient conditions for $\varphi\in \mathcal{E}_\alpha$ were given in Proposition \ref{PertProp}.

\medskip

Given the ergodic measure $\bar{\tilde\mu}$ (\ref{ERm}) of the unperturbed dynamics (i.e., (\ref{PSDEc2}) with $\alpha=0$) and the above setup, the linear response $\Delta \mathbb{F}_{\varphi,\vartheta}^{\bar{\tilde\mu}} (t)$ of the statistical observable ${\mathbb{F}}_{\varphi}^{\bar{\tilde\mu}}(t,\alpha)$ in (\ref{Fobs}) due to a sufficiently small perturbation $\varepsilon \vartheta(t)$ around $\alpha=0$  such that it preserves the $\tau$-periodicity of the unperturbed dynamics is summarised as follows. 

\begin{theorem}[FDT II]\label{FDTII}
Let Assumptions~\ref{A2.1b}-\ref{A2.1A} be satisfied for $\alpha\in \Ac$,  and suppose  that Proposition \ref{PertProp} holds. Assume further that $t\mapsto \hat b(t,\alpha,x)$, $t\mapsto\hat \sigma(t,\alpha,x)$ are $\tau$-periodic for all $(\alpha,x)\in \Ac\times\Rd$. Then,  the following hold:
 \begin{itemize}[leftmargin = 0.8cm]
 \item[(i)]  For every $\alpha\in \Ac$,    $\{\tilde{\mu}_t^\alpha: t\in \Rp\}$ is a family of $\tau$-periodic skew-product measures  induced by $\tilde\Phi^\alpha$ in (\ref{Lift_alph}),  with the uniquely $\tilde{\mathcal{P}}_t^{\alpha*}$-\,ergodic  measure $\bar{\tilde{\mu}}^\alpha$ satisfying (\ref{StaFKP}). 
 
\vspace{.2cm}\item[(ii)]  The map $\alpha\mapsto \bar{\tilde{\mu}}^\alpha(d\tilde{x}) = \bar{\tilde{\rho}}^{\alpha}(\tilde{x})d\tilde{x}$ (with $\bar{\tilde{\rho}}^{\alpha}$ understood in the sense of (\ref{abus})) is weakly differentiable at $\alpha =0$ for all $\tilde{x}\in \Sc\times\Rd$, and the linear response  (\ref{Resp_fun}) associated with perturbations of observables based on the  $\tilde{\mathcal{P}}_t^{*}$-\,ergodic skew-product measure $\bar{\tilde\mu}$ is given by 
\begin{align}\label{linres_invII}
\bar{\mathcal{R}}_{\varphi}(t-r):= {\mathcal{R}}^{\bar{\tilde\mu}}_{\varphi}(t-r,0) = \partial_r\mathcal{K}^{\bar{\tilde\mu}}_{\varphi, \mathbb{W}}(t-r, r),\qquad0 \leqslant r\leqslant t,
\end{align}
for $\varphi\in \D(\tilde{\LG}_\alpha)\cap\mathcal{E}_\alpha$, with the correlation function  $\mathcal{K}^{\bar{\tilde\mu}}_{\varphi, \mathbb{W}}$ in (\ref{corr_fun}) with   $\mathbb{W}\in\mathcal{C}_\infty\big(\BT\big)$ given by 
\begin{align}\label{WW}
\mathbb{W}(\tilde{x}) = \frac{\eta(\tilde{x})}{\bar{\tilde{\rho}}(\tilde{x})}, \quad \textrm{s.t.}\quad \langle\eta, \varphi\rangle := \langle \partial_\alpha\bar{\tilde{\rho}}^\alpha, \varphi\rangle\big|_{\alpha=0},
\end{align} 
with $\partial_\alpha\bar{\tilde{\rho}}^\alpha$ understood in the weak sense.

\end{itemize}
\end{theorem}

\noindent {\it Proof.} Part (i) is a direct consequence of Theorem \ref{Rand_SOL} and Theorem \ref{Ps_erg} given the fact that  Assumptions~\ref{A2.1A}--\ref{A2.1b} hold for $\alpha$ in a proper interval $\mathcal{A}$ containing $\alpha=0$. For Part (ii), we proceed as at the beginning of the the proof of Theorem \ref{Linear_resp_formula}, except that due to Part (i), both the unperturbed and the perturbed measures are $\tau$-periodic. Thus,  for $\vartheta \in \mathcal{C}^1_\infty(\Rp;\R)$ and  $\varepsilon>0$ sufficiently small so that $\varepsilon \vartheta\in \Ac$, and  for all $\varphi\in \mathcal{E}_\alpha$  there exists  $0<C_{\varepsilon,\varphi}<\infty$ such that 
 \begin{align}\label{qq1}
 \langle\varphi, \tilde{\mu}_t^{\varepsilon\vartheta}\rangle- \langle\varphi, \tilde{\mu}_t\rangle &=\int_{\BT} \left(\tilde{\mathcal{P}}^{\varepsilon\vartheta}_t\varphi(\tilde{x})- \tilde{\mathcal{P}}_t\varphi(\tilde{x})\right) \tilde\mu_0(d\tilde{x})\notag\\[.1cm] 
&=\int_{\BT} \left(\int_0^t\varepsilon \vartheta(r)\E\big[ \,\hat{\mathfrak{f}}_{t,\varepsilon}(r,\tilde\Phi(r,\ccdot,\tilde x))\big]dr\right) \tilde{\mu}_0(d\tilde{x})\notag\\[.2cm]
 &\leqslant \varepsilon t\Vert\vartheta\Vert_{\infty}\,C_{\varepsilon,\varphi}, 
\end{align}
where the bound is due to Proposition \ref{PertProp}. Averaging both sides over $t\in \Sc$ we have 
\begin{align}\label{qq2}
 \lim_{\varepsilon\downarrow 0} \frac{1}{\varepsilon} \langle\varphi, \bar{\tilde{\mu}}^{\varepsilon\vartheta} -\bar{\tilde{\mu}}\rangle <\infty \quad \forall \;\varphi\in \mathcal{E}_\alpha. 
\end{align}
Thus, $\bar{\tilde\mu}^\alpha$ is weakly differentiable on $\mathcal{E}_\alpha$.  Next,  by the H\"ormander condition in Assumption~\ref{A2.1b}, we have  $\bar{\tilde{\mu}}^\alpha(d\tilde{x}) = \bar{\tilde{\rho}}^\alpha(\tilde{x})d\tilde{x}$ so that  for any $\vartheta \in \mathcal{C}^1_\infty\big(\Rp,\R\big)$
\begin{align}\label{qq3}
 \lim_{\varepsilon\downarrow 0} \frac{1}{\varepsilon} \langle\varphi, \bar{\tilde{\rho}}^{\varepsilon\vartheta} -\bar{\tilde{\rho}}\rangle 
<\infty \quad \forall\;\varphi\in \mathcal{E}_\alpha,
\end{align}
  which implies that $\bar{\tilde{\rho}}^{\alpha}$ is weakly differentiable at $\alpha=0$ (with $\bar{\tilde{\rho}}^{\alpha}$ understood in the sense of (\ref{abus}) to simplify notation).
 Furthermore, (\ref{StaFKP}) yields
\begin{align}\label{SEqd}
\big\langle \bar{\tilde{\rho}}^\alpha, \tilde{\LG}_\alpha\varphi\big\rangle =0, \qquad \alpha\in \Ac, \;\;\varphi \in \D(\tilde{\LG}_\alpha)\cap\mathcal{E}_\alpha.
\end{align}
Differentiating (\ref{SEqd}) with respect to the parameter $\alpha$ (in the weak sense), we obtain  
\begin{align}\label{SEqd1}
\big\langle\partial_\alpha \bar{\tilde{\rho}}^\alpha,-\tilde{\LG}\varphi\big\rangle\Big|_{\alpha=0} = \big\langle \,\bar{\tilde{\rho}},\,\tilde{\mathcal{V}}\varphi \,\big\rangle=\big\langle \tilde{\mathcal{V}}^*\bar{\tilde{\rho}}, \varphi\big\rangle.
\end{align}
Next, we set $\big\langle \eta, -\tilde{\LG}\varphi\big\rangle:=\langle \partial_\alpha\bar{\tilde{\rho}}^\alpha,-\tilde{\LG}\varphi\rangle|_{\alpha=0}$, and note that $\tilde{\mathcal{P}}_t\varphi \in \D(\tilde{\LG}_\alpha)\cap\mathcal{E}_\alpha$, $t\geqslant 0$,  for any $\varphi \in D(\tilde{\LG}_\alpha)\cap\mathcal{E}_\alpha$  (due to Assumption \ref{A2.1b} and the associated smoothing property of $(\mathcal{P}_{s,t})_{s\leqslant t}$ generating $(\tilde{\mathcal{P}}_{t})_{t\in \Rp}$ in (\ref{tldP}); e.g., Proposition \ref{StrongF} and \cite{Da Prato2}).  Thus, we have 
 \begin{align}\label{KWW}
\mathcal{K}^{\bar{\tilde\mu}}_{\varphi,\mathbb{W}}(t-r,r) = \left\langle \tilde{\mathcal{P}}_r\left(\mathbb{W}\tilde{\mathcal{P}}_{t-r}\varphi\right), \bar{\tilde{\rho}}\right\rangle=\left\langle \tilde{\mathcal{P}}_{t-r}\varphi, \mathbb{W}\bar{\tilde{\rho}}\right\rangle = \left\langle \tilde{\mathcal{P}}_{t-r}\varphi, \eta\right\rangle,
\end{align}
where $\mathbb{W}$ is given in (\ref{WW}), and subsequently 
\begin{align*}
\partial_r\mathcal{K}_{\varphi, \mathbb{W}}(t-r,r) =\partial_r\langle \tilde{\mathcal{P}}_{t-r}\varphi, \eta\rangle = \langle -\tilde{\LG}\tilde{\mathcal{P}}_{t-r}\varphi, \eta\rangle.
\end{align*}
Since $\tilde{\mathcal{P}}_{t-r}\varphi\in \D(\tilde{\mathcal{L}}_\alpha)\cap\mathcal{E}_\alpha$ for $0\leqslant r\leqslant t$, we have by (\ref{SEqd1}) that for $\varphi \in \D(\tilde{\LG}_\alpha)\cap\mathcal{E}_\alpha$
\begin{align*}
\hspace{1cm}\partial_r\mathcal{K}_{\varphi, \mathbb{W}}(t-r,r) = \langle -\tilde{\LG}\tilde{\mathcal{P}}_{t-r}\varphi, \eta\rangle & = \big\langle \tilde{\mathcal{V}}\tilde{\mathcal{P}}_{t-r}\varphi, \bar{\tilde{\rho}}\,\big\rangle\\[.2cm]
&=\int_{\BT}(\tilde{\mathcal{V}}\tilde{\mathcal{P}}_{t-r}\varphi)(\tilde{x})\bar{\tilde{\mu}}(d\tilde{x})= {\mathcal{R}}^{\bar{\tilde\mu}}_{\varphi}(t-r, 0). \hspace{1.5cm}\qed
\end{align*}

\begin{rem}\label{Remark2}\label{entropy_rem}\rm
Note that Theorem~\ref{FDTI} is more general than Theorem~\ref{FDTII} in the sense that  it only requires time-periodicity of the coefficients of the SDE (\ref{PSDE}) and the existence of time-periodic probability measure for the unperturbed dynamics (i.e., for $\alpha=0$ in (\ref{PSDE})) but it does not preclude the perturbed dynamics to have time-periodic measures. Thus, a natural question arises as to the connection between the linear response functions $\bar{\mathcal{R}}_{\varphi}$  in (\ref{linres_inv}) of Theorem~\ref{FDTI} and (\ref{linres_invII}) of Theorem \ref{FDTII}, respectively, in the case when both the unperturbed and the perturbed dynamics (i.e., for $\alpha\in \mathcal{A}$ in (\ref{PSDE})) have time-periodic ergodic measures of period $\tau$.  The desired connection stems from the fact that under Assumption \ref{A2.1b} the identity (\ref{SEqd1}) leads to 
 \begin{align*}
\langle \mathbb{W},-\tilde{\LG}\varphi\rangle_{\bar{\tilde{\mu}}}= \langle \mathbb{W}\bar{\tilde\rho} , -\tilde{\LG}\varphi\rangle   = \langle \eta, -\tilde{\LG}\varphi\rangle = \langle \tilde{\mathcal{V}}^*\bar{\tilde{\rho}}, \varphi\rangle = \Big\langle \frac{\tilde{\mathcal{V}}^*\bar{\tilde{\rho}}}{\bar{\tilde{\rho}}}, \bar{\tilde{\rho}}\hspace{.04cm}\varphi\Big\rangle
= \langle \mathbb{B}, \varphi\rangle_{\bar{\tilde{\mu}}},
\end{align*}
for any $\varphi\in \D(\tilde{\LG}_\alpha)\cap\mathcal{E}_\alpha$,  where $\tilde{\LG}=\partial_s+{\LG}$ is the generator of the one-point motion $\tilde{x}\mapsto\tilde\Phi(t,\om,\tilde{x})$, $t\geqslant 0$,  on the extended state space $\Sc\times\Rd$. Thus, we obtain  $\mathbb{B} = -\tilde{\LG}^{\bar{\tilde{\mu}}*}\mathbb{W}$, where $\tilde{\LG}^{\bar{\tilde{\mu}}*}$ is the $L^2(\bar{\tilde{\mu}})$ dual of $\tilde{\LG}$.  In fact, the above result also implies that $\mathbb{W}$ in (\ref{WW}) of Theorem~\ref{FDTII} is not unique  in contrast to $\mathbb{B}$ in (\ref{linres_inv}) of Theorem~\ref{FDTI}. To see this,  assume that there exists $\tilde{\mathbb{W}}\in L^1(\bar{\tilde{\mu}})$ on $\Sc\times\Rd$ such that 
\begin{align*}
\partial_r\mathcal{K}^{{\bar{\tilde\mu}}}_{\varphi, \mathbb{W}}(t-r,r) = \partial_r\mathcal{K}^{{\bar{\tilde\mu}}}_{\varphi, \tilde{\mathbb{W}}}(t-r,r) \qquad \forall\;\varphi\in \mathcal{C}_c^{\infty}(\Sc\times\Rd),\; \;0\leqslant r\leqslant t, 
\end{align*}
which implies that (see (\ref{KWW}))
\begin{align*}
\frac{d}{dr}\big\langle \tilde{\mathcal{P}}_{t-r}\varphi, \mathbb{W}-\tilde{\mathbb{W}}\big\rangle_{\bar{\tilde\mu}}=0 \qquad \forall \;\varphi\in \mathcal{C}_c^{\infty}(\Sc\times\Rd),\; \;0\leqslant r\leqslant t.
\end{align*}
Since $\mathcal{C}_c^{\infty}(\Sc\times\Rd)\subset \D(\tilde{\LG})$, it follows by the smoothing property of $\tilde{\mathcal{P}}_{t-r}$ (under Assumption~\ref{A2.1b})  and the dominated convergence theorem that 
\begin{align*}
\big\langle \tilde{\LG}\varphi, \mathbb{W}-\tilde{\mathbb{W}}\big\rangle_{\bar{\tilde\mu}} =0.
\end{align*}
Since $\mathbb{W}-\tilde{\mathbb{W}}\in L^1(\bar{\tilde{\mu}})$,  $\mathbb{W}-\tilde{\mathbb{W}}$ is a.e.~constant; hence, $\mathbb{W}$ satisfying Theorem \ref{FDTII}  is not  unique. 
\end{rem}

\begin{exa}[Stochastic Lorenz model with periodic forcing]\label{ex_L_fdt}\rm
We return to the stochastic Lorenz model with time-periodic forcing used in Example~\ref{Lor_ex}, and we consider a simple case of perturbed dynamics in (\ref{Lorenz}) in the form 
\begin{align}\label{perLor}
dv_t^\alpha = \big[-Av_t^\alpha-G(v^\alpha_t) +F(t)+\varepsilon \hspace{0.04cm}\mathfrak{F}(v^\alpha_t)\vartheta(t)\big]dt+\sigma(v_t^\alpha)dW_t,
\end{align} 
where $\vartheta\in \mathcal{C}_{\infty}^1(\Rp,\R)$, $\vartheta(0)=0$, $v\mapsto \mathfrak{F}(v)\in {\mathcal{C}}^\infty_\infty(\Rd)$, and $A,G(v),F(t)$ are defined as in (\ref{Lorenz}) for $v_t^\alpha = (x_t^\alpha,y_t^\alpha,z_t^\alpha)$. Here, the last term in the drift represents a perturbation of the dynamics~(\ref{Lorenz})  due to $\alpha(t) = \varepsilon \vartheta(t)$, which is chosen for simplicity to be in the space-time factorised form (recall, however, that the perturbation in this framework can take a more general form; see Remark~\ref{pert_fact_rem}).

Note that, analogously to (\ref{Lorenz}), the periodically forced dynamics in (\ref{perLor}) satisfies  Assumption~\ref{A2.1b}, and recall that  in appropriate parameter regimes of (\ref{Lorenz}) (or in (\ref{perLor}) with  $\alpha=0$) there  exists a time-periodic ergodic measure $\mu_t\in \PP(\Rd)$ with a smooth density $\rho_t$ with respect to the Lebesgue measure on $\R^3$ for all $t\in \Sc, \,0<\tau<\infty$.

\begin{figure}[t]
\captionsetup{width=1\linewidth}
\centering

\vspace*{-.2cm}\includegraphics[width = 16.2cm]{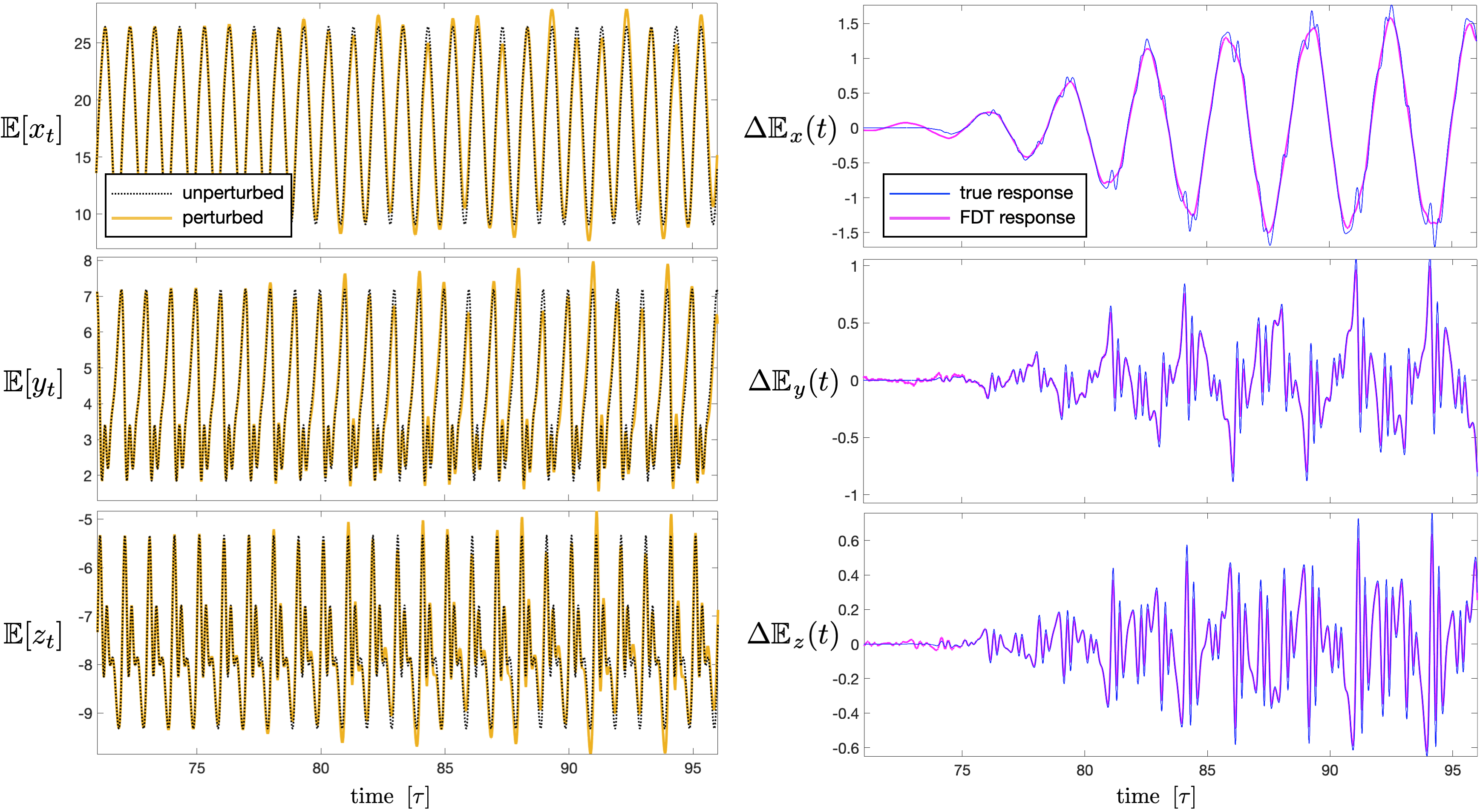}

\vspace{-0.2cm}\caption{\footnotesize The responses $\Delta\mathbb{E}_x(t)$, $\Delta\mathbb{E}_y(t)$, $\Delta\mathbb{E}_z(t)$  of the expectations $\mathbb{E}[x_t]$, $\mathbb{E}[y_t]$, $\mathbb{E}[z_t]$, associated with the dynamics (\ref{Lorenz}) perturbed to~(\ref{perLor}) via the  spatially uniform, time-aperiodic perturbation (\ref{ex_tht}) with $\varepsilon = 0.05$, and $t_0 = 80\tau, \Delta T = 12$. The left column shows the direct simulation of (\ref{perLor}), and the right column shows the  comparison between the exact response and the linear response given by (\ref{LRR})  and the response function given by ~(\ref{manR}).  The parameters of the unperturbed system (\ref{Lorenz}) are $\bar\alpha = 7.3, \bar\beta = 26, \bar\gamma= 7, \bar\varrho = 10$, $\bar f = 100, \bar\delta=0.9, \tau = 1, \bar\sigma=0.2$; see Figure \ref{lor_fig}. Further details are discussed in  the main text of Example \ref{ex_L_fdt}.}\label{fdt_1}
\end{figure}

Assuming that the conditions of Theorem \ref{Linear_resp_formula} hold, the linear response $\Delta \mathbb{F}_{\varphi,\vartheta}^{\tilde\mu_0} (t)$ in  (\ref{LinR})
\begin{align}\label{LRR}
\Delta \mathbb{F}_{\varphi,\vartheta}^{\tilde\mu_0} (t)=\int_0^t\mathcal{R}_{\varphi}^{\tilde\mu_0}(t-r, r)  \vartheta(r)dr, \qquad \vartheta\in \mathcal{C}^1_\infty(\Rp,\R), \quad \vartheta(0)=0,
\end{align}
of  the observable $\varphi$ to the  perturbation $\varepsilon\hspace{0.04cm}\mathfrak{F}(v)\vartheta(t)$ in (\ref{perLor})  in the  random time-periodic ergodic regime is determined by convolving the response function $\mathcal{R}^{\tilde\mu_0}_{\varphi}$ with~$\vartheta$,  where  (see~Theorem~\ref{FDTI})
\begin{align}\label{ex_corr_nper}
\mathcal{R}_{\varphi}^{\tilde\mu_0}(t-r,r) = \int_{\Sc\times\R^3}\partial_v\big(\mathfrak{F}(v)\tilde{\rho}_r(\tilde{v})\big)\tilde{\mathcal{P}}_{t-r}\varphi(\tilde{v})d\tilde v = \mathcal{K}_{\varphi,\mathbb{B}_r}^{\tilde\mu_0}(t-r,r),  \quad0\leqslant r\leqslant t, 
\end{align}
which is solely based on the unperturbed dynamics. In the above expression $\tilde\mu_t(\tilde v) = \tilde\rho_t(\tilde v)d\tilde v$, $\tilde{v} = (s, v)\in [0, \,\tau)\times\R^3$, and $\mathcal{K}_{\varphi,\mathbb{B}_r}^{\tilde\mu_0}(t-r,r)$ is the correlation function (\ref{corr_fun}) of the random variables $\varphi(\tilde{v}_t)$, and $\mathbb{B}_r(\tilde{v}) = -\partial_v(\mathfrak{F}(v)\tilde{\rho}_r(\tilde{v}))/\tilde{\rho}_r(\tilde{v})$, with  $\tilde\mu_t(\tilde v) = \tilde\rho_t(\tilde v)d\tilde v$ understood in the sense of (\ref{abus}). Thus, the response function can be written in a form amenable to computations  as 
\begin{align}\label{manR}
\mathcal{R}_{\varphi}^{\tilde\mu_0}(t-r,r) &=\int_{\Sc\times\Rd} \big(\tilde{\mathcal{V}}^*\rho_r(v)\big)\big(\mathcal{P}_{r,t}\varphi(t,v)\big)d{v}\notag\\[.2cm]
&= \int_{\Sc\times\Rd} \mathbb{B}_r(v)\big(\mathcal{P}_{r,t}\varphi(t,v)\big)\rho_r(v)d{v},  \quad0\leqslant r\leqslant t, 
\end{align}
(see Remark \ref{remark_FDTI}(ii))  where $\tilde{\mathcal{V}}^*$ is defined in (\ref{DUal}), $\mathbb{B}_r$ is given in (\ref{linres_tper}), and $(\mathcal{P}_{r,t})_{t\geqslant r}$ is the family of transition evolutions defined in (\ref{calP}). Furthermore,  (\ref{manR}) can be evaluated in a more practical fashion via  appropriate  ergodic averages, as discussed in Remark~\ref{remark_FDTI}(ii).

For simplicity of the numerical illustration we consider the linear response of the expectation of the solutions to (\ref{Lorenz}) to  the perturbation $\varepsilon\hspace{0.04cm}\mathfrak{F}(v^\alpha)\vartheta(t)$ introduced in the drift coefficient of~(\ref{perLor}). Given the dynamics (\ref{perLor}) and  $\varphi(v)=v$, setting  $p=2$ is sufficient for Assumption \ref{A2.1} and Proposition \ref{PertProp} to hold (i.e., Theorems \ref{Rand_SOL}, \ref{Ps_erg}, \ref{Linear_resp_formula}, \ref{FDTI}  will hold for all time if the perturbation maintains dissipativity, which is the case here). In the examples shown in Figures \ref{fdt_1}, \ref{fdt_2} we denote the expectation of the solutions to  (\ref{perLor})  by  $\mathbb{E}[x_t]$, $\mathbb{E}[y_t]$, $\mathbb{E}[z_t]$, and we consider the response of the expectation  to  a spatially uniform perturbation $\varepsilon\hspace{0.04cm}\mathfrak{F}(v^\alpha)\vartheta(t)$ with $\mathfrak{F}(v^\alpha)=(\bar f,\,0,\,0)^T$ and 
\begin{equation}\label{ex_tht}
\vartheta(t) = \varTheta(t)\big(1+\cos(2\pi/3.3 \,t)\big),
\end{equation} 
where 
\begin{equation}\label{ex_THT}
\varTheta(t) = 
\begin{cases}
\qquad \;\; \;1 & t>t_0+\Delta T,\\
\frac{2}{\Delta T}(t-t_0)-\frac{1}{{\Delta T}^2}(t-t_0)^2 \quad & t_0\leqslant t \leqslant t_0+\Delta T,\\
\frac{2}{\Delta T}(t-t_0)+\frac{1}{{\Delta T}^2}(t-t_0)^2 & t_0-\Delta T\leqslant t \leqslant t_0,\\
\qquad -1 & t<- t_0-\Delta T.
\end{cases}
\end{equation}
 The simulations were performed for (\ref{perLor}) with the same  parameter values as those in Example~\ref{Lor_ex} in the random time-periodic regime,  and the perturbation with $t_0 = 80\tau, \Delta T = 1$ with  the amplitudes set to $\varepsilon = 0.05$ in Figure~\ref{fdt_1}, and $\varepsilon = 0.25$ in Figure \ref{fdt_2}.  The unperturbed initial  time-periodic measure at $t=0$, i.e., $\tilde \mu_0 = \delta_0\otimes \mu_0$, was  approximated from  long-time simulations of an ensemble of solutions to (\ref{perLor}) with $\alpha=0$ at $t = \tau n$, $n\in \mathbb{N}_0$. To simplify the notation, the linear response $\Delta \mathbb{F}_{\varphi,\vartheta}^{\tilde\mu_0} (t)$ in (\ref{LinR}) of the expectations is denoted by, respectively,  $\Delta\E_{\hspace{.03cm}x}(t)$,  $\Delta\E_{\hspace{.03cm}y}(t)$,  $\Delta\E_{\hspace{.03cm} z}(t)$. The linear response was estimated with the help of the fluctuation-dissipation formula (\ref{manR}), where $\mathcal{K}_{\varphi,\mathbb{B}_r}(t-r,r)$ in (\ref{corr_fun}) exploits the statistical correlations in the time-asymptotic dynamics of the unperturbed system (\ref{Lorenz}) via (\ref{corr}). As expected from the theory, the linear response provides a good approximation for a sufficiently small perturbation (Figure \ref{fdt_1}), and it deteriorates with the increasing amplitude of the perturbation (Figure \ref{fdt_2}); the accuracy of the approximation improves still with the decreasing amplitude of the perturbation but we do not show these unsurprising  results. 

\begin{figure}[t]
\captionsetup{width=1\linewidth}
\centering

\vspace{-.2cm}\includegraphics[width = 16.2cm]{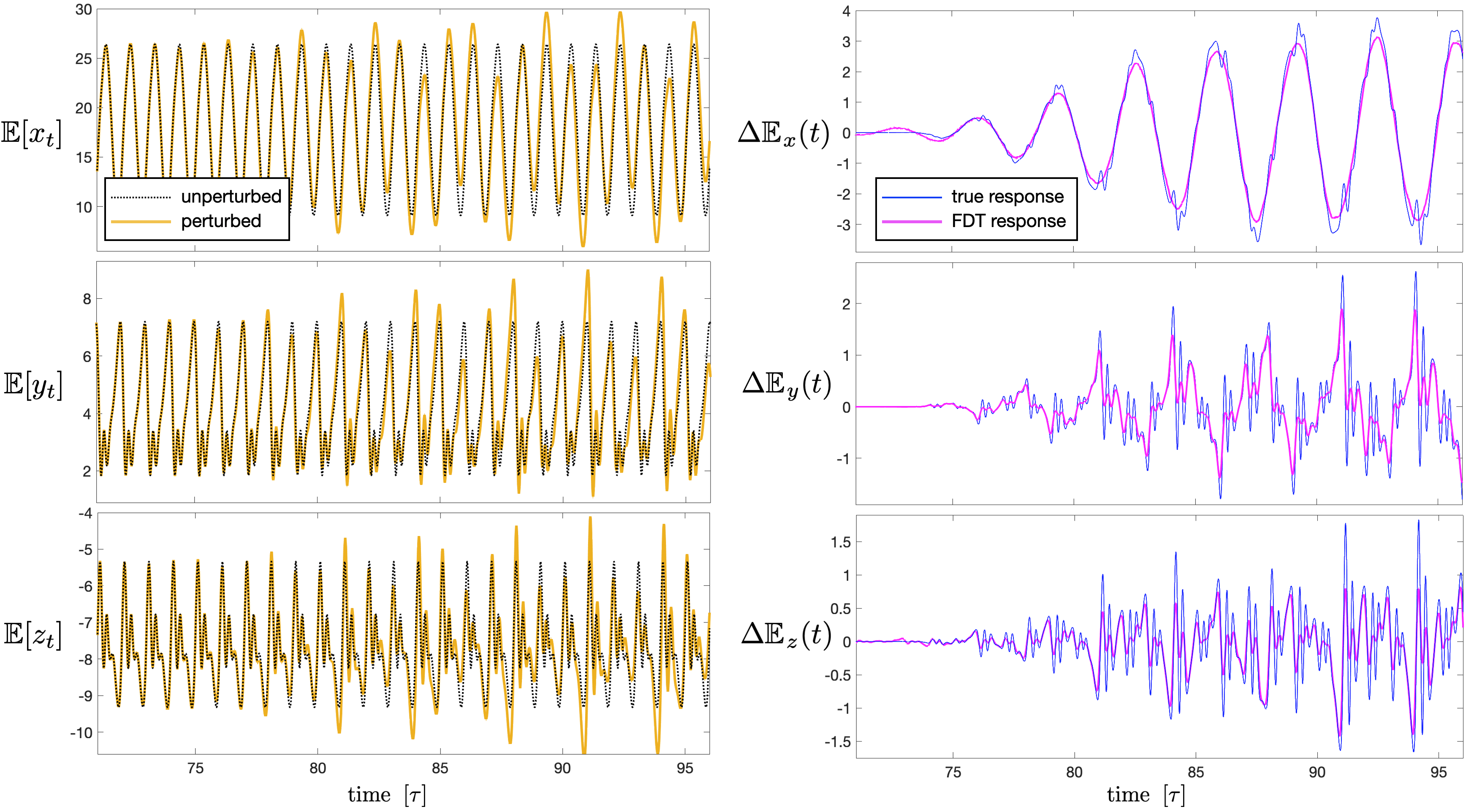}

\vspace{-.2cm}\caption{\footnotesize The responses $\Delta\mathbb{E}_x(t)$, $\Delta\mathbb{E}_y(t)$, $\Delta\mathbb{E}_z(t)$ of $\mathbb{E}[x_t]$, $\mathbb{E}[y_t]$, $\mathbb{E}[z_t]$ associated with the dynamics (\ref{Lorenz}) to a spatially uniform, time-aperiodic perturbation (\ref{ex_tht}) in  (\ref{perLor}) with $\varepsilon = 0.25$; the remaining parameters are as in Figure~\ref{fdt_1}. See the main text in Example \ref{ex_L_fdt} for more details.}\label{fdt_2}
\end{figure}

In the simple example illustrated in Figure \ref{fdt_3}, we consider the response of the expectation of the solution to (\ref{perLor}) in the stable random time-periodic regime to a spatially uniform time-periodic perturbation $\varepsilon\hspace{0.04cm}\mathfrak{F}(v^\alpha)\vartheta(t)$ with $\varepsilon = 0.1$,  $\mathfrak{F}(v^\alpha)=(\bar f,\,0,\,0)^T$   and  
\begin{equation}\label{ex_tht3}
\vartheta(t) = H(t-80.25\tau)\cos^2(2\pi/\tau \,t),
\end{equation}
where $H(t)$ is the Heaviside step function. As in the previous examples, the accuracy of the linear response via the FDT formulas (\ref{ex_corr_nper}) or (\ref{manR}) combined with (\ref{LRR}) improves for decreasing magnitude  of the perturbation but we do not show these unsurprising results. Note that in this case the perturbed measures are also $\tau$-periodic and one could consider the FDT  for the ergodic measure $\bar{\tilde\mu}$ as in Theorem~\ref{FDTII}; such considerations of the linear response are more relevant in the abstract analysis and we do not pursue such a scenario  here (see, however, Remark \ref{entropy_rem}).  

Finally, it needs to be stressed that the direct numerical evaluation of the correlation function $\mathcal{K}_{\varphi,\mathbb{B}_r}(t-r,r)$ in (\ref{ex_corr_nper}) is, in general, very computationally intensive due to the need for estimating the time-dependent density $\rho_t(v)$, $v\in \R^3$, $t\in \Sc$; more practical implementations rely on various approximations (e.g., a Gaussian approximation of the underlying density), and they  were discussed in \cite{Madja10} in the time-periodic setting, and in \cite{ Abramov07, Abramov08,Abramov09, Majda05, Madja10, Madja10b,grit99,grit02,grit07,grit08,majdaqi19} in the stationary setting. These references also consider much more elaborate examples than what we could consider in this  work.
\end{exa}

\begin{figure}
\captionsetup{width=1\linewidth}
\centering
\includegraphics[width = 16.2cm]{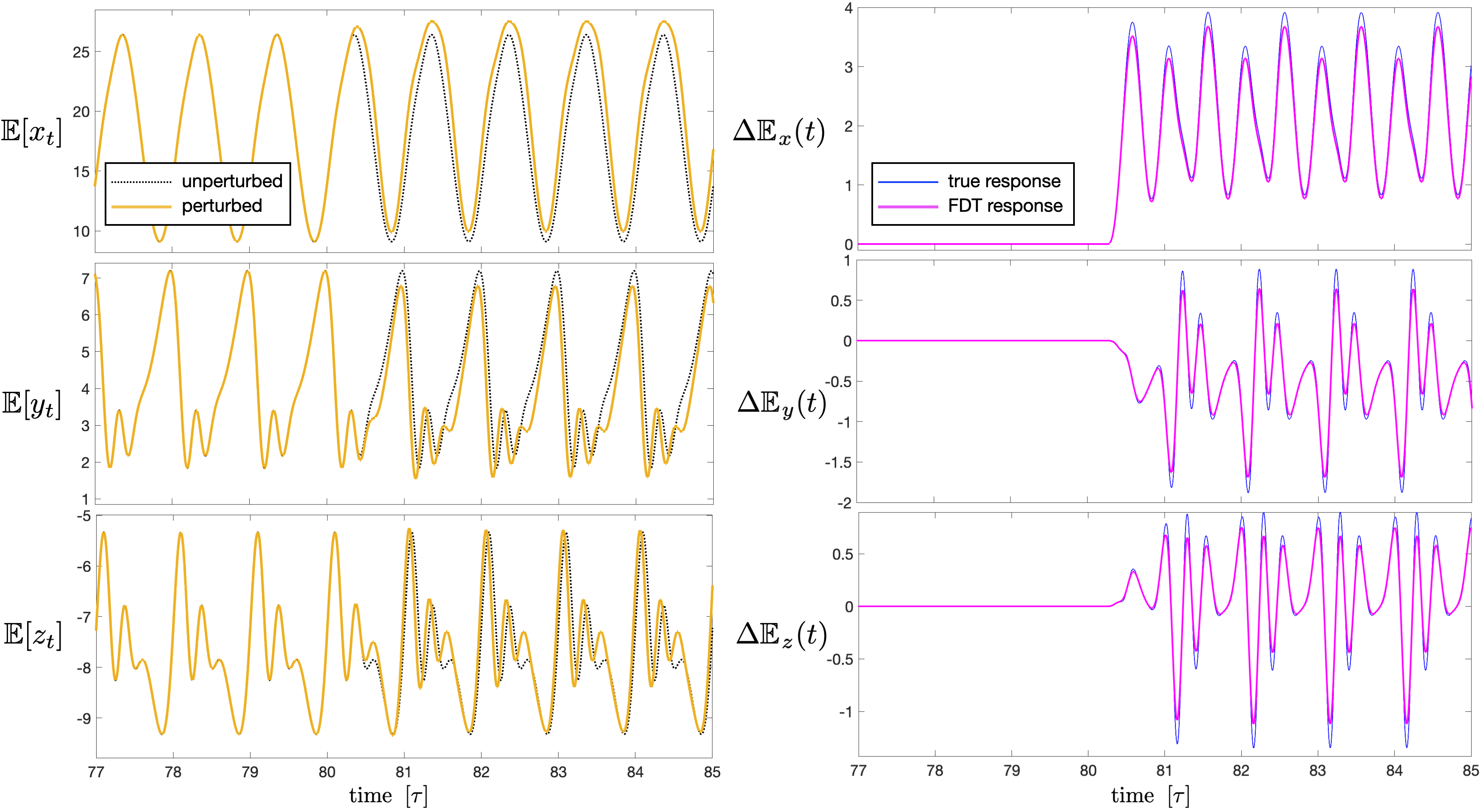}

\vspace{-.2cm}\caption{\footnotesize The responses $\Delta\mathbb{E}_x(t)$, $\Delta\mathbb{E}_y(t)$, $\Delta\mathbb{E}_z(t)$ of $\mathbb{E}[x_t]$, $\mathbb{E}[y_t]$, and $\mathbb{E}[z_t]$ associated with the dynamics (\ref{Lorenz})  to a spatially uniform, time-periodic perturbation (\ref{ex_tht3}) in (\ref{perLor}) with $\varepsilon = 0.1$; the remaining parameters are as in Figure~\ref{fdt_1}. See the main text in Example \ref{ex_L_fdt} for more details.}\label{fdt_3}
\end{figure}

 \begin{center} \textsc{Acknowledgements}\end{center}

{\small The research of M.B.~was supported by the Office of Naval Research grant ONR N00014-15-1-2351 and ONRG N62909-20-1-2037. K.U.~was supported by the first grant as a postdoctoral research fellow. }

\appendix
 \addtocontents{toc}{\protect\setcounter{tocdepth}{0}}
\section{Growth conditions and existence of absolute moments of solutions of non-autonomous SDE's}\label{app_hasminski}

Here, we provide explicit examples of two classes of time-periodic coefficients in the SDE~(\ref{NSDE}) which satisfy Assumption \ref{A2.1}.

\begin{lem}\label{lem_app1}
Let $\{\phi(t,s,\ccdot,\ccdot): \;t\geqslant s\}$ be a stochastic flow generated by the SDE (\ref{NSDE}), and let $V\in \mathcal{C}^{1,2}(\R\times\Rd;\Rp)$ be a Lyapunov function satisfying the first part of  (\ref{Lyap_f}) with some $1<p<\infty$. Assume that the following `dissipative' growth conditions hold on the coefficients of~(\ref{NSDE})
\begin{equation}\label{mbcnd} 
\langle b(t,x),x\rangle\leqslant L_{b_1}(t)-L_{b_2}(t)|x|^2, \quad \|\sigma(t,x)\|^2_{\textsc{hs}}\leqslant L_\sigma(t)(1+|x|^2).
\end{equation}
Suppose further that there exist bounded  functions $L_{b_1}(\cdot), L_{b_2}(\cdot), L_\sigma(\cdot)\in \mathcal{C}_\infty(\R;\Rp)$  such that 
\begin{align*}
\inf_{t\in \R} \Big(L_{b_2}(t)-2^{\frac{p}{2}-1}L_{b_1}-\frot (2^{\frac{p}{2}-1}+1)L_\sigma(t)(p-1)\Big)>0.
\end{align*}
Then, there exists a stochastic flow $\{\phi(t,s,\ccdot,\ccdot):\;s\leqslant t\}$ on $\Rd$ induced by the solutions of (\ref{NSDE}), which has a finite $p$-th absolute moment for all time. Furthermore, the following holds 
\begin{align*}
\limsup_{(t-s)\rightarrow \infty}\,\E\big[ V\big(t,\phi(t,s,x) -x\big)\big]<\infty
\end{align*}
in Assumption \ref{A2.1}(iii).
\end{lem}
\noindent {\it Proof.} \rm It can be established (in a similar way to that in \cite[Theorem 3.4.6]{Kunita}) that the growth conditions (\ref{mbcnd}) lead to existence and uniqueness of global solutions of 
\begin{align}\label{SDE_app}
dX^{s,x}_{t} = b(t,X^{s,x}_{t})dt+\sigma(t,X^{s,x}_{t})dW_{t-s},\quad X^{s,x}_{s} =x,
\end{align}
such that  $X^{s,x}_{t}(\om) = \phi(t,s,\om,x)$ $\p$-a.s., $t\geqslant s$.
In order to prove the main part of the Lemma, consider first $g(x) = |x|^p$, $p\geqslant 2$; then 
\begin{align}\label{Lmb}
 \mathcal{L}_t|x|^p &= p|x|^{p-2}\sum_{i=1}^db_i(t,x)x_i+ \frot p|x|^{p-4}\sum_{i,j=1}^d\Big\{|x|^2\delta_{ij}+(p-2)x_ix_j\Big\}(\sigma\sigma^T)_{ij}(t,x,x)\notag\\
&\leqslant p |x|^{p-2}(L_{b_1}(t)-L_{b_2}(t)|x|^2)+\frot L_\sigma(t)p(p-1)(1+|x|^2)|x|^{p-2}\notag\\[.2cm]
& = p \Big(L_{b_1}(t)+\frot L_\sigma(t)(p-1)\Big)|x|^{p-2}-p\Big(L_{b_2}(t)-\frot L_\sigma(t)(p-1)\Big)|x|^p.
\end{align}
Next, since $|x|^{p-2}\leqslant (1+|x|^2)^{\frac{p}{2}-1}\leqslant (1+|x|^2)^{\frac{p}{2}}\leqslant 2^{\frac{p}{2}-1}(1+|x|^p)$, we get 
\begin{align*}
\mathcal{L}_t|x|^p &\leqslant  2^{\frac{p}{2}-1}p\Big(L_{b_1}(t)+\frot L_\sigma(t)(p-1)\Big)\\&\hspace{4cm}-p\Big(L_{b_2}(t)-2^{\frac{p}{2}-1}L_{b_1}-\frot (2^{\frac{p}{2}-1}+1)L_\sigma(t)(p-1)\Big)|x|^p,
\end{align*}
which can be written as  
\begin{align}\label{mbbnd}
\mathcal{L}_t|x|^p &\leqslant \mathfrak{a}_p- \mathfrak{b}_p|x|^p,
\end{align}
with coefficients 
\begin{align}
\mathfrak{a}_p &= p\,2^{\frac{p}{2}-1}\sup_{t\in \R} \Big(L_{b_1}(t)+\frot L_\sigma(t)(p-1)\Big),\label{aa}\\[.3cm]
\mathfrak{b}_p &= p \inf_{t\in \R} \Big(L_{b_2}(t)-2^{\frac{p}{2}-1}L_{b_1}-\frot (2^{\frac{p}{2}-1}+1)L_\sigma(t)(p-1)\Big).\label{bb}
\end{align}
 
\smallskip
\noindent It turns out that sharper bounds can be obtained for $p=2,3$; these are derived in Proposition~\ref{shrpbnd}.

\medskip
Next, consider  $X^{s,x}_{t}(\om) = \phi(t,s,\om,x)$ solving (\ref{SDE_app}) and $g(x) = |x|^p$, $p\geqslant 2$. Then, 
It\^o's Lemma and the bound~(\ref{mbbnd}) lead to  
\begin{align}\label{mb22}
d\,\E\big[|X^{s,x}_{t}|^p\big] &=\E\big[\mathcal{L}_t |X^{s,x}_{t}|^p\big]dt \leqslant \big(\mathfrak{a}_p-\mathfrak{b}_p \,\E\big[|X^{s,x}_{t}|^p\big]\big)dt.
\end{align}
 
\smallskip
\noindent Therefore, based on the differential form of Gronwall's inequality,   (\ref{mb22}) yields 
\begin{align}\label{moms_app}
 \E |X^{s,x}_{t}|^p\leqslant e^{-\mathfrak{b}_p(t-s)}\E|x|^p+\frac{\mathfrak{a}_p}{\mathfrak{b}_p}\Big(1-e^{-\mathfrak{b}_p(t-s)}\Big).
\end{align}
Consequently, for $p\geqslant 2$, and $L_{b_1}, L_{b_2}, L_\sigma$ such that $\mathfrak{b}_p>0$ in (\ref{bb}), we have   
\begin{align*}
0\leqslant \lim_{(t-s)\rightarrow \infty} \E|\phi(t,s,x)|^p\leqslant \frac{\mathfrak{a}_p}{\mathfrak{b}_p}<\infty.
\end{align*}

\smallskip
\noindent For  $1< p<2$ we use H\"older's inequality and obtain 
\begin{align*}
\E\left[ |X^{s,x}_{t}|^p\right]\leqslant \Big(\E\big[ |X^{s,x}_{t}|^{2p}\big]\Big)^{\frac{1}{2}}.
\end{align*}
Thus, analogous derivations to those in (\ref{Lmb}) onwards can be carried out for $p' = 2p> 2$, with $1< p<2$.

\medskip
Finally, note that for  $Y^{s,x}_{t}(\om) = \phi(t,s,\om,x)-x$ we obtain 
\begin{align*}
dY^{s,x}_{t} = b(t,Y^{s,x}_{t}+x)dt+\sigma(t,Y^{s,x}_{t}+x)dW_t,\quad Y^{s,x}_{s} =0,
\end{align*}
so that analogous calculations in conjunction with Assumption (\ref{A2.1}) lead to 
\begin{align*}
\E\big[V(t,\phi(t,s,x)-x)\big]^p &\leqslant \mathfrak{C}\hspace{.05cm}\E |\phi(t,s,x)-x|^{p} \leqslant \mathfrak{C}\,\frac{\mathfrak{a}_{p}}{\mathfrak{b}_{p}}\Big(1-e^{\mathfrak{b}_{p}(t-s)}\Big), 
\end{align*}
Consequently, for $p> 1$ and $0<L_{b_1}, L_{b_2}, L_\sigma<\infty$ such that $\mathfrak{b}_{p}>0$

\begin{align*}
0\leqslant \lim_{s\rightarrow -\infty} \sup_{s\leqslant t}\E \big[V\big(t,\phi(t,s,x)-x\big)\big]<\infty,
\end{align*}
\begin{align*}
\hspace{3.9cm} \quad 0\leqslant \lim_{t\rightarrow \infty} \sup_{s\leqslant t}\E\big[V\big(\phi(t,s,x)-x\big)\big]<\infty.  \hspace{4.2cm}\qed
 \end{align*}

\begin{prop}[Sharper conditions for existence  of absolute moments $p=2,3$ for SDE's with dissipative growth conditions]\label{shrpbnd} 
\mbox{}

Consider the same setup as in Lemma~\ref{lem_app1}. Suppose further that there exist bounded  functions $L_{b_1}(\cdot), L_{b_2}(\cdot), L_\sigma(\cdot)\in \mathcal{C}_\infty(\R;\Rp)$  such that 
\begin{align*}
\tilde{\mathfrak{b}}_2 &= 2 \inf_{t\in \R} \Big(L_{b_2}(t)-\frot L_\sigma(t)\Big)>0.
\end{align*}
Then, there exists a stochastic flow $\{\phi(t,s,\ccdot,\ccdot):\;s\leqslant t\}$ on $\Rd$ induced by the global solutions of (\ref{NSDE}), which has a finite second  absolute moment for all time. Moreover, if for some $\varkappa>0$
\begin{equation}
\tilde{\mathfrak{b}}_3= 3 \inf_{t\in \R} \Big(L_{b_2}(t)-L_\sigma(t)-\frac{4}{27\varkappa^2}(L_{b_1}(t)+ L_\sigma(t))\Big)>0,
\end{equation}
then the stochastic flow has a finite third moment for all time. 
\end{prop}

\noindent {\it Proof:} For $p=2$, we proceed in a  way similar to (\ref{Lmb}), and we have 
\begin{align}\label{Lmb_2}
 \mathcal{L}_t|x|^2 &\leqslant 2(L_{b_1}(t)-L_{b_2}(t)|x|^2)+L_\sigma(t)(1+|x|^2)\notag\\[.2cm]
 & = 2 \Big(L_{b_1}(t)+\frot L_\sigma(t)\Big)-2\Big(L_{b_2}(t)-\frot L_\sigma(t)\Big)|x|^2\leqslant  \tilde{\mathfrak{a}}_2-\tilde{\mathfrak{b}}_2|x|^3,
\end{align}
where
\vspace{-.2cm}
\begin{align}
\tilde{\mathfrak{a}}_2 &= 2\sup_{t\in \R} \Big(L_{b_1}(t)+\frot L_\sigma(t)\Big),\label{aa_2}\\
\tilde{\mathfrak{b}}_2 &= 2 \inf_{t\in \R} \Big(L_{b_2}(t)-\frot L_\sigma(t)\Big).\label{bb_2}
\end{align}
Thus, based on (\ref{moms_app}) the second absolute moment exists for all time if $\tilde{\mathfrak{b}}_2>0$.

For $p=3$ we have
\begin{align}\label{Lmb_3}
\mathcal{L}_t|x|^3 &\leqslant 3|x|(L_{b_1}(t)-L_{b_2}(t)|x|^2)+3|x| L_\sigma(t)(1+|x|^2)\notag\\[.2cm]
& = 3 \Big(L_{b_1}(t)+ L_\sigma(t)\Big)|x|-3\Big(L_{b_2}(t)-L_\sigma(t)\Big)|x|^3\notag\\
& \leqslant 3\varkappa \Big(L_{b_1}(t)+ L_\sigma(t)\Big)-3\Big(L_{b_2}(t)-L_\sigma(t)-\frac{4}{27\varkappa^2}(L_{b_1}(t)+ L_\sigma(t))\Big)|x|^3\notag\\
& \leqslant  \tilde{\mathfrak{a}}_3-\tilde{\mathfrak{b}}_3|x|^3,
\end{align}
with 
\vspace{-.2cm}
\begin{align}
\tilde{\mathfrak{a}}_3 &= 3\varkappa\sup_{t\in \R} \Big(L_{b_1}(t)+ L_\sigma(t)\Big),\label{aa_3}\\
\tilde{\mathfrak{b}}_3&= 3 \inf_{t\in \R} \Big(L_{b_2}(t)-L_\sigma(t)-\frac{4}{27\varkappa^2}(L_{b_1}(t)+ L_\sigma(t))\Big),\label{bb_3}
\end{align}
where we used the fact that $ |x|\leqslant \varkappa +\frac{4}{27\varkappa^2}|x|^3$, $\varkappa>0$. Thus, based on (\ref{moms_app}) the third  absolute moment exists for all time if $\tilde{\mathfrak{b}}_3>0$. \qed

\begin{rem}\rm
It is worth noting that, for $L_{b_1}, L_{b_2},L_\sigma$ constant, and such that $\tilde{\mathfrak{b}}_3> 0$, the upper bound on the asymptotic moment $\E|\phi(t,s,x)|^p$ for $p=3$ is optimised for 
$\varkappa^2 = 12\left(\frac{L_{b_1}+L_\sigma}{L_{b_2}-L_\sigma}\right)$ 
so that 
\begin{equation}\label{exx_1}
\min_{\varkappa>0}\,\frac{\tilde{\mathfrak{a}}_3}{\tilde{\mathfrak{b}}_3} = \left(27\right)^{1/2}\left(\frac{L_{b_1}+L_\sigma}{L_{b_2}-L_\sigma}\right)^{3/2}. 
\end{equation}
Moreover, 
\begin{equation}\label{exx_2}
\frac{\tilde{\mathfrak{a}}_2}{\tilde{\mathfrak{b}}_2}\leqslant  \left(\frac{L_{b_1}+L_\sigma}{L_{b_2}-L_\sigma}\right) \leqslant \left(\min_{\varkappa>0}\,\frac{\tilde{\mathfrak{a}}_3}{\tilde{\mathfrak{b}}_3}\right)^{2/3} = \left(27\right)^{1/3}\left(\frac{L_{b_1}+L_\sigma}{L_{b_2}-L_\sigma}\right); 
\end{equation}
this fact merely reflects the Jensen's inequality for the second and third absolute moments, i.e., $\E|X_t^{s,x}|^2\leqslant (\E|X_{t}^{s,x}|^3)^{2/3}$, but it is useful in Example \ref{Lor_ex}. 
\end{rem}

\begin{lem}\label{weird_cond}
Let $\{\phi(t,s,\ccdot,\ccdot): t\geqslant s\}$ be a stochastic flow generated by the SDE (\ref{NSDE}) and let $V\in \mathcal{C}^{1,2}(\R\times\Rd;\Rp)$ be a Lyapunov function satisfying the first part of condition (\ref{Lyap_f}). Suppose further that there exist $L_b(\cdot), L_\sigma(\cdot), C(\cdot)\in \mathcal{C}_\infty(\R;\Rp)$ such that 
\begin{align}\label{Lgrowth}
\langle b(t,x), x\rangle\leqslant L_{b}(t)(1+\vert x\vert^2), \quad \Vert \sigma(t,x)\Vert^2_{\textsc{hs}}\leqslant L_\sigma(t)\big( 1+\vert x\vert^2\big),
\end{align}
and 
\begin{align}\label{constr1}
 0<\limsup_{(t-s)\rightarrow \infty}\exp\left(\int_s^tC(u,p)du\right)<\infty,
\end{align}
where $C(t,p) = L_b(t)+\frac{1}{2}(p-1)L_\sigma(t),$ for some $1<p<\infty$. 

Then, for $x\in \Rd,$ there exist global solutions of (\ref{NSDE}) such that the stochastic flow $\phi$ induced by the solutions of (\ref{NSDE}) has a finite $p$-th absolute moment for all time. Furthermore, the following holds 
\begin{align*}
\limsup_{(t-s)\rightarrow \infty}\,\E\big[ V(t,\phi(t,s,x) -x)\big]<\infty
\end{align*}
in Assumption \ref{A2.1}(iii).
\end{lem}

\noindent{\it Proof}.  \rm 
 First, we suppose that $p\geqslant 2$ and set $g(x) = 1+|x|^2$ and $\varphi(x) = g(x)^{\frac{p}{2}};$ then 
\begin{align*}
\mathcal{L}_t\varphi(x) = pg(x)^{\frac{p}{2}-1}\sum_{i=1}^db_i(t,x)x_i+\frac{1}{2}pg(x)^{\frac{p}{2}-2}\sum_{i,j=1}^d\Big\{ g(x)\delta_{ij}+(p-2)x_ix_j\Big\}(\sigma\sigma^T)_{ij}(t,x,x).
\end{align*}
By the Growth conditions (\ref{Lgrowth}) on the coefficients $ b, \sigma$, we obtain
\begin{align*}
\mathcal{L}_t\varphi(x)\leqslant p\,C(t,p)\varphi(x),
\end{align*}
where $C(t,p) = L_b(t)+ c(t)+\frac{1}{2}(p-1)L_\sigma(t).$
Next, let $Y_{t,s}^x(\om) = \phi(t,s,\om,x)-x,$ it follows that $Y_{t,s}^x(x)$ solves the following SDE
\begin{align*}
dY_{t,s}^x = b(t,Y_{t,s}^x+x)dt+\sigma(t,Y_{t,s}^x+x)dW_t,\quad Y_{s,s}^x =0.
\end{align*}

\vspace{.2cm}
\noindent By It\^o's formula, we have 
\begin{align*}
\E\left[\varphi(Y^x_{t,s})\right] &=\varphi(0)+\E\left[\int_{s}^t\mathcal{L}_t\varphi(Y_{u,s}^x)du\right]\\[.2cm]
&\hspace{0cm} \leqslant \varphi(0)+p\int_{s}^t C(u,p)\E\left[\varphi(Y_{u,s}^x)\right]du.
\end{align*}
By Gronwall's inequality, we have 
\begin{align*}
\E\left[\varphi(Y_{t,s}^x)\right] \leqslant \varphi(0)\exp\left(p\int_{s}^tC(u,p)du\right).
\end{align*}

\vspace{.2cm}
\noindent But $Y_{t,s}^x(\om) = \phi(t,s,\om,x)-x$ \,$\p$-a.s., $\varphi(x) = g(x)^{\frac{p}{2}} = \left( 1+|x|^2\right)^{\frac{p}{2}}$ and  $\varphi(0) = 1,$ thus, 
\begin{align*}
\E\Big[\left(1+|\phi(t,s,\om,x)-x|^2\right)^{\frac{p}{2}}\Big]\leqslant \exp\left(p\int_{s}^tC(u,p)du\right).
\end{align*}

\vspace{.3cm}
\noindent Next, note for $p\geqslant 2,$ $|x|^p\leqslant (1+|x|^2)^{\frac{p}{2}}$  and by the assumption that $V(t,x)\leqslant \mathfrak{C}|x|^p,$ we obtain
\begin{align*}
\E\Big[ V(t,\phi(t,s,x)-x)\Big]\leqslant \mathfrak{C} \E\Big[\left(1+|\phi(t,s,x)-x|^2\right)^{\frac{p}{2}}\Big] \\ \leqslant \mathfrak{C} \exp\left(p\int_{s}^tC(u,p)du\right).
\end{align*}
Since $C(\ccdot,p)\in \mathcal{C}(\R;\R)$ such that $\displaystyle \limsup_{s\rightarrow-\infty}\exp\left(\int_s^tC(u,p)du\right)<\infty,$ then for $p\geqslant 2,$ we have
\begin{align}\label{Cas1}
\E\Big[ V(t,\phi(t,s,x)-x)\Big]\leqslant \mathfrak{C}\limsup_{s\rightarrow -\infty}\exp\left(p\int_{s}^tC(u,p)du\right)<\infty.
\end{align}
The case where $1\leqslant p<2,$ we use H\"older's inequality, namely
\begin{align*}
\E\Big[ V(t,\phi(t,s,x)-x)\Big]\leqslant \mathfrak{C}\E\left[ |\phi(t,s,x)-x|^p\right]\leqslant \mathfrak{C}\Big(\E\left[ |\phi(t,s,x)-x|^{2p}\right]\Big)^{\frac{1}{2}},
\end{align*}
and the rest follows, since in this case $2p\geqslant 2$.
\qed

\section{Strong Feller property for flows induced by non-autonomous SDE's}
In order to prove Theorem~\ref{StrongF}, which is a generalisation of standard results to non-autonomous SDE's,  we first outline some basic notions from Malliavin calculus; the actual proof is given in and discussed in Appendix~\ref{SfA}. 

\addtocontents{toc}{\protect\setcounter{tocdepth}{1}}
\subsection{Malliavin calculus estimates}\label{MallSec}
Establishing the strong Feller property of Markov evolutions $(\mathcal{P}_{s,t})_{t\geqslant s}$ in our setting requires some estimates rooted in Malliavin calculus. We recall the main concepts and results on the Wiener space $(\Om, \mathcal{F},\p)$; see  (e.g.,~\cite{Hairer10, Hairer11, Malliavin, Millet, Nulart, Watanabe}) for a comprehensive treatment. 
To this end, consider the Hilbert space $\mathcal{H} = L^2([s,\infty);\R^m)$ equiped with the inner product
 \begin{align*}
 \langle \eta_1, \eta_2\rangle_{\mathcal{H}} = \int_s^{\infty} \eta_1(t)\cdot \eta_2(t)dt.
 \end{align*}
 
 \smallskip
 \noindent For a Hilbert space $E$ and a real number $p\geqslant 1$,  $L^{p}(\Om; E)$ the space of $E$-valued random variable $\xi$ such that $\displaystyle \E(\|\xi\|^p_E):=\int_{\Om}\|\xi\|^p_{E}\,d\hspace{.03cm}\p<\infty.$ Also, we set $
\displaystyle L^{\infty-}(\Om; E) := \bigcap_{1\leqslant p<\infty}L^{p}(\Om; E)$. 

Following the approach due to Malliavin (e.g.,~\cite{Malliavin, Nulart}), we introduce a derivative operator $\D$ for a random variable $G$ on the space $L^{\infty-}(\Om;E).$ We say that $G\in \mathbb{D}^{1,\infty}(E)$ if there exists $\mathcal{D}G\in L^{\infty-}(\Om;\mathcal{H}\times E)$ such that for any $\eta\in \mathcal{H},$
\begin{align*}
\lim_{\varepsilon\rightarrow 0}\E\bigg\| \frac{G\left(\om+\varepsilon\int_s^{.}\eta(\ell)d\ell\right)-G(\om)}{\varepsilon} -\langle \mathcal{D}G, \eta\rangle_{\mathcal{H}}\bigg\Vert^{p}_{E} =0,
\end{align*}

\smallskip
\noindent holds for every $p\geqslant 1.$ In this case, one defines the Malliavin derivative of $G$ in the direction of $\eta\in \mathcal{H}$ by $\mathcal{D}^\eta G:=\langle \mathcal{D}G, \eta\rangle_{\mathcal{H}}.$ For any $p\geqslant 1,$ we define the Sobolev space $\mathbb{D}^{1,p}(E)$ as the completion of $\mathbb{D}^{1,\infty}(E)$ under the norm
\begin{align*}
\| G\|_{1,p,E} = \big(\E\| G\|^p_{\mathcal{H}}\big)^{1/p}+ \big(\E\| \D G \|^p_{\mathcal{H}\times E}\big)^{1/p}.
\end{align*}

\medskip
\noindent We define the $k$-th Malliavin derivative by $\D^kG = \D(\D^{k-1}G),$ which is a random variable with values in $\mathcal{H}^{\otimes k}\times E.$ For any integer $k\geqslant 1,$ the Sobolev space $\mathbb{D}^{k,p}(E)$ is the completion of $\mathbb{D}^{k,\infty}(E)$ under the norm
\begin{align*}
\|G\|_{k,p,E} =\| G\|_{k-1,p,E}+\|\D^{k}G\|_{1,p,\mathcal{H}^{\times k}\times E}.
\end{align*}

\medskip
\noindent It turns out that $\D$ is a closed operator from $L^{p}(\Om; E)$ to $L^{p}(\Om; \mathcal{H}\times E)$. The ajoint $\delta$ of the operator $\D$ called the divergence operator is continuous from $\mathbb{D}^{1,p}(\mathcal{H}\times E)$ to $L^{p}(\Om; E)$ for any $p>1$, with the duality relationship given as 
\begin{align}\label{Div}
\E\big[\langle \D G, u\rangle_{\mathcal{H}\times E}\big] = \E\big[ \langle G, \delta(u)\rangle_{E}\big],
\end{align}

\medskip
for any $G\in \mathbb{D}^{1,p}(\mathcal{H}\times E)$ and $u\in \mathbb{D}^{1,q}(\mathcal{H}\times E),$ with $\frac{1}{p}+\frac{1}{q} =1.$  

\medskip
\noindent Throughout the remaining part of this section we assume the following notation: 
\begin{itemize}[leftmargin = 0.4cm]
\item[-] $C$ is a generic constant which may depend on $T,$ the exponent $p>1,$ the initial point $x$ and fixed element $\eta$ of the Hilbert space $\mathcal{H}= L^2\big([s, \infty);\R^m\big).$
\item[-] ({\bf $H_n$}) denotes a class of coefficients $b, \sigma_k$, $1\leqslant k\leqslant m$, where $\sigma_k$ are columns of $\sigma$, such that $b(t,\ccdot)\in \tilde{\mathcal{C}}^n, \sigma_k(t,\ccdot) \in \bar{\bar{\mathcal{C}}}^{n}\big(\Rd\big)$.
\end{itemize}

\begin{prop}
Suppose the coefficients $b,\sigma$ of the SDE (\ref{NSDE}) are in the class ($H_2$). Then, for any $t\geqslant s,$ we have $\phi(t,s,\ccdot,\ccdot)\in \mathbb{D}^{1,\infty}(\Rd)$ and the Malliavin derivative $\D^\eta \phi(t,s)$ of $\phi(t,s)$ in the direction of $\eta = (\eta^1, \eta^2,\cdots, \eta^m)\in \mathcal{H}$ is the unique solution of the following affine SDE
\begin{align*}
\begin{cases}
\displaystyle d\D^\eta \phi(t,s)= D_xb\big(t,\phi(t,s)\big)\D^\eta \phi(t,s)dt + \displaystyle\sum_{k=1}^m D_x\sigma_k\big(t,\phi(t,s)\big)\D^\eta \phi(t,s) dW_t^k \\[.4cm] \hspace{6.5cm} \displaystyle +\sum_{k=1}^m\sigma_k\big(t,\phi(t,s)\big)\eta^k(t)dt, \quad t>s,\\
\D^\eta \phi(s,s) =0. &
\end{cases}
\end{align*}
\end{prop}
\begin{cor}[Chain rule, cf. \cite{Nulart}]\label{ChainR}\rm
Suppose that condition ({\bf $H_2$}) holds true. Then, for any $\eta\in \mathcal{H},$ $p\geqslant 2$ and for any $\varphi\in \mathcal{C}_\infty^2(\Rd)$, we have 
\begin{align*}
\lim_{\varepsilon\rightarrow 0}\E\left\vert \frac{\varphi\big(\phi^{\varepsilon\eta}(t,s)\big) - \varphi\big(\phi(t,s)\big)}{\varepsilon} - D_x\varphi\big(\phi(t,s)\big)\D^\eta \phi(t,s)\right\vert^p=0,
\end{align*} 
where $\phi^{\varepsilon \eta}(t,s), \; t\geqslant s,\; \varepsilon\in (0, 1)$ is the solution of the following perturbed SDE
\begin{align*}
\begin{cases} \displaystyle d\phi^{\varepsilon \eta}(t,s) = b\big(t,\phi^{\varepsilon\eta}(t,s)\big)dt +\sum_{k=1}^m\sigma_k\big(t, \phi^{\varepsilon \eta}(t,s)\big)dW_t^k+\varepsilon\sum_{k=1}^m\sigma_k\big(t, \phi^{\varepsilon\eta}(t,s)\big)\eta^k(t)dt,\\
\phi^{\varepsilon\eta}(s,s)= x\in \Rd.
\end{cases}
\end{align*}
Moreover, $\varphi(\phi(t,s))\in \mathbb{D}^{1,\infty}(\R)$ and $\D \varphi(\phi(t,s)) = D_x\varphi(\phi(t,s))\D \phi(t,s).$ 
\end{cor}
\begin{definition}[Mean square gradient]\rm Let $G(x):\Om\rightarrow\Rd$ be a measurable function for all $x\in \Rd$ and $i\in F$. We say that the mean square gradient of $G(x)$ with respect to $x$ exists if there is a linear map $A(x): \Om\rightarrow \R^{d\times d}$ such that for any $v\in \Rd$, 
\begin{align*}
\lim_{\varepsilon\rightarrow 0}\E\left\vert \frac{G(x+\varepsilon v) - G(x)}{\varepsilon} -A(x)v\right\vert^2 =0.
\end{align*}
We denote the mean square gradient matrix $A(x)$ by $D_xG(x).$
\end{definition}
\begin{theorem}[e.g., \cite{Kunita, Malliavin, Nulart}]
Assume the condition ({\bf $H_2$}) holds true. Let $\phi(t,s,\om,x),\; t\geqslant s $, be the solution of the SDE (\ref{NSDE}). Then, the mean square gradient of $\phi(t,s,\ccdot,x)$ with respect to $x$ exists. If we define $J_{t,s} = D_x\phi(t,s,\ccdot,x),$ then 
\begin{align}\label{DfLow}
\begin{cases}\displaystyle dJ_{t,s} = D_xb(t, \phi(t,s,\ccdot,x))J_{t,s}+\sum_{k=1}^m D_x\sigma_k(t, \phi(t,s,\ccdot,x))J_{t,s}dW_t^k, \quad t\geqslant s,\\
J_{s,s} = I,
\end{cases}
\end{align}
where $I$ is a $d\times d$ identity matrix. Moreover, the inverse $J_{t,s}^{-1}$ of $J_{t,s}$ exists and satisfy 
\begin{align}
\begin{cases} \displaystyle dJ_{t,s}^{-1} = -J_{t,s}^{-1}\Big( D_xb\big(t,\phi(t,s,\ccdot,x)\big)-\sum_{k=1}^m D_x\sigma_k\big(t, \phi(t,s,\ccdot,x)\big)D_x\sigma_k\big(t, \phi(t,s,\ccdot,x)\big)\Big) dt\\
 \hspace{8.8cm} \displaystyle -\sum_{k=1}^m J_{t,s}^{-1}D_x\sigma_k\big(t, \phi(t,s,\ccdot,x)\big)dW_t^k,\\
J_{s,s}^{-1} = I.
\end{cases}
\end{align}
\end{theorem}

We shall refer to the mean square gradient $\{J_{t,s} = D_x\phi(t,s,\ccdot,\ccdot):\;s\leqslant t\}$ as derivative flow of $\{\phi(t,s,\ccdot, \ccdot):\;s\leqslant t\}$. Next, we provide a crucial $L^{p}$ bound for the derivative flow flow $J_{t,s}$ and that of its inverse $J_{t,s}^{-1}.$  

\begin{lem}\label{3,6}
Suppose the condition ({\bf $H_2$}) holds. Then for any $p\geqslant 2,$ there exists a positive constant $C = C(T, p)$ such that 
\begin{align}
\E\left(\sup_{s\leqslant t,u\leqslant s+T}| J_{t,u}|^p\right)\leqslant C \quad  \text{and}\quad \E\left(\sup_{s\leqslant t, u\leqslant s+T}| J^{-1}_{t,u}|^p\right)\leqslant C.
\end{align}
\end{lem}

\bigskip
Now, let $\D_r\phi(t,s)$ be the solution of the following SDE:
\begin{align}\label{MalliD}
\begin{cases} & \D_u\phi(t,s) \displaystyle = \sigma(u, \phi(u,s))+\int_u^t\D_xb(\ell, \phi(\ell,u))\D_u\phi(\ell,u)d\ell \\  & \hspace{4.5cm}+ \displaystyle \sum_{k=1}^m\int_u^tD_x\sigma_k(\ell, \phi(\ell,u))\D_u\phi(\ell,u)dW_{\ell}^k,\quad \text{for}\; t\geqslant u,\\
&\D_u \phi_{t,s} =0,\quad \text{for}\; s\leqslant t<u.
\end{cases}
\end{align}
Comparing the  SDE's (\ref{DfLow}) and (\ref{MalliD}), we obtain the following by the variation of parameters formula 
\begin{align*}
\begin{cases} \D_u\phi(t,s) = J_{t,u}\sigma(\phi(t,s)), \quad s\leqslant u\leqslant t\leqslant s+T,\\ \D_u \phi(t,s) =0, \quad u>t.
\end{cases}
\end{align*}

Next, we recall a result on the Malliavin differentiability of the derivative flow $J_{t,s}, \; t\geqslant s.$ To this end, let's denote by $\D_u^\ell$ the Malliavin derivative with respect to the $\ell$-th component of the Brownian motion $W$ at time $u.$

\begin{lem}\label{MalliDflow}
Suppose that the condition ({\bf $H_3$}) holds. Then for all $s\leqslant t\leqslant s+T,\; J_{t,s}\in \mathbb{D}^{1,\infty}(\Rd\times\Rd)$ and for any $p\geqslant 2,$ there exists a positive constant $C=C(T,p,x),$ such that for all $j=1,\cdots, m$ and $u\in [s, s+T]$,
\begin{align*}
\E\left[ \sup_{s\leqslant t\leqslant s+T}|\D^j_u J_{t,s}|^p\right]\leqslant C.
\end{align*}
Moreover, for any $t\leqslant s+T,\; X(t,s)\in \mathbb{D}^{2,\infty}(\Rd)$ and for any $p\geqslant 2,$ there exists a positive constant $C= C(T,p,x)$ such that for all $j, l=1,\cdots, m$ and $\varsigma, u\leqslant t$
\begin{align*}
\E\vert \D_u^{j}(\D_\varsigma^l \phi(t,s))|^p\leqslant C.
\end{align*}
\end{lem}

\begin{rem}\rm If ({\bf $H_{\infty}$}) holds true, then $\phi(t,s)\in \mathbb{D}^{\infty}(\Rd)$ and $J_{t,s}\in \mathbb{D}^{\infty}(\Rd\times\Rd).$
\end{rem}

Denote by $(\D \phi(t,s))^T$ the transpose of the Malliavin derivative $\D \phi(t,s).$ From the relationship  $\D_u \phi(t,s) = J_{t,u} \sigma(u, \phi(u,s)),$ we have $(\D_u \phi(t,s))^T = \sigma(u, \phi(u,s))^T J_{t,u}^T$.
\begin{definition}[Malliavin covariance; e.g.,~\cite{Malliavin, Nulart}]\rm
The Malliavin covariance $M_{t,s}$ of the random vector $\phi(t,s)$ is defined by 
\begin{align*}
M_{t,s} = \langle \D \phi(t,s), (\D \phi(t,s))^T\rangle_{\mathcal{H}} & = \int_s^tJ_{t,\ell}\sigma(\ell, \phi(\ell, s)) \sigma(\ell, \phi(\ell,s))^T J_{t, \ell}^Td\ell\\
&=J_{t,s}\left[\int_s^t J_{\ell,s}^{-1} \sigma(\ell, \phi(\ell, s)) \sigma(\ell, \phi(\ell,s))^T (J_{\ell, s}^{-1})^Td\ell\right] J_{t,s}^{T}\\
&= J_{t,s}C_{t,s} J_{t,s}^T,
\end{align*}
where $C_{t,s}$ is defined by 
\begin{align*}
C_{t,s} = \int_s^t J_{\ell,s}^{-1} \sigma(\ell, \phi(\ell, s)) \sigma(\ell, \phi(\ell,s))^T (J_{\ell, s}^{-1})^Td\ell
\end{align*}
is the so-called reduced Malliavin covariance of $\phi(t,s).$
\end{definition}

We conclude this section by elucidating the invertibility of the Malliavin covariance almost surely and its integrability of all negative orders. 
\begin{prop}[e.g.,~\cite{Hairer10,  Nulart, Watanabe}]\label{SMOOTH} Suppose Assumption \ref{Homa} holds. Then, for every $t \geqslant s,$ the Malliavin covariance matrix $M_{t,s}$ of the random vector $\phi(t,s)$ is invertible $\p$\,-\,a.s., and
$\E\left[ \det(M_{t,u}^{-p})\right]<\infty,$
 for every $t,u\in [s, s+T], \; T>0$, and $p>1.$
  Moreover, for any $x\in \Rd$, $s\leqslant t$, the law of $\phi(t,s,\ccdot,x)$ is absolutely continuous with respect to the Lebesgue measure on $\Rd$ and the probability density is smooth.
\end{prop}

\subsection{Strong Feller property for non-autonomous dynamics.}\label{SfA}
A transition evolution denoted by $(\Pt_{s,t})_{t\geqslant s}$ (\ref{calP}) and induced by a stochastic flow $\{\phi(t,s,\ccdot,\ccdot)\!:\; s\leqslant t\}$ has the strong Feller property (i.e., $\Pt_{s,t}\varphi \in \mathcal{C}_\infty(\Rd)$ for any $\varphi\in \mathbb{M}_\infty(\Rd)$) if and only if
\begin{itemize}[leftmargin =0.9cm]
\item[(a)] $(\Pt_{s,t})_{t\geqslant t}$ is a Feller semigroup; i.e., $\Pt_{s,t}: \mathcal{C}_\infty(\Rd)\rightarrow\mathcal{C}_\infty(\Rd)$.
\item[(b)] For any $\varphi\in \mathcal{C}_\infty(\Rd)$ the family $(\Pt_{s,t}\varphi)_{t\geqslant s}$ is equicontinuous.
\end{itemize}  
The first condition follows from the existence of the stochastic flow (see, e.g.,~\cite{Kunita, Has12}); here, we are concerned with flows associated with solutions of the non-autonomous SDE~(\ref{NSDE}).  Thus,  shall only derive the second item.

Intuitively, the strong Feller property states that for  sufficiently close  initial data $x,\,y$ and any realisation $\om$ of the past driving noise, one can construct a coupling between two solutions $\phi(t,s,\om,x)$ and $\phi(t,s,\om,y)$ such that with probability close to $1$ as $x\rightarrow y,$ one has $\phi(t,s,\om,x) = \phi(t,s,\om,y),$ for $t\geqslant s$ (e.g.,~\cite{Hairer11, Hairer11a}). One way of achieving such a coupling (e.g.,~\cite{Hairer11, HairerD}) is via a change of measure on the driving process for one of the two solutions such that the noises $W_t^x$ and $W_t^y$ driving the solutions $\phi(t,s,\om,x)$ and $\phi(t,s,\om,y)$,  are related by 
 \begin{align*}
  dW_t^x = dW_t^{y}+\eta_t^{x,y}dt,
 \end{align*}
 where $\eta_t^{x,y}$ is a control process that steers the solution $\phi(t,s,\om,x)$ towards the solution $\phi(t,s,\om,y)$. If one sets $y = x+\varepsilon\eta$ and looks for a control of the form $\eta^{x,y}_t = \varepsilon \eta,$ then in the limit as $\varepsilon\rightarrow 0,$ the scheme will induce a deformation onto the solution $\phi(t,s,\om,x)$ after time $t$ in the form of {\it Malliavin derivative} of $\phi(t,s,\om,x)$ at $\om\in \Om$ in the direction of $\eta\in \mathcal{H}$, $\mathcal{H} = L^2([s, \infty);\Rd)$;  i.e., 
  \begin{align*}
 \langle \mathcal{D}\phi(t,s,\om,x),\eta\rangle_{\mathcal{H}} = \mathcal{D}^\eta \phi(t,s,\om,x).
  \end{align*}
 On the other hand, the effect of the perturbation of initial condition by $v$ is given by the directional derivative of the solution $\phi(t,s,\om,x)$ at $x$ along $v$; i.e.,
  \begin{align*}
  D_x\phi(t,s,\om,x)v = J_{s,t}(\om,x)v.
  \end{align*}
  In order to assert the strong Feller property, one has to find a control $\eta^v$ (e.g.,~\cite{Hairer11, Hairer11a}) such~that 
  \begin{align}\label{Task}
 \langle\D \phi(t,s,\om,x), \eta^v\rangle_{\mathcal{H}} = J_{t,s}(\om,x)v, 
  \end{align}
 where for brevity of notation we skip the explicit dependence on $\om$ and $x$.
\begin{theorem}\label{AStrongF} Suppose that Assumption \ref{Homa} hold true. Then for any $t\in [s, s+T],$ there exist $C_T>0$ such that, for any $x,y\in \Rd$ and any $\varphi\in \mathcal{C}_\infty(\Rd),$ we have
  \begin{align*}
  \vert \mathcal{P}_{s,t}\varphi(x)-\mathcal{P}_{s,t}\varphi(y)\vert \leqslant C_T\Vert \varphi\Vert_{\infty}\vert x-y\vert.
\end{align*}   
  \end{theorem}
  \noindent {\it Proof.} 
First, we find a control satisfying (\ref{Task}). To this end, for any $v\in \Rd$ with $\vert v\vert =1,$ let $\eta^v = (\D \phi(t,s))^{T}M^{-1}_{t,s}J_{t,s}v$.  Then,
\begin{align*}
 \langle \D \phi(t,s), \eta^v\rangle &= \langle \D \phi(t,s), (\D \phi(t,s))^TM_{t,s}^{-1}J_{t,s}v\rangle_{\mathcal{H}}= \langle \D \phi(t,s), (\D \phi(t,s))^T\rangle_{\mathcal{H}}M_{t,s}^{-1}J_{t,s}v= J_{t,s}v.
 \end{align*}

\smallskip
\noindent Next, we show that $\eta^v \in \mathbb{D}^{1,p}(\mathcal{H})$ for any $p\geqslant 2$. In fact, by chain rule of differentiation,
 \begin{align*}
 \D_\varsigma^k\eta^v &= (\D^k_\varsigma(\D \phi(t,s))^T)M_{t,s}^{-1}J_{t,s}v+(\D \phi(t,s))^TM_{t,s}^{-1}(\D^k_\varsigma J_{t,s})v+(\D \phi(t,s))^T(\D^k_\varsigma M_{t,s}^{-1})J_{t,s}v\\[.2cm]
  &= (\D^k_\varsigma(\D \phi(t,s))^T)M_{t,s}^{-1}J_{t,s}v+(\D \phi(t,s))^TM_{t,s}^{-1}(\D^\ell_\varsigma J_{t,s})v\\[.2cm]
  & \hspace{4.cm} -(\D \phi(t,s))^TM_{t,s}^{-1}\big[ \langle \D_\varsigma^k(\D \phi(t,s)), (\D \phi(t,s))^T\rangle_{\mathcal{H}}\\[.2cm]
  & \hspace{7.5cm}+\langle \D \phi(t,s), \D_\varsigma^k(\D \phi(t,s))^T\rangle_{\mathcal{H}}\big] M_{t,s}^{-1} J_{t,s} v.
  \end{align*}
  By Lemmas \ref{3,6}, \ref{MalliDflow}, and Proposition \ref{SMOOTH} , we arrive at
  \begin{align}\label{SKEst}
  \E\|\eta^v\|^p_{\mathcal{H}}+  \E\|\D\eta^v\|^p_{\mathcal{H}\times\mathcal{H}}\leqslant \E\|\eta^v\|_{\mathcal{H}}^p+\sum_{k=1}^m\E\left[ \int_s^t\| \D_\varsigma^k\eta^v\|^p_{\mathcal{H}}d\varsigma\right]<\infty.
  \end{align}
 Recalling that $\eta^v_{\varsigma,s} = \sigma(\varsigma, X(\varsigma,s))^TJ_{t,\varsigma}^TM_{t,s}^{-1}J_{t,s} v.$ Then, for $\varphi \in \mathcal{C}_\infty^1(\Rd),$ we have
 \vspace{.2cm} 
\begin{align}\label{Crucial}
D_x(\mathcal{P}_{t,s}\varphi)(x)v &= \E\left[D_x[\varphi(\phi(t,s,x))]v\right]= \E\left[ D_x\varphi(\phi(t,s,x))J_{t,s}(\ccdot,x)v\right] \notag \\[.2cm]
    &= \E\left[\langle D_x\varphi(\phi(t,s,x))\D X_{t}^{s,x}(\om), \eta^v\rangle_{\mathcal{H}}\right] = \E\left[\langle\D\varphi(\phi(t,s,x)), \eta^v\rangle_{\mathcal{H}}\right]\notag\\[.2cm]
    &= \E\left[\langle\D\varphi(\phi(t,s,x)),\eta^v\rangle_{\mathcal{H}}\right] = \E\left[\varphi(\phi(t,s,x))\delta(\eta^v)\right] \notag\\[.1cm]   
    & = \E\left[\varphi(\phi(t,s,x))\int_s^t\sigma(\varsigma, \phi(\varsigma,s))^TJ^{T}_{t,\varsigma}M_{t,s}^{-1}J_{t,s}v\star dW_\varsigma\right] \notag\\[.2cm]
    &= \E\left[ \varphi(\phi(t,s,x))\int_s^t\eta_{\varsigma,s}^v\star dW_\varsigma\right],
 \end{align}
where in the first and second lines, we applied chain rule for mean square gradient and Malliavin derivative respectively, and third line is Malliavin integration by parts formula (e.g.,~\cite{Nulart}) and, the stochastic integral in the fourth or fifth line is interpreted in the sense of Skorokhod, i.e., \\[.4cm]$\displaystyle \int_s^t\eta^v_{\varsigma,s}\star dW_{\varsigma}$ is the divergence of the process $\{\eta^v_{\varsigma, s}\I_{[s, t]}(\varsigma): \varsigma\geqslant s\}$ (see equation (\ref{Div})).
  
 Next, since  $\mathcal{C}_\infty^1(\Rd)$ is dense in $\mathcal{C}_\infty(\Rd),$ we have $(\varphi_n)_{n\in \N_1}\subset \mathcal{C}_\infty^1(\Rd)$, $\varphi_n\underset{n\rightarrow\infty}{\longrightarrow} \varphi\in \mathcal{C}_\infty(\Rd)$  so that 
 \begin{align}\label{Density1}
 \begin{cases} \displaystyle \lim_{n\rightarrow\infty}\mathcal{P}_{s,t}\varphi_n(x) = \mathcal{P}_{s,t}\varphi(x), \\[.4cm]
 \displaystyle \lim_{n\rightarrow\infty} D_x(\mathcal{P}_{s,t}\varphi_n)(x)v =  \E\left[ \varphi(\phi(t,s,x))\int_s^t\eta_{\varsigma,s}^v\star dW_\varsigma\right].
 \end{cases}
 \end{align}
 
 \noindent On the other hand, by Proposition \ref{SMOOTH}, there exists a function $0<p_{s,t}\in \mathcal{C}_\infty^\infty(\Rd\times\Rd)$ such that, formally, $\p\big(\{\om: \phi(t,s,\om,x)\in dy\}\big)= P(s, x; t, dy)=p_{s,t}(x, y)dy$. This implies that 
 \begin{align}\label{Density2}
\notag  \lim_{n\rightarrow\infty}D_x(\mathcal{P}_{s,t}\varphi_n)(x)v &= \lim_{n\rightarrow\infty}\int_{\Rd}\varphi_n(y)D_xp_{s,t}(x,y)vdy\\[.2cm]
 &= \int_{\Rd}\varphi(y)D_xp_{s,t}(x,y)v dy= D_x(\mathcal{P}_{s,t}\varphi)(x)v.
 \end{align}
  Comparing (\ref{Density1}) and (\ref{Density2}), we have that (\ref{Crucial}) holds for all $\varphi\in \mathcal{C}_\infty(\Rd)$. 
 
 \bigskip
 Next, the Cauchy--Schwartz inequality yields 
 \begin{align}\label{CWin}
 | D_x(\mathcal{P}_{s,t}\varphi)(x)v| &\leqslant \sqrt{(\mathcal{P}_{s,t}\varphi^2)(x)}\left(\E\left\vert \int_s^t\eta^v_{\varsigma,s}\star dW_\varsigma\right\vert^2\right)^{1/2}, \quad \varphi\in \mathcal{C}_\infty(\Rd).
 \end{align}
 By generalised It\^o isometry (cf.~\cite{Nulart}), we have 
 \begin{align*}
 \E\left\vert \int_s^t\eta^v_{\varsigma,s}\star dW_\varsigma\right\vert^2&= \E\left(\int_s^t\vert \eta^v_{\varsigma,s}\vert^2 d\varsigma\right)+\E\left(\int_s^t\int_s^t\langle\D_\xi\eta^v_{\xi,s},\D_{\varsigma}\eta^v_{\varsigma,s}\rangle_{\Rd\times\Rd} d\xi d\varsigma \right)\\[.4cm]
 &\hspace{0cm} \leqslant \E\left(\int_s^t\vert \eta^v_{\varsigma,s}\vert^2 d\varsigma\right)+ \E\left(\int_s^t\int_s^t\|\D_\xi\eta^v_{\varsigma,s}\|^2_{\Rd\times\Rd} d\xi d\varsigma \right)\\[.4cm]
 &\hspace{0cm} = \E\|\eta^v\|^2_{\mathcal{H}}+\sum_{k=1}^m\E\left(\int_s^t \|\D_\xi^k\eta^v\|^2_{\mathcal{H}}d\xi\right).
 \end{align*}
 
 \medskip
 \noindent Then, by the inequality (\ref{SKEst})  and (\ref{CWin}), there exists $C_{T}>0$ such that 
 \begin{align}\label{ESTF}
 |D_x(\mathcal{P}_{s,t}\varphi)(x)v|\leqslant C_{T}\|\varphi\|_{\infty}|v| \quad \forall \;x,v\in \Rd, \; \varphi\in \mathcal{C}_\infty(\Rd).
 \end{align}
 
 \medskip
 \noindent Finally, let $z_\ell = \ell x+(1-\ell)y, \;\; \ell\in [0, 1]$ and set $v = x-y.$ Then, by the mean value theorem and inequality (\ref{ESTF}), we have 
 \begin{align*}
 \left\vert\mathcal{P}_{s,t}\varphi(x)-\mathcal{P}_{s,t}\varphi(y)\right\vert\leqslant \int_0^1| D_x(\mathcal{P}_{s,t}\varphi)(z_\ell)v|d\ell \leqslant C_{T}\|\varphi\|_{\infty}| x-y|.
 \end{align*}\qed
 
\newpage
{\footnotesize \setstretch{.95}
\bibliographystyle{plainnat}

}

\end{document}